\documentclass[english,twoside,reqno]{amsart}
\usepackage[square,comma]{natbib}
\usepackage[english]{babel} 
\usepackage{exscale}   
\usepackage{amsmath}
\usepackage{amsfonts}
\usepackage{amssymb}
\usepackage{amsxtra}
\usepackage{pstricks}
\usepackage{stmaryrd}
\usepackage{xspace}
\usepackage{epsfig}
\usepackage{mathrsfs}
\numberwithin{equation}{section}
\newtheorem{Theorem}{Theorem}[section]
\newtheorem{Lemma}[Theorem]{Lemma}
\newtheorem{Corollary}[Theorem]{Corollary}
\newtheorem{Proposition}[Theorem]{Proposition}
\newtheorem{Definition}[Theorem]{Definition}
\newtheorem{Remark}[Theorem]{Remark}
\newtheorem{Remarks}[Theorem]{Remarks}
\newtheorem{Example}[Theorem]{Example}
\newtheorem{Notation}[Theorem]{Notation}
\newtheorem{Question}[Theorem]{Question}

\newcommand{\Bset}{\mathbb{B}}
\newcommand{\Cset}{\mathbb{C}}

\newcommand{\Nset}{\mathbb{N}}

\newcommand{\Rset}{\mathbb{R}}
\newcommand{\Sset}{\mathbb{S}}

\newcommand{\Zset}{\mathbb{Z}} 
\newcommand{\cA}{\ensuremath{{\mathcal A}}\xspace}         % Caligr-A
\newcommand{\cB}{\ensuremath{{\mathcal B}}\xspace}         % Caligr-B
\newcommand{\cC}{\ensuremath{{\mathcal C}}\xspace}         % Caligr-C
\newcommand{\cF}{\ensuremath{{\mathcal F}}\xspace}         % Caligr-F
\newcommand{\cM}{\ensuremath{{\mathcal M}}\xspace}         % Caligr-M
\newcommand{\cN}{\ensuremath{{\mathcal N}}\xspace}         % Caligr-N
\newcommand{\cO}{\ensuremath{{\mathcal O}}\xspace}         % Caligr-O
\newcommand{\cR}{\ensuremath{{\mathcal R}}\xspace}         % Caligr-R
\newcommand{\ctM}{\ensuremath{{\widetilde{\mathcal M}}}\xspace} % Caligr-~M
\newcommand{\scrB}{\ensuremath{{\mathscr{B}}}\xspace} %
\newcommand{\scrI}{\ensuremath{{\mathscr{I}}}\xspace} %
\newcommand{\scrtI}{\ensuremath{{\widetilde{\mathscr{I}}}}\xspace} %
\newcommand{\scrM}{\ensuremath{{\mathscr{M}}}\xspace} %
\newcommand{\scrN}{\ensuremath{{\mathscr{N}}}\xspace} %
\newcommand{\aaa}{\mathbf{a}}
\newcommand{\bbb}{\mathbf{b}}
\newcommand{\ii}{\mathbf{i}}
\newcommand{\jj}{\mathbf{j}}

\newcommand{\1}{\ensuremath{{\rm 1\kern-.25em l}}\xspace}  % identity
\newcommand{\alg}{{\operatorname{alg}}} 
\newcommand{\ball}{\operatorname{B}}                    % (Unit) Ball   
\newcommand{\id}{\operatorname{id}}                        % id 

\newcommand{\trace}{\operatorname{tr}}        % normalized trace
\newcommand{\tail}{\ensuremath{\mathrm{tail}}}                 
\newcommand{\Aut}[1]{\ensuremath{\operatorname{Aut}(#1)}\xspace}  % Aut
\newcommand{\Mor}[1]{\ensuremath{\operatorname{Mor}(#1)}\xspace}  % Mor
\newcommand{\set}[2]{\mathopen{\{}#1\mathop{|}#2\mathclose{\}}}%set
\newcommand{\distr}{\operatorname{distr}}   
\newcommand{\sotlim}{\ensuremath{\textsc{sot-}\lim}\xspace}   % strong-top
\newcommand{\wotlim}{\ensuremath{\textsc{wot-}\lim}\xspace}   % weak-top
\newcommand{\sot}{\ensuremath{\textsc{sot}}\xspace}           % strong-top
\newcommand{\wot}{\ensuremath{\textsc{wot}}\xspace}           % weak-top
\begin{document}
\title[Noncommutative extended De Finetti theorem]{A noncommutative extended De Finetti theorem}
\author[C. K\"ostler]{Claus K\"ostler}
\address{University of Illinois at Urbana-Champaign, Department of Mathematics,
Altgeld Hall, 1409 West Green Street, Urbana, 61801, USA}
\email{koestler@uiuc.edu} 
\subjclass[2000]{Primary 46L53; Secondary 60G09, 47A53}
\keywords{Noncommutative de Finetti theorem, distributional symmetries, exchangeability, 
spreadability, noncommutative conditional independence, mean ergodic theorem, noncommutative
Kolmogorov Zero-One Law, noncommutative Bernoulli shifts} 
\date{June 22, 2008}
\begin{abstract}
The extended de Finetti theorem characterizes exchangeable infinite random sequences 
as conditionally i.i.d.~and shows that the apparently weaker distributional symmetry
of spreadability is equivalent to exchangeability. Our main result is a noncommutative 
version of this theorem.  
  
In contrast to the classical result of Ryll-Nadzewski \cite{RyNa57a}, exchangeability turns 
out to be stronger than spreadability for infinite noncommutative random sequences. Out of our 
investigations emerges noncommutative conditional independence in terms of a von Neumann 
algebraic structure closely related to Popa's notion of commuting squares \cite{Popa83a} and K\"ummerer's 
generalized Bernoulli shifts \cite{Kuem88a}. Our main result is 
applicable to classical probability, quantum probability, in particular free probability  
\cite{VDN92a}, braid group representations and Jones subfactors \cite{GHJ89a,GoKo08a}.    
\end{abstract}
 
\maketitle

%%%%%%%%%%%%%%%%%%%%%%%%  S E C T I O N  %%%%%%%%%%%%%%%%%%%%%%%%%%%%%%
% \tableofcontents    %   (looks ugly)      Inhaltsverzeichnis
%%%%%%%%%%%%%%%%%%%%%%%%%%%%%%%%%%%%%%%%%%%%%%%%%%%%%%%%%%%%%%%%%%%%%%%
%%%%%%%%%   Handmade Table of Contents  follows                 %%%%%%%
%%%%%%%%%%%%%%%%%%%%%%%%%%%%%%%%%%%%%%%%%%%%%%%%%%%%%%%%%%%%%%%%%%%%%%%
{%
\def\widedotfill{\leaders\hbox to 10pt{\hfil.\hfil}\hfill}
\def\pg#1{\widedotfill\rlap{\hbox to 15pt{\hfill{\small#1}}}\par}
\rightskip=15pt\leftskip=10pt%
\newcommand{\sct}[2]{\noindent\llap{\hbox to%
    10pt{{#1}\hfill}}~\mbox{#2~}}
\newcommand{\separ}{\vspace{.2em}}

\hrule
%\vspace{12pt}
\section*{Contents}
\sct{}{Introduction and main result}%
\pg{\pageref{section:intro}}%
\separ
\sct{1}{Preliminaries and terminology}%
\pg{\pageref{section:preliminaries}}%
\separ
\sct{2}{Noncommutative stationary processes}%
\pg{\pageref{section:ncsp}}%
\separ
\sct{3}{Conditional independence and conditional factorizability}%
\pg{\pageref{section:cond-indy-fact}}%
\separ
\sct{4}{Stationarity and conditional independence/factorizability}%
\pg{\pageref{section:ci-cf-stationary}}%
\separ
\sct{5}{Noncommutative i.i.d.~sequences may be non-stationary}%
\pg{\pageref{section:non-stationary}}%
\separ
\sct{6}{Stationarity with strong mixing and noncommutative Bernoulli shifts}%
\pg{\pageref{section:kol-0-1}}%
\separ
\sct{7}{Spreadability implies conditional order independence}%
\pg{\pageref{section:spread2cond-indep}}%
\separ
\sct{8}{Spreadability implies conditional full independence}%
\pg{\pageref{section:spread2indep}}%
\separ
\sct{9}{Some applications and outlook}%
\pg{\pageref{section:applications}}%
\separ
\sct{}{References}%
\pg{\pageref{section:bibliography}}
}
\vspace{12pt}
\hrule
%%%%%%%%%%%%%%%%%%%%%%%%%%%%%%%%%%%%%%%%%%%%%%%%%%%%%%%%%%%%%%%%%%%%%%%
%%%%%%%%%   Handmade Table of Contents  end                     %%%%%%%
%%%%%%%%%%%%%%%%%%%%%%%%%%%%%%%%%%%%%%%%%%%%%%%%%%%%%%%%%%%%%%%%%%%%%%%

%%%%%%%%%%%%%%%%%%%%%%%%  S E C T I O N  %%%%%%%%%%%%%%%%%%%%%%%%%%%%%%
\section*{Introduction and main result}
\label{section:intro}
%%%%%%%%%%%%%%%%%%%%%%%%%%%%%%%%%%%%%%%%%%%%%%%%%%%%%%%%%%%%%%%%%%%%%%%
The characterization of random objects with distributional symmetries is
of major interest in modern probability theory and Kallenberg's recent 
monograph \cite{Kalle05a} provides an impressive account on this subject. 
Already in the early 1930s, de Finetti showed that infinite exchangeable 
random sequences are conditionally i.i.d.\ or, more intuitively formulated, 
mixtures of i.i.d.~random variables \cite{Fine31a,CiRe96a}. An early version 
of his celebrated characterization is that for every infinite sequence of 
exchangeable $\{0,1\}$-valued random variables $X \equiv (X_1, X_2, X_3, \ldots)$, 
there exists a probability measure $\nu$ on $[0,1]$ such that the law 
$\mathcal{L}(X)$ is given by 
$$
\mathcal{L}(X) = \int_{[0,1]} m(p) d\nu(p).
$$
Here $m(p)$ denotes the infinite product of the measure with Bernoulli 
distribution $(p, 1-p)$. An extension of this result from the set $\{0,1\}$ 
to any compact Hausdorff space $\Omega$ goes back to Hewitt and Savage 
\cite{HeSa55a} and soon after it was realized by Ryll-Nadzewski \cite{RyNa57a} 
that the apparently weaker distributional symmetry of spreadability is 
equivalent to exchangeability for infinite random sequences. A further 
extension to standard Borel spaces is provided by Aldous in his monograph 
on exchangeability \cite{Aldo85a}. Let us mention that `spreadability' is 
also known as `contractability' in the probability theory and shares common 
ground with `subsymmetric sequences' in Banach space theory. 

Our main result is a noncommutative version of the following extended 
de Finetti theorem. We have adapted its formulation in 
\cite[Theorem 1.1]{Kalle05a} to the purposes of this paper:
%%%%%%%%%%%%%%%%%%%%%%%%%%%%%%%%%%%%%%%%%%%%%%%%%%%%%%%%%%%%%%%%%%%%%%
\begin{Theorem}\label{thm:definetti-kallenberg}
%%%%%%%%%%%%%%%%%%%%%%%%%%%%%%%%%%%%%%%%%%%%%%%%%%%%%%%%%%%%%%%%%%%%%%
Let $ X \equiv (X_n)_{n \in \Nset}\colon (\Omega, \Sigma, \mu) \to 
(\Omega_0, \Sigma_0)$ be a sequence of random variables, where 
$(\Omega, \Sigma)$ and $(\Omega_0, \Sigma_0)$ are standard Borel spaces 
and $\mu$ is a probability measure. Then the following conditions are 
equivalent:
\begin{enumerate}
\item[(a)] $X$ is exchangeable;
\item[(b)] $X$ is spreadable;
\item[(c)] $X$ is conditionally i.i.d.
\end{enumerate}    
\end{Theorem}
Here the conditioning is with respect to the tail $\sigma$-field of the 
random sequence $X$. Three different proofs of this result can be found 
in \cite{Kalle05a} and it is worthwhile to point out that the two 
implications (a) $\Rightarrow$ (b) and (c) $\Rightarrow$ (a) are fairly 
clear; the main work rests on proving the implication (b) $\Rightarrow$ (c). 

An early noncommutative version of de Finetti's theorem was given by 
St{\o}rmer for exchangeable states on the infinite tensor product of 
C*-algebras \cite{Stoe69a}. His pioneering work stimulated further results 
in quantum statistical physics and quantum probability, with focus on bosonic 
systems \cite{HuMo76a,Huds81a,FLV88a}. A quite general noncommutative analogue 
of de Finetti's theorem is obtained by Accardi and Lu in a C*-algebraic setting, 
where only the tail algebra (generated by the exchangeable infinite 
noncommutative random sequence) is required to be commutative  \cite{AcLu93a}. 
Quite recently, inspired by Good's formula and Speicher's free cumulants 
\cite{Spei98a}, a combinatorial approach by Lehner unifies cumulant techniques 
in a *-algebraic setting of exchangeability systems 
\cite{Lehn04a,Lehn03a,Lehn05a,Lehn06a}. He shows that exchangeability entails 
properties of cumulants, as they are known in classical probability to be 
characterizing for (conditional) independence. Presently, no results on 
noncommutative versions of de Finetti's theorem seem to be available in the 
literature beyond the case of commutative tail algebras and Lehner's 
combinatorial results for exchangeability systems; and no results at all are 
present in the noncommutative realm for the extended de Finetti theorem, 
Theorem \ref{thm:definetti-kallenberg}.

Our framework towards a noncommutative version of the extended de Finetti 
theorem needs to be capable to efficiently handle tail events. This suggests 
to deal right from the beginning with W*-algebraic probability spaces. 
We will work with noncommutative probability spaces $(\cM,\psi)$ which 
consist of a von Neumann algebra $\cM$ (with separable predual) and a 
faithful normal state $\psi$ on $\cM$. A noncommutative random variable 
$\iota$ from $(\cA_0,\varphi_0)$ to $(\cM,\psi)$ is given by an injective 
*-homomorphism $\iota\colon \cA_0 \to \cM$ such that $\varphi_0= \psi \circ \iota$.
Furthermore, we require that the $\psi$-preserving conditional expectation from 
$\cM$ onto $\iota(\cA_0)$ exists (see Section \ref{section:preliminaries} for 
further details). 

Here we will constrain our investigations to a random sequence $\scrI$, given 
by an infinite sequence of \emph{identically distributed} random variables 
$\iota\equiv (\iota_n)_{n \in \Nset_0}$ from $(\cA_0,\varphi_0)$ to $(\cM,\psi)$. 
This simplification improves the transparency of our approach,  since it allows 
us to realize $\iota$ as \emph{injective} mappings from the \emph{single} 
probability space $(\cA_0,\varphi_0)$. A treatment beyond identically distributed 
random variables is possible and of course of interest; it would start with a 
probability space $(\cM,\psi)$ and a sequence of (not necessarily injective) 
normal *-homomorphisms from a von Neumann algebra $\cA$ into $\cM$. Since the 
distributional symmetries considered herein will lead anyway to stationarity 
(which implies identical distributions), we omit this primary technical 
generalization.     

We recall that $\cM$ is of the form $L^\infty(\Omega, \Sigma, \mu)$ for some standard 
Borel space $(\Omega, \Sigma, \mu)$ as soon as $\cM$ is commutative; and then 
one has $\psi = \int_{\Omega} \cdot\, d\mu$. In this case a random variable 
$X\colon (\Omega, \Sigma, \mu) \to (\Omega_0, \Sigma_0)$ reappears as an 
injective *-homomorphism $\iota\colon L^\infty(\Omega_0, \Sigma_0, \mu_X) \to  
L^\infty(\Omega, \Sigma,\mu)$ with $\iota (f):= f \circ X$ (the measure $\mu_X$ 
is the distribution of $X$). Given a sequence of random variables 
$(X_n)_{n \in \Nset_0}$, the constraint of identically distributed $X_n$'s ensures 
that we can identify all image measures $\mu_{X_n}$ with the single measure 
$\mu_{X_0}$. Note that this approach is free of any conditions on the existence 
of moments of the $X_n$'s.

Throughout we will work with a noncommutative notion of conditional 
independence which, by our main result, can actually seen to emerge 
out of the transfer of the extended de Finetti theorem to noncommutative 
probability. It encompasses Popa's notion of `commuting squares' 
in subfactor theory \cite{Popa83a,GHJ89a} as well as  Voiculescu's 
freeness with amalgamation \cite{VDN92a}, aside of tensor 
independence and many other examples coming from generalized 
Gaussian random variables \cite{BKS97a,GuMa02a}. 

Consider the random sequence $\scrI$ which generates the von Neumann 
subalgebras $$\cM_I := \bigvee_{i \in I} \iota_i(\cA_0)$$ for 
subsets $I \subset \Nset_0$ and the tail algebra 
$$
\cM^\tail:= \bigcap_{n \in \Nset_0}  \bigvee_{k \ge n}\iota_k(\cA_0).
$$
Let $E_{\cM^\tail}$ denote the $\psi$-preserving conditional 
expectation from $\cM$ onto $\cM^\tail$. Then we say that $\scrI$ is 
\emph{full $\cM^\tail$-independent} if
$$
E_{\cM^\tail}(xy) = E_{\cM^\tail}(x) E_{\cM^\tail}(y)
$$          
for all $x \in \cM^\tail \vee \cM_I$ and $y \in \cM^\tail \vee \cM_J$ 
with $I \cap J = \emptyset$. We will also meet a weaker notion of 
independence, called  \emph{order $\cM^\tail$-independence}, which 
requires the (ordered) sets $I$ and $J$ to satisfy $I <J$ or $I > J$, 
instead of disjointness. 

These two notions of conditional independence do \emph{not} require 
$\cM^\tail \subset \cM_I$ and allow a noncommutative dual formulation 
of random measures as they are necessary in the context of de Finetti's 
theorem. Interesting on its own, the paradigm of an infinite sequence 
$X$ of exchangeable $\{0,1\}$-valued random variables clearly 
illustrates that, in its algebraic reformulation, stipulating the 
inclusion $\cM^\tail \subset \cM_I$ implies the triviality 
$\cM^\tail \simeq \Cset$ and thus forces $X$ to be i.i.d. Thus it is 
crucial to allow $\cM^\tail \not \subset \cM_I$ if one is interested 
in transferring results on distributional symmetries to a noncommutative 
setting.

In order to state our main result, we informally introduce the relevant 
distributional symmetries. Given the two random sequences $\scrI$ and $\tilde{\scrI}$
with random variables $\iota$ resp.~ $\tilde{\iota}$, both from $(\cA_0,\varphi_0)$ to 
$(\cM,\psi)$, we write
$$
(\iota_0, \iota_1, \iota_2, \ldots ) 
\stackrel{\distr}{=} 
(\tilde{\iota}_0, \tilde{\iota}_1, \tilde{\iota}_2, \ldots)
$$
if $\scrI$ and $\tilde{\scrI}$ have the same distribution: 
$$
\psi\big(\iota_{\ii(1)}(a_1) \iota_{\ii(2)}(a_2) 
                    \cdots   \iota_{\ii(n)}(a_n)\big) 
= 
\psi\big(\tilde{\iota}_{\ii(1)}(a_1) \tilde{\iota}_{\ii(2)}(a_2) 
                           \cdots  \tilde{\iota}_{\ii(n)}(a_n)\big) 
$$
for all $n$-tuples $\ii\colon \{1, 2, \ldots, n\} \to \Nset_0$, 
$(a_1, \ldots, a_n) \in \cA_0^n$ and $n \in \Nset$. Now a random sequence 
$\scrI$ is said to be \emph{exchangeable} if its distribution is 
invariant under permutations: 
$$
 (\iota_0,\iota_1, \iota_2,  \ldots ) 
\stackrel{\distr}{=} 
(\iota_{\pi(0)},\iota_{\pi(1)}, \iota_{\pi(2)},  \ldots ) 
$$
for any finite permutation $\pi \in \Sset_\infty$ of $\Nset_0$.  We say 
that random sequence $\scrI$ is \emph{spreadable} if every subsequence 
has the same distribution:
$$
(\iota_0,\iota_1, \iota_2, \ldots )  
\stackrel{\distr}{=}  
(\iota_{n_0},\iota_{n_1}, \iota_{n_2}, \ldots ) 
$$ 
for any subsequence $(n_0, n_1, n_2,\ldots)$ of $(0,1,2,\ldots)$. Finally, 
$\scrI$ is \emph{stationary} if the distribution is 
shift-invariant:
$$
 (\iota_0,\iota_1, \iota_2,  \ldots )  
\stackrel{\distr}{=} 
(\iota_{k},\iota_{k+1}, \iota_{k+2},\ldots ) 
$$  
for all $k\in \Nset$.

We are ready to formulate our main result, a noncommutative dual version 
of Theorem \ref{thm:definetti-kallenberg}.
%%%%%%%%%%%%%%%%%%%%%%%%%%%%%%%%%%%%%%%%%%%%%%%%%%%%%%%%%%%%%%%%%%%%%%
\begin{Theorem}\label{thm:main}
%%%%%%%%%%%%%%%%%%%%%%%%%%%%%%%%%%%%%%%%%%%%%%%%%%%%%%%%%%%%%%%%%%%%%% 
Let $\scrI$ be a random sequence with (identically distributed) random 
variables
$$
\iota\equiv (\iota_i)_{i \in \Nset_0} 
\colon 
(\cA_0, \varphi_0) \to (\cM, \psi)
$$
and consider the following conditions:
\begin{enumerate}
\item[(a)]
$\scrI$ is exchangeable;
\item[(b)]
$\scrI$ is spreadable;
\item[(c)]
$\scrI$ is stationary and full $\cM^\tail$-independent; 
\item[(d)]
$\scrI$ is full $\cM^\tail$-independent;
\item[(c$_{\text{\normalfont{o}}}$)]
$\scrI$ is stationary and order $\cM^\tail$-independent; 
\item[(d$_{\text{\normalfont{o}}}$)]
$\scrI$ is order $\cM^\tail$-independent.
\end{enumerate}
Then we have the implications:
\[
\begin{matrix}
                  \text{\normalfont{(a)}} 
 & \Rightarrow  & \text{\normalfont{(b)}} 
 & \Rightarrow  & \text{\normalfont{(c)}} 
 & \Rightarrow  & \text{\normalfont{(d)}} \\
 &              &  
 &              & \Downarrow 
 &              & \Downarrow \\
 &              &  
 &              & \text{\normalfont{(c$_{\text{\normalfont{o}}}$)}} 
 & \Rightarrow & \text{\normalfont{(d$_{\text{\normalfont{o}}}$)}} 
\end{matrix}
\]
Moreover, there are examples such that \normalfont{(a)} 
$\not \Leftarrow$ \normalfont{(b)} $\not \Leftarrow$ 
\normalfont{(c)}  $\not \Leftarrow$ \normalfont{(d)}
and \normalfont{(c$_{\text{\normalfont{o}}}$)} 
$\not \Leftarrow$ \normalfont{(d$_{\text{\normalfont{o}}}$)}.
\end{Theorem}
%%%%%%%%%%%%%%%%%%%%%%%%%%%%%%%%%%%%%%%%%%%%%%%%%%%%%%%%%%%%%%%%%%%%%%
Similar to the classical case, the hard part of the proof is that
spreadability implies conditional \emph{full} independence. This is 
done by means from noncommutative ergodic theory.  

One might object that a noncommutative version of the extended de Finetti 
theorem should provide an equivalence of these conditions. But our 
investigations show that such a common folklore understanding would be 
conceptually misleading in the noncommutative world. The crucial 
implications from distributional symmetries to conditional (full/order) 
independence are still valid. All listed converse implications fail due 
to deep structural reasons, and the others are presently open in the 
generality of our setting.  

The failure of the implication `(b) $\Rightarrow$ (a)' relies on the fact that, 
in the noncommutative realm, spreadability of infinite random sequences goes 
beyond the representation theory of the symmetric group. As developed in
\cite{GoKo08a}, braid group representations with infinitely many generators 
lead to \emph{braidability}, a new symmetry intermediate to exchangeability 
and spreadability. This braidability extends exchangeability and provides a 
rich source of spreadable noncommutative random sequences such that the reverse 
implication `(b) $\Rightarrow$ (a)' fails. Some of these `counter-examples' are 
known in subfactor theory as vertex models from quantum statistical physics. 
The inequivalence of exchangeability and spreadability is a familiar phenomena 
for random arrays \cite{Kalle05a}. Since this phenomena already occurs for 
infinite sequences in the noncommutative setting, it provides another facet of 
the common folklore result that $(d+1)$-dimensional classical models correspond to 
$d$-dimensional quantum models \cite{EvKa98a}.  

Examples for the failure of the implication `(c) $\Rightarrow$ (b)' are also 
available in the context of braid group representations. It is shown in 
\cite{GoKo08a} that an appropriate cocycle perturbation of the unilateral 
shift of a stationary random sequence may obstruct spreadability without effecting 
the structure of conditional full independence.  Again, related `counter-examples' 
arise in the most natural manner. For example, the symbolic shift on the Artin generators 
of the braid group $\Bset_\infty$ induces an endomorphism of the braid group von 
Neumann algebra $L(\Bset_\infty)$ such that, when acting on the subalgebra $L(\Bset_2)$, 
the resulting stationary random sequence exhibits the failure of `(c) $\Rightarrow$ (b)'. 
Such phenomena are impossible in the classical case by Theorem \ref{thm:definetti-kallenberg}.    

Finally, one can not expect in the noncommutative realm that i.i.d.~random 
sequences are automatically stationary. The failure of the implication  
`(d) $\Rightarrow$ (c)', and thus of `(d$_{\text{\normalfont{o}}}$) $\Rightarrow$ 
(c$_{\text{\normalfont{o}}}$)', is closely related to the fact that our notion 
of noncommutative conditional independence is more general than (conditioned 
versions of) tensor independence or free independence. The latter two
notions of independence enjoy universality properties \cite{Spei97a,BeSc02a,NiSp06a} 
which immediately entail stationarity of an i.i.d.\ random sequence. In particular, they
are rigid with respect to certain `local perturbations' of noncommutative random sequences. 
But we will see that, starting with a stationary (conditionally full/order) independent 
random sequence, our more general notion of independence is non-rigid with respect to 
such `local perturbations'. Related examples arise again in the context of braid group 
representations or noncommutative Gaussian random variables. Thus stationarity plays are 
more distinguished role in the quantum setting and cannot simply be deduced from independence 
properties as it is the case for classical probability or Voiculescu's free probability. 

A closer look at Theorem \ref{thm:main} reveals that it is  `dual' to the 
usual formulations of de Finetti's theorem. In terms of quantum physics, our theorem 
is formulated in the Heisenberg picture, whereas the usual formulations  use the 
Schr{\"o}dinger picture. Or equivalently phrased: our result is on the level of the 
von Neumann algebra, whereas the latter identify the geometry of exchangeable states 
in the predual of the von Neumann algebra. Using the theory of noncommutative $L^1$-spaces 
it would be of interest to examine in detail the geometry of exchangeable, spreadable or 
`conditionally independent' subspaces.

We summarize the content of this paper.

Section \ref{section:preliminaries} introduces our setting of noncommutative probability
spaces, random sequences and distributional symmetries. It closes with the proof of
some of the elementary implications of Theorem \ref{thm:main}. 

Section \ref{section:ncsp} provides the needed background results on noncommutative stationary 
processes and their endomorphisms. Since spreadability immediately implies stationarity, most 
parts of the proof of Theorem \ref{thm:main} will be carried out in an equivalent framework 
of stationary processes.
  
We introduce in Section \ref{section:cond-indy-fact} two noncommutative versions
of classical conditional independence, called  `conditional independence' (CI) and 
`conditional factorizability' (CF). Both notions are equivalent if the conditioning 
is trivial or appropriate additional algebraic structure is supposed. But (CF) is a 
priori weaker than (CI) and more easily to control in applications. Their definition 
reflects that the conditioning is with respect to a von Neumann algebra which may 
\emph{not} be contained in the image of two random variables. Further we relate 
`conditional independence' to Popa's `commuting squares' of von Neumann algebras.      

The main result of Section \ref{section:ci-cf-stationary} is that (CI) and (CF) are 
equivalent for a stationary random sequence if the conditioning is with respect to 
a subalgebra of the fixed point algebra of the corresponding endomorphism, see Theorem 
\ref{thm:fact=indep}. Moreover, we introduce the notions of `conditional order 
independence' (CI$_{\text{o}}$) and `conditional order factorizability' (CF$_{\text{o}}$). 
These two notions are apparently weaker and reflect that the index set $\Nset_0$ of the 
random sequence is considered as an ordered set. Already (CF$_{\text{o}}$), the weakest 
of the four properties, will suffice to establish mixing properties of stationary 
processes. Finally, we illustrate (CI) and (CF) by the algebraic reformulation of de Finetti's 
original example, an infinite sequence of exchangeable $\{0,1\}$-valued random variables. 

Section \ref{section:non-stationary} focuses on appropriate `local perturbations' of $\Cset$-independent
stationary random sequences. Our main result is that a noncommutative i.i.d.~random sequence may be
non-stationary. We provide related examples and disprove the implications `(d) $\Rightarrow$ (c)' 
and `(d$_{\text{\normalfont{o}}}$) $\Rightarrow$ (c$_{\text{\normalfont{o}}}$)' of Theorem \ref{thm:main}.

Section \ref{section:kol-0-1} provides a noncommutative generalization of Kolmogorov's zero-one-law 
for a random sequence with (CF$_{\text{o}}$). Further we show in Theorem \ref{thm:mixing} that 
(CF$_{\text{o}}$) and stationarity imply strong mixing over the tail algebra and fixed point characterization
results. We coin in this section also the notion of a noncommutative Bernoulli shift, as it is suggested 
by our results on distributional symmetries and inspired by K\"ummerer's notion of a generalized Bernoulli shift.
These shifts can be recognized as the unilateral `discrete' version of noncommutative continuous Bernoulli 
shifts from \cite{HKK04a}. 

Section \ref{section:spread2cond-indep} is devoted to an integral part of the noncommutative extended de Finetti 
theorem, the proof that spreadability of a random sequence yields conditional order independence (CI$_\text{o}$). 
Here the conditioning is shown to be with respect to the tail algebra of the random sequence. 

Section \ref{section:spread2indep} upgrades the results of the previous section. Our main result is Theorem
\ref{thm:spread2ind} which provides the proof of the crucial part of Theorem \ref{thm:main}: spreadability
implies conditional full independence (CI) of a random sequence. An important tool within its proof is a local
version of the mean ergodic theorem, Theorem \ref{thm:local-mean-ergodic}. It will allow us to perform mean ergodic
approximations in a spreadability preserving manner.     

Applications and an outview are contained in Section \ref{section:applications}. We cite results from 
\cite{GoKo08a} on braidability and on the failure of the implications `(a) $\Leftarrow$ (b)' and `(b) $\Leftarrow$ (c)'
of the noncommutative extended de Finetti theorem, Theorem \ref{thm:main}. Moreover, we present a
general central limit theorem for spreadable random sequences which can be regarded to be the noncommutative
prototype of a conditioned central limit theorem. We also discuss briefly its potential connections
to interacting Fock spaces. Finally, we give immediate applications of Theorem \ref{thm:main} to inequalities in noncommutative $L^1$-spaces, as they appear in the work of Junge and Xu.

\subsection*{Acknowledgments} 
The present paper took its origin from work with Rolf Gohm on one of the most 
simple examples coming from the Jones fundamental construction \cite{GoKo08a}, 
and joint work with Roland Speicher on the structure of noncommutative white noises 
\cite{KoSp07a}. At both occasions we found `spreadability' without being aware of it. 
We are indebted to Marius Junge and Wojciech Jaworski who independently pointed out 
possible connections to distributional symmetries and initiated the author's investigations 
resulting in the present paper. We are thankful to several helpful discussions with 
Benoit Collins, Rolf Gohm, Marius Junge, James Mingo and Roland Speicher in the course 
of writing this paper.  

%%%%%%%%%%%%%%%%%%%%%%%%  S E C T I O N  %%%%%%%%%%%%%%%%%%%%%%%%%%%%%%
\section{Preliminaries and Terminology}
\label{section:preliminaries}
%%%%%%%%%%%%%%%%%%%%%%%%%%%%%%%%%%%%%%%%%%%%%%%%%%%%%%%%%%%%%%%%%%%%%%%
Noncommutative notions of probability spaces have in common that they
consist of an algebra which is equipped with a linear functional. Here we 
shall work with the W*-algebraic version of such spaces, since they allow
us to capture probabilistic tail events of random sequences. We refer the 
reader to \cite{AFL82a,Kuem85a,KuMa98a,VDN92a,NiSp06a} for further 
information on noncommutative probability spaces, in particular *-algebraic 
or C*-algebraic settings.       
%%%%%%%%%%%%%%%%%%%%%%%%%%%%%%%%%%%%%%%%%%%%%%%%%%%%%%%%%%%%%%%%%%%%%
\begin{Definition}\normalfont \label{def:ncps}
%%%%%%%%%%%%%%%%%%%%%%%%%%%%%%%%%%%%%%%%%%%%%%%%%%%%%%%%%%%%%%%%%%%%%
A \emph{probability space} $(\cM,\psi)$ consists of  a von Neumann 
subalgebra $\cM$ with separable predual and a faithful normal state
$\psi$ on $\cM$. A von Neumann algebra $\cM_0$ of $\cM$ is said to be 
\emph{$\psi$-conditioned} if the $\psi$-preserving conditional expectation 
$E_{\cM_0}$ from $\cM$ onto $\cM_0$ exists. %The resulting triple 
%$(\cM,\psi, \cM_0)$ will also be called a \emph{conditioned probability space}.
Two probability spaces  $(\cM_1,\psi_1)$ and $(\cM_2, \psi_2)$ are said to be
\emph{isomorphic} if there exists an isomorphism $\Pi \colon \cM_1 \to \cM_2$ such
that $\psi_1 = \psi_2 \circ \Pi$. The $\psi$-preserving automorphisms of $\cM$
will be denoted by $\Aut{\cM,\psi}$. 
\end{Definition}
%%%%%%%%%%%%%%%%%%%%%%%%%%%%%%%%%%%%%%%%%%%%%%%%%%%%%%%%%%%%%%%%%%%%%
By Takesaki's theorem, the 
$\psi$-preserving conditional expectation $E_{\cM_0}$ exists
if and only if $\sigma_t^\psi(\cM_0) = \cM_0$ for all $t\in \Rset$
\cite[IX, Theorem 4.2]{Take03b}. Here $\sigma_t^\psi$ denotes the modular 
automorphism group associated to $(\cM,\psi)$. Thus the existence of such a 
conditional expectation is automatic if $\psi$ is a trace, 
i.e.~ $\psi(xy) = \psi(yx)$ for all $x,y \in \cM$. 

The noncommutative generalization of random variables is casted 
in terms of *-homomorphisms \cite{AFL82a}. For the purpose
of this paper the following definition of a random variable will be 
sufficient.

\begin{Definition}\normalfont\label{def:rv}
Let $(\cA_0, \varphi_0)$ and $(\cM,\psi)$ be two probability spaces. 
A \emph{ random variable} is an injective *-homomorphism 
$\iota\colon \cA_0 \to \cM$ satisfying two additional properties:
\begin{enumerate}
\item $\varphi_0 = \psi \circ \iota$; 
\item $\iota(\cA_0)$ is $\psi$-conditioned.
\end{enumerate}    
A random variable will also be addressed as the mapping 
$\iota\colon  
(\cA_0, \varphi_0) \to (\cM,\psi).
$
\end{Definition}

Every classical random variable in the context of standard measure spaces yields 
by algebraisation a random variable in the sense of Definition \ref{def:rv}. 
Conversely, if the von Neumann algebra $\cM$ is commutative then the usual notion of a 
random variable on standard probability spaces can be recovered from Definition \ref{def:rv}.
Note that our assumption of injectivity is no restriction if a single random variable
is considered. 

\begin{Remark}\normalfont\label{rem:conditioned}
Assertion (ii) in above definition is superfluous if $\psi$ is a trace. Note also that this
assertion has equivalent formulations:
\begin{enumerate}
\item[(iii)] 
$\iota$ intertwines the modular automorphism groups of $(\cA_0,\varphi_0)$ and $(\cM,\psi)$;
\item[(iv)]
There exists a (unique) unital completely positive map $\iota^+\colon \cM \to \cA_0$ 
satisfying $\psi(x \iota(a))= \varphi_0(\iota^+(x)a)$ for all $x \in \cM$ and $a\in \cA_0$. 
\end{enumerate}
The map $\iota^+$ is also called the adjoint of $\iota$. We refer the reader to 
\cite{AcCe82a,GrKu82a,HKK04a,AnDe06a} for further details and background results 
on the equivalences of (ii) to (iv).
\end{Remark} \normalfont
\begin{Remark} \normalfont
Commonly (selfadjoint) operators in the von Neumann algebra $\cM$ (or, more generally, 
its noncommutative $L^p$-spaces) are also denoted as `noncommutative random variables'
in the literature. Such approaches require assumptions on the existence of higher moments 
of a random variable. The framework of injective *-homomorphisms has the advantage that 
it is free of any assumptions on the existence of moments. Of course, 
we can easily produce a random variable in the operator sense from our setting by 
considering $\iota(x)$ for some fixed $x \in \cA_0$.   
\end{Remark}
We are interested in sequences of random variables. 
%%%%%%%%%%%%%%%%%%%%%%%%%%%%%%%%%%%%%%%%%%%%%%%%%%%%%%%%%%%%%%%%%%%%%%
\begin{Notation}\normalfont
%%%%%%%%%%%%%%%%%%%%%%%%%%%%%%%%%%%%%%%%%%%%%%%%%%%%%%%%%%%%%%%%%%%%%%
We write $I < J$ for two subsets $I, J \subset \Nset_0$ if $i<j$ for all 
$i \in I $ and $j \in J$. The cardinality of $I$ is denoted by $|I|$. For 
$N \in \Nset$, we denote by $I + N$ the shifted set $\set{i + N}{i \in I}$.
\end{Notation}
%%%%%%%%%%%%%%%%%%%%%%%%%%%%%%%%%%%%%%%%%%%%%%%%%%%%%%%%%%%%%%%%%%%%%%
\begin{Definition}\normalfont \label{def:rvs-neu}
%%%%%%%%%%%%%%%%%%%%%%%%%%%%%%%%%%%%%%%%%%%%%%%%%%%%%%%%%%%%%%%%%%%%%%
An  \emph{(identically distributed) random sequence $\mathscr{I}$} is 
a sequence of random variables 
$$
\iota \equiv (\iota_i)_{i\in \Nset_0} \colon (\cA_0,\varphi_0) \to (\cM,\psi).
$$
The family $(\cA_I)_{I \subset \Nset_0}$, with von Neumann subalgebras
\[
\cA_I = \bigvee_{i \in I} \iota_i(\cA_0),
\] 
is called the \emph{canonical filtration (generated by $\mathscr{I}$)} 
and $\mathscr{I}$ is said to be \emph{minimal} if $\cA_{\Nset_0} = \cM$. 
The von Neumann subalgebra
\[
\cA^{\tail} := \bigcap_{n \in \Nset_0}\bigvee_{k \ge n} \iota_k(\cA_0)
\]
is called the \emph{tail algebra of $\scrI$}.    

Suppose a second random sequence $\scrtI$ is defined by the random variables
$
\widetilde{\iota}\equiv (\widetilde{\iota}_i)_{i \in \Nset_0} \colon 
(\cA_0, \varphi_0) \to (\widetilde{\cM}, \widetilde{\psi}).
$ 
Then $\scrI$ and $\scrtI$ are \emph{isomorphic} 
if there exists an isomorphism $\Pi\colon \widetilde{\cM} \to \cM$ such that 
$
\psi \circ \Pi = \widetilde{\psi}
$
and
$
\Pi \circ 
\widetilde{\iota}_n =\iota_n
$
for all $n \in \Nset_0$. 
\end{Definition}
%%%%%%%%%%%%%%%%%%%%%%%%%%%%%%%%%%%%%%%%%%%%%%%%%%%%%%%%%%%%%%%%%%%%%%
Whenever it is convenient, we may turn a random sequence into a 
minimal one by restricting the probability space $(\cM,\psi)$ to 
$(\cA_{\Nset_0}, \psi|_{\cA_{\Nset_0}})$.
We have already introduced distributional symmetries in the introduction. 
Here we present equivalent definitions which are less intuitive, but more 
convenient within our proofs.  
%%%%%%%%%%%%%%%%%%%%%%%%%%%%%%%%%%%%%%%%%%%%%%%%%%%%%%%%%%%%%%%%%%%%%%
\begin{Notation}\normalfont
%%%%%%%%%%%%%%%%%%%%%%%%%%%%%%%%%%%%%%%%%%%%%%%%%%%%%%%%%%%%%%%%%%%%%%
The group $\Sset_\infty$ is the inductive limit of the symmetric groups 
$\Sset_n$, $n \ge 2$, where $\Sset_n$ is generated on $\Nset_0$ by the 
transpositions 
$$
\pi_{i}\colon (i-1,i) \to (i,i-1)
$$ 
with $1 \le i <n$. By $[n]$ we denote the linearly ordered set $\{1,2,\ldots,n\}$. 
\end{Notation}
%%%%%%%%%%%%%%%%%%%%%%%%%%%%%%%%%%%%%%%%%%%%%%%%%%%%%%%%%%%%%%%%%%%%%%
\begin{Definition}\normalfont \label{def:relations}
%%%%%%%%%%%%%%%%%%%%%%%%%%%%%%%%%%%%%%%%%%%%%%%%%%%%%%%%%%%%%%%%%%%%%%
Let $\ii,\jj\colon [n] \to \Nset_0$ be two $n$-tuples. 
\begin{enumerate}
\item
$\ii$ and $\jj$ are \emph{translation equivalent}, in symbols:
$\ii \sim_\theta \jj$, if there exists $k \in \Nset_0$ such that 
$$
\ii = \theta^k \circ \jj \qquad \text{or}  \qquad  \theta^k \circ \ii = \jj .
$$
Here denotes $\theta$ the right translation $m \mapsto m+1$ on $\Nset_0$. 
\item
$\ii$ and $\jj$ are \emph{order equivalent}, in symbols: $\ii \sim_o \jj$,
if there exists a permutation $\pi \in \Sset_\infty$ such that
$$
\ii = \pi \circ \jj  \qquad  \text{and} \qquad \pi|_{\jj([n])} \text{ is order preserving.} 
$$ 
\item
$\ii$ and $\jj$ 
are \emph{symmetric equivalent}, in symbols: $\ii \sim_\pi \jj$, if there exists a 
permutation $\pi \in \Sset_\infty$ such that 
$$
\ii  =  \pi \circ \jj.  
$$  
\end{enumerate}
\end{Definition}
%%%%%%%%%%%%%%%%%%%%%%%%%%%%%%%%%%%%%%%%%%%%%%%%%%%%%%%%%%%%%%%%%%%%%%
We have the implications 
$ (\ii \sim_\theta \jj) \Rightarrow (\ii \sim_o \jj) \Rightarrow (\ii \sim_\pi \jj)$.  
%%%%%%%%%%%%%%%%%%%%%%%%%%%%%%%%%%%%%%%%%%%%%%%%%%%%%%%%%%%%%%%%%%%%%
\begin{Remark}\normalfont \label{rem:partial-shift}
%%%%%%%%%%%%%%%%%%%%%%%%%%%%%%%%%%%%%%%%%%%%%%%%%%%%%%%%%%%%%%%%%%%%%
Order equivalence of two $n$-tuples $\ii$ and $\jj$ can also be 
expressed with the help of the partial shifts $(\theta_N)_{N\ge 0} 
\colon \Nset_0 \to \Nset_0$, where 
\[
\theta_N (n) = 
\begin{cases}
n & \text{if $n < N$};\\
n+1 & \text{ if $n \ge N$}.
\end{cases}
\]
Each $\theta_N$ is order-preserving and it is easy to see that 
$\ii \sim_o \jj$ if and only if there exist partial shifts 
$\theta_{N_1}, \theta_{N_2},\ldots,\theta_{N_k}$ such that 
$\theta_{N_1}\circ \theta_{N_2} \circ\cdots \circ \theta_{N_k} 
\circ \ii = \jj$. Note also that any subsequence $(n_0, n_1, 
n_2, \cdots)$ of the infinite sequence $(0,1,2,3, \ldots)$ 
can be approximated via actions of the subshifts 
$(\theta_{N})_{N \ge 0}$.
\end{Remark}
%%%%%%%%%%%%%%%%%%%%%%%%%%%%%%%%%%%%%%%%%%%%%%%%%%%%%%%%%%%%%%%%%%%%%
\begin{Remark}\normalfont
%%%%%%%%%%%%%%%%%%%%%%%%%%%%%%%%%%%%%%%%%%%%%%%%%%%%%%%%%%%%%%%%%%%%%
Order equivalence is used in the context of a general limit theorem 
in \cite{KoSp07a} and our present formulation is an equivalent one. 
\end{Remark}
%%%%%%%%%%%%%%%%%%%%%%%%%%%%%%%%%%%%%%%%%%%%%%%%%%%%%%%%%%%%%%%%%%%%%
For the notation of mixed higher moments of random variables,
it is convenient to use Speicher's notation of multilinear maps. 
%%%%%%%%%%%%%%%%%%%%%%%%%%%%%%%%%%%%%%%%%%%%%%%%%%%%%%%%%%%%%%%%%%%%%
\begin{Notation} \normalfont \label{notation:mlm}
%%%%%%%%%%%%%%%%%%%%%%%%%%%%%%%%%%%%%%%%%%%%%%%%%%%%%%%%%%%%%%%%%%%%%
Let the random sequence $\scrI$ be given by 
$$
\iota\equiv (\iota_i)_{i \in \Nset_0} \colon (\cA_0, \varphi_0) \to (\cM, \psi).
$$
We put, for $\ii\colon [n] \to \Nset_0$, 
$\mathbf{a} = (a_1, \ldots, a_{n}) \in \cA_0^n$ and $n \in \Nset$, 
\begin{eqnarray*}
\iota[\ii; \mathbf{a}] 
&:=& \iota_{\ii(1)}(a_1) \iota_{\ii(2)}(a_2) \cdots \iota_{\ii(n)}(a_{n}),  \label{eqn:iota1}\\
\psi_\iota[\ii; \mathbf{a}] 
&:=& \psi\big(\iota[\ii; \mathbf{a}]  \big)
\label{eqn:iota2}. 
\end{eqnarray*}
\end{Notation}
%%%%%%%%%%%%%%%%%%%%%%%%%%%%%%%%%%%%%%%%%%%%%%%%%%%%%%%%%%%%%%%%%%%%%%
We are ready to introduce distributional symmetries in  terms of the 
mixed moments of a random sequence. 
%%%%%%%%%%%%%%%%%%%%%%%%%%%%%%%%%%%%%%%%%%%%%%%%%%%%%%%%%%%%%%%%%%%%%%
\begin{Definition}\label{def:ds}
%%%%%%%%%%%%%%%%%%%%%%%%%%%%%%%%%%%%%%%%%%%%%%%%%%%%%%%%%%%%%%%%%%%%%%
A random sequence $\scrI$ is 
\begin{enumerate}
\item
\emph{exchangeable} if, for any $n \in \Nset$,  
$$\psi_\iota[\ii; \cdot\,] = \psi_\iota[\jj; \cdot\,] 
\quad\text{whenever} \quad \ii \sim_\pi \jj;$$ 
\item
\emph{spreadable} if, for any $n \in \Nset$, 
$$
\psi_\iota[\ii; \cdot\,] = \psi_\iota[\jj; \cdot\,] 
\quad\text{whenever} \quad \ii \sim_o \jj;$$
\item
\emph{stationary} if, for any $n \in \Nset$,  
$$
\psi_\iota[\ii; \cdot\,] = \psi_\iota[\jj; \cdot\,]
\quad\text{whenever} \quad \ii \sim_\theta \jj.$$       
\end{enumerate}
\end{Definition}
We close this section with the proof of the obvious implications in
the noncommutative extended de Finetti theorem.
%%%%%%%%%%%%%%%%%%%%%%%%%%%%%%%%%%%%%%%%%%%%%%%%%%%%%%%%%%%%%%%%%%%%%%%%%%
\begin{proof}[Proof of Theorem \ref{thm:main}, elementary parts]
It is evident from Definition \ref{def:relations} and Definition 
\ref{def:ds} that exchangeability implies spreadability, and that
spreadability implies stationarity. This shows the implication 
`(a) $\Rightarrow$ (b)' and the elementary parts on stationarity 
of the implications `(b) $\Rightarrow$ (c)' and  
`(b) $\Rightarrow$ (c$_{\text{0}}$)'. The implications 
`(c) $\Rightarrow$ (d)' and `(c$_{\text{0}}$) $\Rightarrow$ (d$_{\text{0}}$)'
are trivial. 
\end{proof}  
%%%%%%%%%%%%%%%%%%%%%%%%%%%%%%%%%%%%%%%%%%%%%%%%%%%%%%%%%%%%%%%%%%%%%%%%%%

%%%%%%%%%%%%%%%%%%%%%%%%%%%%%%%%%%%%%%%%%%%%%%%%%%%%%%%%%%%%%%%%%%%%%%%%%%%
%%%%%%%%%%%%%%%%%%%%%%%  S E C T I O N  %%%%%%%%%%%%%%%%%%%%%%%%%%%%%%%%%%%
%%%%%%%%%%%%%%%%%%%%%%%%%%%%%%%%%%%%%%%%%%%%%%%%%%%%%%%%%%%%%%%%%%%%%%%%%%%
\section{Noncommutative stationary processes}
\label{section:ncsp}
%%%%%%%%%%%%%%%%%%%%%%%%%%%%%%%%%%%%%%%%%%%%%%%%%%%%%%%%%%%%%%%%%%%%%%%%%%%
Exchangeable or spreadable random sequences are stationary and can thus be 
expressed as stationary processes. Since the remaining sections of this paper
will rest on this well known connection, we provide more in detail some of
their specific properties in this section. We will introduce stationary processes 
such that they are in a canonical correspondence to stationary random sequences 
(see Definition \ref{def:ds}). Their notion is very closely related to K\"ummerer's 
approach in \cite{Kuem93UP, Kuem03a} (see also \cite[Section 2.1]{Gohm04a}).
Moreover, we present a result from \cite{Kuem88b} which shows that a unilateral
stationary process (as introduced next) extends to a bilateral stationary process. 
%%%%%%%%%%%%%%%%%%%%%%%%%%%%%%%%%%%%%%%%%%%%%%%%%%%%%%%%%%%%%%%%%%%%%%
\begin{Definition}\normalfont \label{def:stat-process}
%%%%%%%%%%%%%%%%%%%%%%%%%%%%%%%%%%%%%%%%%%%%%%%%%%%%%%%%%%%%%%%%%%%%%%
A \emph{(unilateral) stationary process} $\scrM$ consists of a
probability space $(\cM,\psi)$, a $\psi$-conditioned subalgebra 
$\cM_0 \subset \cM$ and an endomorphism $\alpha$ of $\cM$ satisfying
\begin{enumerate}
\item unitality: $\alpha(\1) =\1$;
\item stationarity: $\psi \circ \alpha = \psi$;
\item conditioning: 
$\alpha$ and the modular automorphism group $\sigma^\psi_t$ commute for all $t\in \Rset$. 
\end{enumerate}
The stationary process $\scrM$ is also denoted by the quadruple $(\cM,\psi,\alpha, \cM_0)$ and 
\[
\iota^\alpha \equiv (\iota_i^\alpha)_{i \in \Nset_0}\colon (\cM_0, \psi_0) \to (\cM,\psi),  
\qquad 
\iota_i^\alpha := \alpha^i|_{\cM_0}, 
\]
is called the \emph{random sequence associated to $\scrM$}, for brevity also denoted by 
$\scrI^\alpha$. 

The family of von Neumann subalgebras $(\cM_I)_{I \subset \Nset_0}$, with
$$
\cM_I:= \bigvee_{i \in I} \alpha^i(\cM_0),
$$
 is called the \emph{canonical filtration 
(generated by $\scrM$)}. The von Neumann subalgebra
\[
\cM^\tail = \bigcap_{n \in \Nset_0}  \alpha^n(\cM)  
\]
is called the \emph{tail algebra of $\scrM$}. We denote by $\cM^\alpha$ the \emph{fixed point
algebra} of the endomorphism $\alpha$.  

Finally, two minimal stationary processes $\scrM$ and $\widetilde{\scrM}$ are 
\emph{isomorphic} if there exists an isomorphism $\Pi \colon \widetilde{\cM} 
\to \cM$ such that
\[
\psi \circ \Pi 
= \widetilde{\psi}, \quad\Pi \circ  \widetilde{\alpha}  
={\alpha} \circ \Pi , \quad \Pi(\widetilde{\cM}_0) = \cM_0. 
\]

\end{Definition}

The modular condition (iii) is needed for a non-tracial state $\psi$ to ensure:    
\begin{itemize}
\item[$\diamond$]
Von Neumann algebras generated by the $\alpha^n(\cM_0)$'s and their intersections are 
$\psi$-conditioned (see Remark \ref{rem:def-stat-process}).
\item[$\diamond$]
Stationary processes and stationary random sequences are in correspondence 
(see Lemma \ref{lem:corres-rv-endo}).     
\item[$\diamond$]
A unilateral stationary process extends to a bilateral stationary process 
(see Theorem \ref{thm:unilateral-bilateral}). 
\end{itemize}

\begin{Remark}\normalfont \label{rem:def-stat-process}
Condition (iii) of Definition \ref{def:stat-process} entails that the 
$\psi$-preserving conditional expectation $E_{\cM_I}$ from $\cM$ onto 
$\cM_I$ exists: $\cM_0$ is globally $\sigma^\psi_t$-invariant and now 
condition (iii) implies that $\alpha(\cM_0)$ and, more generally, $\cM_I$ 
are globally $\sigma^\psi_t$-invariant. Thus Takesaki's theorem on the 
existence of $\psi$-preserving conditional expectations applies. Of course, 
the condition (iii) can be dropped if $\psi$ is a trace. We are indebted to 
K\"ummerer for simple examples on hyperfinite $III_\lambda$ factors such that 
$\alpha(\cA_0)$ fails to be globally $\sigma_t^\psi$-invariant without 
condition (iii) \cite{KuemPC}. 
\end{Remark}

To avoid reiterations throughout this paper we shall use the following convention 
for properties of a stationary process.

\begin{Definition}\normalfont
The stationary process $\scrM$ is said to have property `A' if its
associated random sequence $\scrI^\alpha$ has property `A'. 
For example, $\scrM$ is minimal if its associated random sequence $\scrI^\alpha$ is minimal. 
\end{Definition}

The canonical filtrations of a stationary process $\scrM$ and its associated random sequence
$\scrI^\alpha$ always coincide. But the tail algebra $\scrM^\tail$ of $\scrM$ may be larger 
than the tail algebra of $\scrI^\alpha$,
\[
\scrM^{\scrI\tail}
= \bigcap_{n \in \Nset_0} \bigvee_{k \ge n} \iota^{(\alpha)}_k(\cM_0) 
= \bigcap_{n \in \Nset_0} \bigvee_{k \ge n} \alpha^k(\cM_0). 
\]
%%%%%%%%%%%%%%%%%%%%%%%%%%%%%%%%%%%%%%%%%%%%%%%%%%%%%%%%%%%%%%%%%%%%
\begin{Lemma}
%%%%%%%%%%%%%%%%%%%%%%%%%%%%%%%%%%%%%%%%%%%%%%%%%%%%%%%%%%%%%%%%%%%%
If $\scrM$ is minimal, then $\cM^{\scrI\tail} = \cM^\tail$.
\end{Lemma}
%%%%%%%%%%%%%%%%%%%%%%%%%%%%%%%%%%%%%%%%%%%%%%%%%%%%%%%%%%%%%%%%%%%%
\begin{proof}
This is easily concluded from 
\[
\bigvee_{k \ge n} \iota_k^{(\alpha)}(\cM_0) 
= \bigvee_{k \ge n} \alpha^k(\cM_0)
= \alpha^n \bigvee_{k \ge 0}\alpha^k(\cM_0) \subseteq  \alpha^n (\cM)
\]
and minimality.    
\end{proof}
%%%%%%%%%%%%%%%%%%%%%%%%%%%%%%%%%%%%%%%%%%%%%%%%%%%%%%%%%%%%%%%%%%%%% 

We continue with the correspondence between stationary processes and
stationary random sequences under the condition of minimality. We 
include this well known result for reasons of transparency since
the proof of the noncommutative version of the extended de Finetti  
theorem makes heavily use of it.      
%%%%%%%%%%%%%%%%%%%%%%%%%%%%%%%%%%%%%%%%%%%%%%%%%%%%%%%%%%%%%%%%%%%%% 
\begin{Lemma}\label{lem:corres-rv-endo}
%%%%%%%%%%%%%%%%%%%%%%%%%%%%%%%%%%%%%%%%%%%%%%%%%%%%%%%%%%%%%%%%%%%%%
There is a one-to-one correspondence between (equivalence classes of)
\begin{enumerate}
\item[(a)]
minimal stationary processes $\scrM = (\cM, \psi, \alpha, \cM_0)$;
\item[(b)]
minimal stationary random sequences $\scrI$ with random variables
\[
(\iota_n)_{n \ge 0}\colon (\cA_0, \varphi_0) \to (\cM,\psi) .
\]
 
\end{enumerate}
The correspondence from (a) to (b) is given by
\[
(\cA_0,\varphi_0) := (\cM_0, \psi|_{\cM_0}) 
\quad\text{and}\quad
\iota_n := \alpha^n|_{\cM_0}.  
\]
The correspondence from (b) to (a) is established via
\[
\cM_0 := \iota_0(\cA_0) 
\quad\text{and}\quad 
\alpha(\iota[\ii; \mathbf{a}]):= \iota[\Theta \circ \ii; \mathbf{a}]
\] 
for all $n \in \Nset$, $n$-tuples $\ii \colon [n] \to \Nset_0$ and 
$\mathbf{a} \in \cA_0^n$.
\end{Lemma}
%%%%%%%%%%%%%%%%%%%%%%%%%%%%%%%%%%%%%%%%%%%%%%%%%%%%%%%%%%%%%%%%%%%%%
\begin{proof}
%%%%%%%%%%%%%%%%%%%%%%%%%%%%%%%%%%%%%%%%%%%%%%%%%%%%%%%%%%%%%%%%%%%%%
We omit all fairly clear parts of the proof and only show that 
the properties of $\scrI$ imply the modular condition 
$\alpha \sigma_t^\psi = \sigma_t^\psi \alpha$. Since the von Neumann 
algebras $\iota_n(\cA_0)$ are $\psi$-conditioned, the random variables 
$\iota_n$ intertwine $\sigma_t^{\varphi_0}$ and $\sigma_t^\psi$, the 
modular automorphism groups of $(\cA_0,\varphi_0)$ and $(\cM,\psi)$ 
(see Remark \ref{rem:conditioned} and \cite[Lemma 2.5]{AnDe06a}). Thus
\begin{eqnarray*}
\sigma_t^\psi \circ \alpha(\iota[\ii; \mathbf{a}])
&=& \sigma_t^\psi \iota[\Theta \circ \ii;\mathbf{a}]
= \iota[\Theta \circ \ii;\sigma_t^{\varphi_0}(\mathbf{a})]
=  \alpha(\iota[\ii; \sigma_t^{\varphi_0}(\mathbf{a})])\\
&=&  \alpha \circ \sigma_t^\psi (\iota[\ii;\mathbf{a}]) 
\end{eqnarray*}
establishes $\alpha \sigma_t^\psi = \sigma_t^\psi \alpha$ on a weak*-total
subset of $\cM$. Here $\sigma_t^{\varphi_0}(\aaa)$ denotes the $n$-tuple 
$\big(\sigma_t^{\varphi_0}(a_1), \ldots \sigma_t^{\varphi_0}(a_n)\big)$.
\end{proof}
%%%%%%%%%%%%%%%%%%%%%%%%%%%%%%%%%%%%%%%%%%%%%%%%%%%%%%%%%%%%%%%%%%%%%
We will need the next theorem for approximations in the proof of Theorem 
\ref{thm:fact=indep}. 
%%%%%%%%%%%%%%%%%%%%%%%%%%%%%%%%%%%%%%%%%%%%%%%%%%%%%%%%%%%%%%%%%%%%%
\begin{Definition}\normalfont
%%%%%%%%%%%%%%%%%%%%%%%%%%%%%%%%%%%%%%%%%%%%%%%%%%%%%%%%%%%%%%%%%%%%%
A stationary process $\hat{\scrM}=(\hat{\cM},\hat{\psi}, \hat{\alpha}, 
\hat{\cM_0})$ is said to be \emph{bilateral} if the endomorphism 
$\hat{\alpha}$ is an automorphism of $\hat{\cM}$. A bilateral stationary 
process $\hat{\scrM}$ is \emph{minimal} if $\hat{\cM}= 
\bigvee_{n \in \Zset}\hat{\alpha}^n(\hat{\cM_0})$.
\end{Definition}
%%%%%%%%%%%%%%%%%%%%%%%%%%%%%%%%%%%%%%%%%%%%%%%%%%%%%%%%%%%%%%%%%%%%%%
%%%%%%%%%%%%%%%%%%%%%%%%%%%%%%%%%%%%%%%%%%%%%%%%%%%%%%%%%%%%%%%%%%%%%%
\begin{Theorem} \label{thm:unilateral-bilateral}
%%%%%%%%%%%%%%%%%%%%%%%%%%%%%%%%%%%%%%%%%%%%%%%%%%%%%%%%%%%%%%%%%%%%%%
A unilateral stationary process $\scrM =(\cM,\psi, \alpha,\cM_0)$ extends 
to a bilateral stationary process  $\hat{\scrM}=(\hat{\cM},\hat{\psi}, 
\hat{\alpha},\hat{\cM_0})$. In other words, there exists a random variable
$
j \colon (\cM,\psi) \to (\hat{\cM},\hat{\psi}) 
$
such that 
$$
j(\cM_0)= \hat{\cM}_0 
\quad\text{and} \quad 
j \alpha^n  = \hat{\alpha}^n j \qquad (n \in \Nset_0).
$$
If $\hat{\scrM}$ is minimal, then $\hat{\cM}^{\hat{\alpha}}= j(\cM^\alpha)$.  
\end{Theorem}
%%%%%%%%%%%%%%%%%%%%%%%%%%%%%%%%%%%%%%%%%%%%%%%%%%%%%%%%%%%%%%%%%%%%%
This theorem is immediate from K\"ummerer's work on state-preserving Markov 
dilations. We provide some results from \cite{Kuem88b} which are essential 
for its proof. 

Let $(\cA, \varphi)$ and $(\cB, \psi)$ be two probability spaces. A  
morphism $T \colon (\cA,\varphi) \to (\cB,\psi)$ is a unital completely
positive map $T\colon \cA \to \cB$ satisfying $\varphi = \psi \circ T$.  
The morphisms from $(\cA,\varphi)$ into itself are denoted by $\Mor{\cA,\varphi}$.  
%%%%%%%%%%%%%%%%%%%%%%%%%%%%%%%%%%%%%%%%%%%%%%%%%%%%%%%%%%%%%%%%%%%%%
\begin{Definition}\normalfont
%%%%%%%%%%%%%%%%%%%%%%%%%%%%%%%%%%%%%%%%%%%%%%%%%%%%%%%%%%%%%%%%%%%%%
A morphism $T\in \Mor{\cA,\varphi}$ admits a \emph{state-preserving 
dilation} if there exists a probability space $(\hat{\cA}, \hat{\varphi})$, 
an automorphism $\hat{T} \in \Aut{\hat{\cA},\hat{\varphi}}$, two morphisms 
$j \colon (\cA,\varphi) \to (\hat{\cA},\hat{\varphi})$ and $Q \colon 
(\hat{\cA},\hat{\varphi}) \to (\cA,\varphi)$ such that $ T^n = Q \,\hat{T}^n j$ 
for all $n \in \Nset_0$. A state-preserving dilation is \emph{minimal} if 
$\hat{\cA} = \bigvee_{n \in \Zset}\hat{T}^n j(\cA)$.    
\end{Definition}
%%%%%%%%%%%%%%%%%%%%%%%%%%%%%%%%%%%%%%%%%%%%%%%%%%%%%%%%%%%%%%%%%%%%%%
Note in above definition that  $T^n = Q \,\hat{T}^n j$ reads as 
$\id_{\cA}=Q j$ for $n=0$. This implies that $j$ is a random variable 
from $(\cA,\varphi)$ to $(\cB,\psi)$ and the composition $j Q$ is the 
$\psi$-preserving conditional expectation from $\cB$ onto $j(\cA)$. We 
refer the reader to \cite{Kuem85a} for further details.    
%%%%%%%%%%%%%%%%%%%%%%%%%%%%%%%%%%%%%%%%%%%%%%%%%%%%%%%%%%%%%%%%%%%%%%
\begin{Proposition}[\cite{Kuem88b}]\label{prop:dil-condition}
%%%%%%%%%%%%%%%%%%%%%%%%%%%%%%%%%%%%%%%%%%%%%%%%%%%%%%%%%%%%%%%%%%%%%%
Let $(\cA,\varphi)$ be a probability space and suppose $\alpha$ is a 
$\varphi$-preserving unital endomorphism of $\cA$. Then the following 
are equivalent:
\begin{enumerate}
\item[(a)] $\alpha$ admits a state-preserving dilation. 
\item[(b)] $\alpha$ commutes with the modular automorphism 
           group $\sigma_t^\varphi$.
\end{enumerate}
\end{Proposition}
%%%%%%%%%%%%%%%%%%%%%%%%%%%%%%%%%%%%%%%%%%%%%%%%%%%%%%%%%%%%%%%%%%%%%%
We include the proof from \cite{Kuem88b} for the convenience of the reader. 
It uses inductive limits of C*-algebras (see for example 
\cite[Subsection 1.23]{Saka71a}).  
%%%%%%%%%%%%%%%%%%%%%%%%%%%%%%%%%%%%%%%%%%%%%%%%%%%%%%%%%%%%%%%%%%%%%%
\begin{proof}
%%%%%%%%%%%%%%%%%%%%%%%%%%%%%%%%%%%%%%%%%%%%%%%%%%%%%%%%%%%%%%%%%%%%%%
The implication `(a) $\Rightarrow$ (b)' is shown in \cite[2.1.8]{Kuem85a}. 
So it remains to prove the converse.

For $n \in \Nset_0$ put $\cA_n := \cA$ and for $n \ge m$ interpret $\alpha^{n-m}$ 
as a *-isomorphism of $\cA_m$ into $\cA_n$. Define $\tilde{\cA}$ as the C*-inductive
limit of $\set{\cA_m; \alpha^{n-m}}{(n,m)\in \Nset_0 \times \Nset_0, n \ge m}$. 
Moreover, putting $\varphi_n := \varphi$ for $n \in \Nset_0$, one has 
$\varphi_n(\alpha^{n-m}(x)) = \varphi_n(x)$ for $x\in \cA_m = \cA = \cA_n$ with 
$n \ge m$. The state $\tilde{\varphi}$ on $\tilde{\cA}$ is introduced as the inductive 
limit of  $\set{\varphi_m; \alpha^{n-m}}{(n,m) \in \Nset_0\times \Nset_0, n \ge m}$ 
\cite[1.23.10]{Saka71a}.
Identify $\tilde{\cA}$ with the norm closure of $\bigcup_{n \ge 0} \cA_n$. Thus 
$\alpha(\cA_n)$ is identified with $\cA_{n-1}$ ($n\ge 1$). Consequently $\alpha$ 
extends to  $\bigcup_{n \ge 0} \cA_n$ and then to its norm closure $\tilde{\cA}$. 
Doing so we obtain a $\tilde{\varphi}$-preserving automorphism $\tilde{\alpha}$ of 
$\tilde{\cA}$. We define an injection $\tilde{j}$ by identifying $\cA$ with 
$\cA_0$.   

Since $\alpha$ commutes with the modular automorphism group $\sigma_t^\varphi$ the 
subalgebra $\alpha^n(\cA) \subset \cA$ is globally $\sigma_t^\varphi$-invariant and 
the $\varphi$-preserving conditional expectation from $\cA$ onto $\cA_n$ exists for 
all $n \in \Nset$ (see \cite[Theorem]{Take03b}). Correspondingly, for each 
$n \in \Nset$, we find a completely positive map $Q_n \colon \cA_n \to \cA$ such that, 
for $m \le n$,  
\begin{eqnarray*}
\varphi_n = \varphi \circ Q_n, \qquad Q_n \circ \tilde{j} = \id_\cA, \qquad Q_n|_{\cA_m} = Q_m.      
\end{eqnarray*}
By continuity this leads to a completely positive map $\tilde{Q}\colon \tilde{\cA} \to \cA$ such 
that, for $n \ge 1$,
\begin{eqnarray*}
\tilde{\varphi} = \varphi \circ \tilde{Q}, \qquad \tilde{Q} \circ \tilde{j} = \id_\cA, 
\qquad \tilde{Q}|_{\cA_n} = Q_n.   
\end{eqnarray*}
Let $\sigma_t^{\varphi_m}$ be the modular automorphism group associated to $(\cA_m,\varphi_m)$. 
It follows $\sigma_t^{\varphi_n}(x) = \sigma_t^{\varphi_m}(x)$ for $x \in \cA_m$, $n \ge m$. 
Therefore the modular groups on $\cA_n$ extend to a group $\sigma_t$ on $\tilde{\cA}$ such that
$\tilde{\varphi}$ satisfies the KMS condition with respect to $\sigma_t$ (see \cite[8.12.3]{Pede79a}). 
Hence $\tilde{\varphi}$ extends to a faithful normal state $\hat{\varphi}$ on 
$\hat{\cA}:= \Pi_{\tilde{\varphi}}(\tilde{\cA})^{\prime\prime}$ \cite[8.14.4]{Pede79a}.    

Now it is routine to show that $\tilde{\alpha}$ extends to the 
$\hat{\varphi}$-preserving automorphism $\hat{\alpha}$ of $\hat{\cA}$, the map  
$\tilde{Q}$ to the completely positive map $Q \colon \hat{\cA} \to \cA$ satisfying 
$\varphi \circ Q = \hat{\varphi}$ and the injection $\tilde{j}$ to an injective 
*-homomorphism $j \colon \cA \to \hat{\cA}$ such that $\varphi = \hat{\varphi} \circ j$
and $j(\cA) = \Pi_{\tilde{\varphi}}(\cA_0)^{\prime\prime}$. Finally, 
$\alpha^n =  Q \hat{\alpha}^n  j$ is immediately verified for $n \in \Nset_0$.           
\end{proof}
%%%%%%%%%%%%%%%%%%%%%%%%%%%%%%%%%%%%%%%%%%%%%%%%%%%%%%%%%%%%%%%%%%%%%
%%%%%%%%%%%%%%%%%%%%%%%%%%%%%%%%%%%%%%%%%%%%%%%%%%%%%%%%%%%%%%%%%%%%%
\begin{proof}[Proof of Theorem \ref{thm:unilateral-bilateral}]
%%%%%%%%%%%%%%%%%%%%%%%%%%%%%%%%%%%%%%%%%%%%%%%%%%%%%%%%%%%%%%%%%%%%%
The endomorphism $\alpha$ of $\scrM$ satisfies the condition (b) of Proposition 
\ref{prop:dil-condition}. Thus there exists a probability space $(\hat{\cM},\hat{\psi})$,
a $\hat{\psi}$-preserving automorphism $\hat{\alpha}$ of $\hat{\cM}$, and a random variable
$j\colon (\cM,\psi) \to (\hat{\cM},\hat{\psi})$ such that $j \alpha^n = \hat{\alpha^n} j$ 
for all $n \in \Nset_0$. Clearly $\hat{\cM_0}:= j(\cM_0)$ is a $\hat{\psi}$-conditioned subalgebra
of $\hat{\cM}$. Finally, $\hat{\cM}^{\hat{\alpha}}= j(\cM^\alpha)$ is the content of
\cite[Corollary 3.1.4]{Kuem85a}.
\end{proof}
%%%%%%%%%%%%%%%%%%%%%%%%%%%%%%%%%%%%%%%%%%%%%%%%%%%%%%%%%%%%%%%%%%%%%

%%%%%%%%%%%%%%%%%%%%%%%%%%%%%%%%%%%%%%%%%%%%%%%%%%%%%%%%%%%%%%%%%%%%%
%%%%%%%%%%%%%%%%%%%%%%%%%%%%%%%%%%%%%%%%%%%%%%%%%%%%%%%%%%%%%%%%%%%%%
\section{Conditional independence and conditional factorizability}
\label{section:cond-indy-fact}
%%%%%%%%%%%%%%%%%%%%%%%%%%%%%%%%%%%%%%%%%%%%%%%%%%%%%%%%%%%%%%%%%%%%%
From our investigations of distributional symmetries emerge two closely
related noncommutative generalizations of classical conditional independence. 
Here we concentrate on the case of two random variables; the more 
general setting of random sequences is covered in the consecutive 
section where we will meet a further ramification of these two notions. 
%%%%%%%%%%%%%%%%%%%%%%%%%%%%%%%%%%%%%%%%%%%%%%%%%%%%%%%%%%%%%%%%%%%%%%
\begin{Definition}\normalfont \label{def:independence}
%%%%%%%%%%%%%%%%%%%%%%%%%%%%%%%%%%%%%%%%%%%%%%%%%%%%%%%%%%%%%%%%%%%%%%
Let $(\cM,\psi)$ be a probability space with three $\psi$-conditioned 
von Neumann subalgebras $\cM_0, \cM_1$ and $\cM_2$. Then $\cM_1$ and 
$\cM_2$ are said to be
\begin{enumerate}
\item
\emph{$\cM_0$-independent} or \emph{conditionally independent} if 
\[
E_{\cM_0}(xy) = E_{\cM_0}(x) E_{\cM_0}(y)
\]
for all $x \in \cM_1 \vee \cM_0$ and $y\in \cM_2 \vee \cM_0$; 
\item
\emph{$\cM_0$-factorizable} or \emph{conditionally factorizable} if 
\[
E_{\cM_0}(xy) = E_{\cM_0}(x) E_{\cM_0}(y)
\]
for all $x \in \cM_1$ and $y\in \cM_2$. 
\end{enumerate}
\end{Definition}
%%%%%%%%%%%%%%%%%%%%%%%%%%%%%%%%%%%%%%%%%%%%%%%%%%%%%%%%%%%%%%%%%%%%%
This definition does \emph{not} assume the inclusion $\cM_0 \subset 
\cM_1 \cap \cM_2$. It is open if conditional factorizability implies 
conditional independence and thus the equivalence of these two notions. 
But this is of course the case if $\cM_0 \simeq \Cset$, and we will 
state in Lemma \ref{lem:fact+cond} further conditions under which 
$\cM_0$-factorizability implies $\cM_0$-independence. 
%%%%%%%%%%%%%%%%%%%%%%%%%%%%%%%%%%%%%%%%%%%%%%%%%%%%%%%%%%%%%%%%%%%%%
\begin{Remark}\normalfont
%%%%%%%%%%%%%%%%%%%%%%%%%%%%%%%%%%%%%%%%%%%%%%%%%%%%%%%%%%%%%%%%%%%%%
An alternative formulation of $\cM_0$-independence is sometimes easier 
to verify in applications. Under the assertions of Definition 
\ref{def:independence}, the following are equivalent: 
\begin{enumerate}
\item[(a)] $\cM_1$ and $\cM_2$ are $\cM_0$-independent;
\item[(b)] there exist $\cM_0$-independent von Neumann subalgebras 
$\ctM_1$ and $\ctM_2$ of $\cM$ such that $\cM_0 \subset \ctM_1 \cap 
\ctM_2$ and  $\cM_i \subset \ctM_i$ ($i=1,2$). 
\end{enumerate}
Since this equivalence is fairly clear, we omit its proof.  
\end{Remark}  
%%%%%%%%%%%%%%%%%%%%%%%%%%%%%%%%%%%%%%%%%%%%%%%%%%%%%%%%%%%%%%%%%%%%%
\begin{Remark}\normalfont \label{rem:kuem}
%%%%%%%%%%%%%%%%%%%%%%%%%%%%%%%%%%%%%%%%%%%%%%%%%%%%%%%%%%%%%%%%%%%%%
If $\cM_0 \simeq \Cset$, we will also write $\Cset$-independence 
instead of $\cM_0$-independence. Note that $\cM_1$ and $\cM_2$ are 
$\Cset$-independent if and only if  
$
\psi(xy)= \psi(x) \psi(y)
$
for all $x \in \cM_1$ and $y \in \cM_2$ \cite{Kuem88a}.  
\end{Remark}\normalfont
%%%%%%%%%%%%%%%%%%%%%%%%%%%%%%%%%%%%%%%%%%%%%%%%%%%%%%%%%%%%%%%%%%%%%
The failure of the inclusion $\cM_0 \subset \cM_1 \cap \cM_2$ happens 
frequently in the context of distributional symmetries and is, in 
classical probability, intimately related to random probability 
measures. We  illustrate this by the most simple example which may be 
taken from classical probability (just choose $\cA \simeq \Cset^2 
\otimes \Cset^2$ in Example \ref{exa:fact-indep}). 
%%%%%%%%%%%%%%%%%%%%%%%%%%%%%%%%%%%%%%%%%%%%%%%%%%%%%%%%%%%%%%%%%%%%
\begin{Example}\normalfont \label{exa:fact-indep}
Let $\cA_1$ and $\cA_2$ be two $\Cset$-independent von Neumann 
subalgebras of the probability space $(\cA,\varphi)$. We define the 
probability space $(\cM, \psi)$ by $\cM := \cA\oplus \cA$ and 
$\psi := \frac12(\varphi \oplus \varphi)$.  For $i =1,2$, the 
embeddings $ \cA_i \ni x \to x \oplus x \in \cM$ define the 
von Neumann subalgebras $\cM_1$ and $\cM_2$, respectively. 
Furthermore, we put $\cM_0 = \Cset \1_\cA \oplus \Cset \1_\cA 
\simeq \Cset^2$. One has $\cM_i \vee \cM_0 = \cA_i \oplus \cA_i$ for 
$i=1,2$ and calculates
$$
E_{\cM_0}(xy) = E_{\cM_0}(x) E_{\cM_0}(y)
$$   
for all $x \in \cM_1 \vee \cM_0$ and $y \in \cM_2 \vee \cM_0$. Thus 
$\cM_1$ and $\cM_2$ are $\cM_0$-independent. But  $\cM_1 \cap \cM_2 
\simeq \Cset$, so $\cM_0 \not\subset \cM_1 \cap \cM_2$.  
\end{Example}
%%%%%%%%%%%%%%%%%%%%%%%%%%%%%%%%%%%%%%%%%%%%%%%%%%%%%%%%%%%%%%%%%%%%
\begin{Remark}\normalfont
%%%%%%%%%%%%%%%%%%%%%%%%%%%%%%%%%%%%%%%%%%%%%%%%%%%%%%%%%%%%%%%%%%%%
Another calculation shows in the above example that $\cM_1$ and 
$\cM_2$ are $\Cset$-independent. But this is rather an accident 
because we have chosen identical states on each component of the 
direct sum. 
\end{Remark}
%%%%%%%%%%%%%%%%%%%%%%%%%%%%%%%%%%%%%%%%%%%%%%%%%%%%%%%%%%%%%%%%%%%% 
Evidently, $\cM_0$-independence implies $\cM_0$-factorizability. But 
it is open if $\cM_0$-factorizability implies $\cM_0$-independence. 
Frequently this can be concluded if additional algebraic structures 
are available (see also Theorem \ref{thm:fact=indep}). All presently 
known examples (within our setting) satisfy at least one of the 
following conditions.   
%%%%%%%%%%%%%%%%%%%%%%%%%%%%%%%%%%%%%%%%%%%%%%%%%%%%%%%%%%%%%%%%%%%%
\begin{Lemma} \label{lem:fact+cond}
%%%%%%%%%%%%%%%%%%%%%%%%%%%%%%%%%%%%%%%%%%%%%%%%%%%%%%%%%%%%%%%%%%%%
$\cM_0$-factorizability and $\cM_0$-independence are equivalent 
under each of the following additional assertions:
\begin{enumerate}
\item
(trivial conditioning) $\cM_0 \simeq \Cset$;
\item
(central conditioning) $\cM_0 \subset \cM \cap \cM^\prime$;
\item               
(classical probability) $\cM = \cM^\prime$;
\item
(relative commutants) $\cM_0 \subset {\cM_1}^\prime \cap {\cM_2}^\prime$;
\item
(commuting squares) $\cM_0 \subset \cM_1 \cap \cM_2$.
\end{enumerate}
\end{Lemma}  
%%%%%%%%%%%%%%%%%%%%%%%%%%%%%%%%%%%%%%%%%%%%%%%%%%%%%%%%%%%%%%%%%%%%
\begin{proof} Each of the assertions (i) to (iv) implies that the 
vector spaces 
$$
\set{a x}{a \in \cM_0, x \in \cM_1} \quad \text{and} \quad 
\set{yb}{b \in \cM_0, y \in \cM_2}
$$ 
are weak* total in $\cM_0 \vee \cM_1$ and $\cM_0 \vee \cM_2$, 
respectively. Thus the module property of conditional expectations and 
$\cM_0$-factorizability imply
$$
E_{\cM_0}(axyb) = a E_{\cM_0}(xy)b = a E_{\cM_0}(x) E_{\cM_0}(y)b = 
E_{\cM_0}(ax) E_{\cM_0}(yb). 
$$
This equalities extend bilinearly and an approximation argument 
completes the proof in the cases (i) to (iv). The proof under the 
assertion (v) is trivial.   
\end{proof}
%%%%%%%%%%%%%%%%%%%%%%%%%%%%%%%%%%%%%%%%%%%%%%%%%%%%%%%%%%%%%%%%%%%%
Our notion of conditional independence is in close contact with 
Popa's notion of \emph{commuting squares} \cite{Popa83b,Popa83a,
PiPo86a}. Detailed information on their role in subfactor theory is 
provided in \cite{JoSu97a,GHJ89a}. We will make frequent use of some 
of their properties. Note that these assertions do \emph{not} apply 
for conditional factorizability. 
%%%%%%%%%%%%%%%%%%%%%%%%%%%%%%%%%%%%%%%%%%%%%%%%%%%%%%%%%%%%%%%%%%%%
\begin{Proposition}
%%%%%%%%%%%%%%%%%%%%%%%%%%%%%%%%%%%%%%%%%%%%%%%%%%%%%%%%%%%%%%%%%%%%
Suppose $\cM_0 \subset \cM_1 \cap \cM_2$, in addition to the 
assertions of Definition \ref{def:independence}. Then the following 
are equivalent:
\begin{enumerate}
\item 
$\cM_1$ and $\cM_2$ are $\cM_0$-independent;
\item 
$E_{\cM_1}(\cM_2)= \cM_0$;
\item 
$E_{\cM_1} E_{\cM_2} = E_{\cM_0}$
\item 
$E_{\cM_1}E_{\cM_2} = E_{\cM_2} E_{\cM_1}$ 
and $\cM_1 \cap \cM_2 = \cM_0$.
\end{enumerate}
In particular, it holds $\cM_0 = \cM_1 \cap \cM_2$ if one and thus all 
of the four assertions are satisfied. 
\end{Proposition}
%%%%%%%%%%%%%%%%%%%%%%%%%%%%%%%%%%%%%%%%%%%%%%%%%%%%%%%%%%%%%%%%%%%%
\begin{proof}
%%%%%%%%%%%%%%%%%%%%%%%%%%%%%%%%%%%%%%%%%%%%%%%%%%%%%%%%%%%%%%%%%%%%
The tracial case for $\psi$ is proved in \cite[Prop. 4.2.1]{GHJ89a}. 
The non-tracial case follows from this, after some minor modifications 
of the arguments therein.  
\end{proof}
%%%%%%%%%%%%%%%%%%%%%%%%%%%%%%%%%%%%%%%%%%%%%%%%%%%%%%%%%%%%%%%%%%%%
We close this section with some remarks on examples and references 
which are closely related to conditional independence in our 
noncommutative setting. The author is presently not aware of 
published examples in the quantum setting beyond the assertions 
stated in Lemma \ref{lem:fact+cond}. It would be of interest to find 
examples of von Neumann algebras which are conditionally factorizable,
but not conditionally independent, if possible at all.   
%%%%%%%%%%%%%%%%%%%%%%%%%%%%%%%%%%%%%%%%%%%%%%%%%%%%%%%%%%%%%%%%%%%%
\begin{Remarks}\normalfont \label{rem:independence}
%%%%%%%%%%%%%%%%%%%%%%%%%%%%%%%%%%%%%%%%%%%%%%%%%%%%%%%%%%%%%%%%%%%%
(1) $\Cset$-independence emerged from investigations of K\"ummerer 
on the structure of stationary quantum Markov processes \cite{Kuem85a,
Kuem88b,Kuem88a,Kuem93UP}. Its generalization to commuting squares is 
explored further from the perspective of noncommutative probability in 
\cite{Rupp95a,Koes00a,Koes03a,HKK04a,KoSp07a}.\\
(2) Examples for $\Cset$-independence are classical independence, tensor 
independence and free independence. Further examples originate from 
pioneering work of Bo\.zejko and Speicher \cite{Bosp91a,BoSp94a} and 
are given by generalized or noncommutative Gaussian random variables 
\cite{BKS97a,GuMa02a,Krol02a}. The most well known among them are 
$q$-Gaussian random variables. Crucial for the appearance of 
$\Cset$-independence are the presence of white noise functors 
\cite{Kuem96a,GuMa02a} and a vacuum vector of the underlying deformed 
Fock space which is separating for the considered von Neumann algebras.\\ 
(3) Sources of examples for $\cM_0$-independence are, of course, 
conditional independence in probability theory and random probability 
measures on standard Borel probability spaces. Further examples, 
satisfying the inclusion $\cM_0 \subset \cM_1 \cap \cM_2$, arise from 
amplifications of examples for $\Cset$-independence by tensor product 
constructions. Freeness with amalgamation as well as commuting squares 
from subfactor theory are further sources of $\cM_0$-independence 
(with $\cM_0 \subset \cM_1 \cap \cM_2$). We refer to \cite{HKK04a} for 
a more detailed treatment of some of these examples.\\
(4) $\cM_0$-independence appears, also under the assumption $\cM_0 
\subset \cM_1 \cap \cM_2$, in the work of Junge and Xu on noncommutative 
Rosenthal inequalities \cite{JuXu03a} and within Junge's quantum 
probabilistic approach to embedding Pisier's operator Hilbert space 
$OH$ into the predual of the hyperfinite $III_1$-factor \cite{Jung06a}. 
\end{Remarks}
%%%%%%%%%%%%%%%%%%%%%%%%%%%%%%%%%%%%%%%%%%%%%%%%%%%%%%%%%%%%%%%%%%%%%
%%%%%%%%%%%%%%%%%%%%%%%%%%%%%%%%%%%%%%%%%%%%%%%%%%%%%%%%%%%%%%%%%%%%%
\section{Stationarity and conditional independence/factorizability}
\label{section:ci-cf-stationary}
%%%%%%%%%%%%%%%%%%%%%%%%%%%%%%%%%%%%%%%%%%%%%%%%%%%%%%%%%%%%%%%%%%%%%
This section is devoted to show in Theorem \ref{thm:fact=indep} that 
conditional factorizability implies conditional independence in the 
context of stationarity and under a certain conditioning. 
We close with an illustration of conditional independence and 
conditional factorizability by an algebraic treatment of an infinite 
sequence of exchangeable $\{0,1\}$-valued random variables.
 
Due to the noncommutativity of our setting, there are (at least) 
two natural ways to extend the notions of conditional independence and 
conditional factorizablility (see Definition \ref{def:independence}) from 
two random variables to random sequences indexed by $\Nset_0$. One may 
regard $\Nset_0$ as a set, or as an ordered set (with its natural order).  
%%%%%%%%%%%%%%%%%%%%%%%%%%%%%%%%%%%%%%%%%%%%%%%%%%%%%%%%%%%%%%%%%%%%%
\begin{Definition}\normalfont \label{def:rvs}
%%%%%%%%%%%%%%%%%%%%%%%%%%%%%%%%%%%%%%%%%%%%%%%%%%%%%%%%%%%%%%%%%%%%%
The (identically distributed) random sequence $\mathscr{I}$, given by 
$$
\iota \equiv (\iota_i)_{i\in \Nset_0}\colon (\cA_0,\varphi_0) \to (\cM,\psi),
$$
with canonical filtration $(\cA_I)_{I \subset \Nset_0}$, is said to be 
\begin{enumerate}
\item[(CI)]
\emph{full $\cN$-independent} or \emph{conditionally full inde\-pendent}, 
if $\cA_I$ and $\cA_J$ are $\cN$-inde\-pendent for all 
$I , J \subset\Nset_0$ with $I \cap J = \emptyset$; 
\item[(CI$_\text{o}$)]
\emph{order $\cN$-independent} or \emph{conditionally order independent}, 
if $\cA_I$ and $\cA_J$ are $\cN$-inde\-pendent for all 
$I , J \subset \Nset_0$ with $I < J$ or $I > J$; 
\end{enumerate}
We say that $\scrI$ is 
\begin{enumerate}
\item[(CF)]
\emph{full $\cN$-factorizable} or \emph{conditionally full 
factorizable}, if $\cA_I$ and $\cA_J$ are $\cN$-factor\-izable 
for all 
$I , J \in \subset \Nset_0$ with $I \cap J = \emptyset$;
\item[(CF$_\text{o}$)]
\emph{order $\cN$-factorizable} or \emph{conditionally order factorizable}, 
if $\cA_I$ and $\cA_J$ are $\cN$-factorizable for all 
$I , J \subset \Nset_0$ with $I < J$ or $I > J$; 
\end{enumerate}
\end{Definition}
%%%%%%%%%%%%%%%%%%%%%%%%%%%%%%%%%%%%%%%%%%%%%%%%%%%%%%%%%%%%%%%%%%%%%
We will deliberately drop the attributes `full' or 'order' if we want 
to address conditional independence or conditional factorizability 
only on the informal level or if it is clear from the context whether 
the index set $\Nset_0$ is regarded with order structure or without it. 
Obviously we have the following implications:
\[
\begin{matrix}
 \text{\normalfont{(CI)}} & \Rightarrow &   \text{\normalfont{(CF)}}  \\
              \Downarrow &             & \Downarrow \\
\text{\normalfont{(CI$_\text{o}$)}}   & \Rightarrow & \text{\normalfont{(CF$_\text{o}$)}}. 
\end{matrix}
\]
We record that this gives the following implications in the 
noncommutative extended de Finetti theorem.
%%%%%%%%%%%%%%%%%%%%%%%%%%%%%%%%%%%%%%%%%%%%%%%%%%%%%%%%%%%%%%%%%%%%%%
\begin{proof}[Proof of Theorem \ref{thm:main}, (c) $\Rightarrow$ (c$_\text{o}$)
and (d) $\Rightarrow$ (d$_\text{o}$)]
This is obvious for $\cN = \cM^\tail$ from Definition \ref{def:rvs} and above diagram.
\end{proof}
%%%%%%%%%%%%%%%%%%%%%%%%%%%%%%%%%%%%%%%%%%%%%%%%%%%%%%%%%%%%%%%%%%%%%%

A natural question is to ask if the converse implications in above diagram
are also valid. Actually, we do not know an answer in this generality. But an 
affirmative answer is available for the equivalence of conditional 
independence and conditional factorizability if the random sequence 
$\scrI$ is stationary and $\cN$ contained in the fixed point algebra 
of the corresponding stationary process (see Lemma \ref{lem:corres-rv-endo} 
for this correspondence).
%%%%%%%%%%%%%%%%%%%%%%%%%%%%%%%%%%%%%%%%%%%%%%%%%%%%%%%%%%%%%%%%%%%%%
\begin{Theorem}\label{thm:fact=indep} 
%%%%%%%%%%%%%%%%%%%%%%%%%%%%%%%%%%%%%%%%%%%%%%%%%%%%%%%%%%%%%%%%%%%%% 
Let $\scrM$ be a minimal stationary process and suppose the 
$\psi$-condi\-tioned von Neumann subalgebra $\cN$ satisfies 
$\cN \subset \cM^\alpha$. Then the following are equivalent:
\begin{enumerate}
\item[(CI)] $\scrM$ is full $\cN$-independent;
\item[(CF)] $\scrM$ is full $\cN$-factorizable. 
\end{enumerate}
Furthermore, the following are equivalent under the same assertions:
\begin{enumerate}
\item[(CI$_\text{o}$)] $\scrM$ is order $\cN$-independent;
\item[(CF$_\text{o}$)] $\scrM$ is order $\cN$-factorizable. 
\end{enumerate}
\end{Theorem}
%%%%%%%%%%%%%%%%%%%%%%%%%%%%%%%%%%%%%%%%%%%%%%%%%%%%%%%%%%%%%%%%%%%%%
We will see in Section \ref{section:kol-0-1} that conditional order 
factorizability (CF$_\text{o}$), the weakest of the four properties,
already suffices to identify $\cN$ as the fixed point algebra of 
the endomorphism $\alpha$ which equals moreover the tail algebra.
There it will suffice, due to Theorem \ref{thm:fact=indep}, to 
establish these fixed point characterization results of Kolmogorov type 
on the level of conditional factorizability. Moreover we will benefit 
from this simplification in Section \ref{section:spread2cond-indep} and 
Section \ref{section:spread2indep} when showing that spreadability implies 
conditional independence.  
%%%%%%%%%%%%%%%%%%%%%%%%%%%%%%%%%%%%%%%%%%%%%%%%%%%%%%%%%%%%%%%%%%%%%

%%%%%%%%%%%%%%%%%%%%%%%%%%%%%%%%%%%%%%%%%%%%%%%%%%%%%%%%%%%%%%%%%%%%%
We prepare the proof of Theorem \ref{thm:fact=indep}
by some well known results on approximations.
%%%%%%%%%%%%%%%%%%%%%%%%%%%%%%%%%%%%%%%%%%%%%%%%%%%%%%%%%%%%%%%%%%%%%
\begin{Lemma}\label{lem:continuity}
%%%%%%%%%%%%%%%%%%%%%%%%%%%%%%%%%%%%%%%%%%%%%%%%%%%%%%%%%%%%%%%%%%%%%
Let $x_1,\ldots, x_p \in \ball_1(\cM)$, the unit ball of $\cM$. Suppose
further that each $x_i$ is approximated by a sequence 
$(x_{i,n})_{n \in \Nset} \subset \ball_1(\cM)$ in the strong operator
topology. Then 
\[
x_1 x_2 \cdots x_p = \sotlim_{n \to \infty} x_{1,n}x_{2,n} \cdots x_{p,n}.
\] 
\end{Lemma}
%%%%%%%%%%%%%%%%%%%%%%%%%%%%%%%%%%%%%%%%%%%%%%%%%%%%%%%%%%%%%%%%%%%%%
\begin{proof}
This is evident from the definition of the strong operator topology in
the case $p=2$, since $\|x_{1,n}\| \le 1$ and
\[
x_{1,n}x_{2,n} - x_1 x_2 = x_{1,n}(x_{2,n}-x_2) + (x_{1,n} -x_1)x_2.
\]
The more general case $p>2$ is concluded by induction. 
\end{proof}
%%%%%%%%%%%%%%%%%%%%%%%%%%%%%%%%%%%%%%%%%%%%%%%%%%%%%%%%%%%%%%%%%%%%%
\begin{Lemma}\label{lem:approx-tail}
%%%%%%%%%%%%%%%%%%%%%%%%%%%%%%%%%%%%%%%%%%%%%%%%%%%%%%%%%%%%%%%%%%%%%
Suppose $\scrM$ is a minimal bilateral stationary process
and let the function $N \colon \Nset \to \Zset$ be given.
Then every $a \in \cM^\alpha$ is approximated by a sequence 
$(a_n)_{n\in \Nset} \subset \cM$ in the strong operator topology 
such that 
\[
a_{n} \in \cM_{\{0,1, \ldots, n-1\}+N(n)}  
\quad \text{and} \quad \|a_n\| \le \|a\|.
\]
\end{Lemma}
\begin{proof}
We assume without loss of generality that all $a$ is in the unit ball 
$\ball_1(\cM) \cap \cM^\alpha$.
The *-algebra $\cM^\alg_{\Nset_0}$ is weak*-dense in $\cM$. Thus, by 
Kaplansky's density theorem, $a  \in \cN$ is approximated by 
a sequence $(b_{n})_n\subset \cM^\alg_{\Nset_0} \cap \ball_1(\cM)$ 
in the strong operator topology. Put
\[
a_{n} := \alpha^{N(n)}E_{\cM_{[0,n-1]}} (b_{n}) \in \cM_{\{0, 1,\ldots,n-1\}+N(n)}.
\]
and note that $\alpha$ is an automorphism of $\cM$, since we 
are working in the bilateral setting. We claim that 
\begin{eqnarray}\label{eq:fact=indep-1}
\sotlim_n  a_{n} = a.
\end{eqnarray}
Indeed, the sequence $(E_{\cM_{[0,n-1]}})_n$ is norm bounded and converges 
to $\id_\cM$ in the pointwise strong operator topology; this is clear on 
$\cM_{\Nset_0}^\alg$ and an $\frac{\varepsilon}{2}$-argument gives the 
general case. Thus $\big(E_{\cM_{[0,n-1]}}(b_{n})\big)_n$ converges to 
$a$ in the strong operator topology. We use next the $\psi$-topology which
is induced by the maps $\cM \ni x \mapsto \psi(x^*x)^{1/2}$. Since the strong 
operator topology and the $\psi$-topology coincide on bounded sets, 
\begin{eqnarray*}
\|a_{n} -a\|_\psi  
&=& \|\alpha^{N(n)}\big(E_{\cM_{[0,n-1]}} (b_{n}) -a\big)\|_\psi\\
& = &\|E_{\cM_{[0,n-1]}}(b_{n})-a\|_\psi
\end{eqnarray*}
completes the proof. 
\end{proof}
%%%%%%%%%%%%%%%%%%%%%%%%%%%%%%%%%%%%%%%%%%%%%%%%%%%%%%%%%%%%%%%%%%%%% 
\begin{proof}[Proof of Theorem \ref{thm:fact=indep}]
Only the implications `(CF$_\text{o}$) $\Longrightarrow$ (CI$_\text{o}$)' 
and `(CF) $\Longrightarrow$ (CI)' require a proof, since their reverse 
implications are trivial. We can assume by Theorem \ref{thm:unilateral-bilateral} 
that $\scrM = (\cM,\psi, \alpha, \cM_0\big)$ is a minimal bilateral stationary 
process. This will allow us to approximate elements of $\scrN$ in an appropriate 
manner. Note that full (resp.~ order) $\cN$-factorizability of the family 
$(\cM_I)_{I\subset \Nset_0}$ implies immediately full (resp.\ order)
$\cN$-factorizability of $(\cM_I)_{I \subset \Zset}$ by stationarity; this 
is clear for finite sets $I$ and the general case is done by approximation.   

We need to show that full (resp.\ order) $\cN$-factorizability of 
$(\cM_I)_{I \subset \Nset}$ implies 
\[
E_{\cN}(xy) = E_{\cN}(x) E_{\cN}(y)
\]
for all $x \in \cM_I \vee \cN$ and $y \in \cM_J \vee \cN$ with $I \cap J =\emptyset$
(resp.\ $I <J$). 

For this purpose, we start with bounded sets $I,J \subset\Nset_0$ and 
consider monomials of the form 
\[
x= z_1 a_1 \cdots z_p a_p
\quad \text{and} \quad 
y = z_{p+1} a_{p+1} \cdots z_{2p}a_{2p},
\] 
with $z_i \in \cM_I$, $z_{p+i} \in \cM_J$ and $a_{i}, a_{i+p} \in \cN$ 
($i=1, \ldots, p$). We approximate all $a_{i}$'s in the strong operator 
topology and can assume without loss of generality that all $z_i$'s and 
$a_i$'s are in the unit ball $\ball_1(\cM)$. Let 
$N_i\colon \Nset \to \Zset$ be given function which will be specified 
later. By Lemma \ref{lem:approx-tail}, there exist sequences 
$(a_{i,n})_{n \in \Nset} \subset \ball_1(\cM)$ satisfying 
\begin{eqnarray*}
a_i &=& \sotlim_{n \to \infty} a_{i,n},\\
a_{i,n} &\in& \cM_{\{0,1,\ldots, n-1\}+N_(n)}.
\end{eqnarray*}
Let
\begin{eqnarray*}
x_n &:=& z_1 a_{1,n}z_2a_{2,n}\cdots z_{p}a_{p,n}, \\
y_n &:=& z_{p+1} a_{p+1,n}z_{p+2}a_{p+2,n}\cdots z_{2p}a_{2p,n}.
\end{eqnarray*}
We specify next the choice of the functions
$N_i$. Let $N_i(n) := -n$ and  $N_{p+i}(n):= N$ for $i=1, \ldots, p$, where 
$N > \max{I \cup J}$. 
Note that the sets 
\begin{eqnarray*}
I_n&:=&I \cup \{-n, n+1, \ldots, -1\},\\ 
J_n&:=&J \cup \{N, N+1, \ldots, N+n-1\} 
\end{eqnarray*}
are disjoint if $I$ and $J$ are disjoint; and that $I_n < J_n$ if $I < J$.  
Since $x_{n} \in \cM_{I_n}$ and $y_{n}\in \cM_{J_n}$ we conclude from 
order (resp. full) $\cN$-factorizability that 
\[
 E_{\cN}\big(x_n y_n \big)=
  E_{\cN}\big(x_n \big)   E_{\cN}(y_n\big),
\] 
which entails
\begin{align*}
&E_{\cN}\big(x  y \big) - E_{\cN}(x) E_\cN(y)\\
&=  E_{\cN}\big(x  y  - x_n y_n\big)
+  E_{\cN}\big(x_n)   E_{\cN}\big(y_n - y \big) 
+  E_{\cN}\big(x_n - x \big)   E_{\cN}(y). 
\end{align*}   
We infer from Lemma \ref{lem:continuity} and the $\sot$-$\sot$-continuity
of conditional expectations that the right hand side of this equation 
vanishes for $n\to \infty$ in the strong operator topology. Thus full 
(resp.\ order) $\cN$-factorizability implies, for each $p \in \Nset$,
\begin{align}\label{eq:fact-indep-i}
E_{\cN}&(z_1 a_1 \cdots z_p a_p z_{p+1}a_{p+1} \cdots z_{2p}a_{2p})\notag \\
&= E_{\cN}(z_1 a_1 \cdots z_p a_p) E_{\cN}(z_{p+1}a_{p+1} \cdots z_{2p}a_{2p})
\end{align}
for any $z_1,\ldots ,z_p \in \cM_I$, $z_{p+1}, \cdots z_{2p} \in \cM_J$ and 
$a_1, \cdots a_{2p} \in \cN$ whenever $I$ and $J$ are disjoint (resp.\ ordered)
and bounded. This equality extends by  bilinearity to the 
*-algebras $\cM_I \cup \cN$ and $\cM_J \cup \cN$. (By filling in additional 
factors $\1_\cM$ if necessary we can always achieve that monomials have the 
same number of factors.) Since $\cM_I \cup \cN$ and $\cM_J \cup \cN$ are 
weak* dense in $\cM_I \vee \cN$ and $\cM_J \vee \cN$, the equality 
\eqref{eq:fact-indep-i} extends further to the weak* closure, using Kaplansky's 
density theorem and arguments similar to that in the proof of Theorem 
\ref{thm:mixing}. Finally, another density argument extends the validity of  
\eqref{eq:fact-indep-i} from bounded disjoint sets $I$ and $J$ to possibly 
unbounded disjoint sets. 
\end{proof}

%%%%%%%%%%%%%%%%%%%%%%%%%%%%%%%%%%%%%%%%%%%%%%%%%%%%%%%%%%%%%%%%%%%%%
\begin{Remark}\normalfont
At the time of this writing and in the generality of our setting, we 
have no information about the validity of the remaining implications 
(CI$_\text{o}$) $\Rightarrow$ (CI) and (CF$_\text{o}$) $\Rightarrow$ 
(CF), even under the assumptions of stationarity  and $\cN \simeq \Cset$.
In particular, we do not know if an infinite stationary random sequence exists
which is conditionally order independent, but fails to be conditionally
full independent.
\end{Remark}
%%%%%%%%%%%%%%%%%%%%%%%%%%%%%%%%%%%%%%%%%%%%%%%%%%%%%%%%%%%%%%%%%%%%%
%%%%%%%%%%%%%%%%%%%%%%%%%%%%%%%%%%%%%%%%%%%%%%%%%%%%%%%%%%%%%%%%%%%%%
We continue with an illustration of above concepts of conditional 
independence and conditional factorizability for stationary random 
sequences. The example is the von Neumann algebraic reformulation of 
infinite sequence of zero-one-valued random variables, as they have 
been the subject of de Finetti's pioneering investigations on 
exchangeability \cite{Fine31a}. We will observe in this example why it 
is too restrictive to assume that $\cN$ is contained in the image of 
random variables.
%%%%%%%%%%%%%%%%%%%%%%%%%%%%%%%%%%%%%%%%%%%%%%%%%%%%%%%%%%%%%%%%%%%%%
\begin{Example}\normalfont  \label{example:coin-tosses} 
Let $(\cA_0, \varphi_0)$ be given by 
\[
\cA_0 = \Cset^2 \quad \text{and}  \quad \varphi_0 = \trace_{p}
\] 
with $\trace_{p}((a_1,a_2)) = p a_1 + (1-p) a_2$ for some fixed 
$p \in (0,1)$.  We realize the probability space $(\cM, \varphi)$ 
as a mixture of infinite coin tosses with respect to some probability 
measure $\nu$ on the standard measurable space $([0,1], \Sigma)$, 
assuming $\nu(\{0\}) = \nu(\{1\})=0$ and $\nu(\{p\}) < 1$ for any 
$p \in (0,1)$:   
\begin{align*}
\cM  &= \int^\oplus_{[0,1]} \cM(p)  d\nu(p),                
& \cM(p) &= \bigotimes_{n \in \Nset_0} \Cset_2, \\
\psi &= \int^\oplus_{[0,1]} \psi(p) d \nu(p),          
& \psi(p) &= \bigotimes_{n \in \Nset_0} \trace_p.  
\end{align*}
Here denotes $\cM(p)$ the infinite von Neumann algebraic tensor product 
of $\Cset_2$ with respect to the infinite tensor product state on $\psi(p)$ 
which are obtained by passing through the GNS construction starting from 
the *-algebra $\bigcup_{k \in \Nset} \otimes_{k=0}^n \Cset_2$ equipped 
with the product state $\bigcup_{k \in \Nset} \otimes_{k=0}^n \trace_p$.
We refer the reader to \cite{Take79a} for further information on direct 
integrals of von Neumann algebras and states.\\   
The random variable $\iota_i\colon (\cA_0,\varphi_0) \to (\cM, \psi)$, 
with $i \in \Nset_0$, is defined by the constant embedding  of 
$a \in \Cset_2$ into the $i$-th factor of each fiber of the direct integral:
\[
\iota_i(a) 
= \int^\oplus_{[0,1]} \underbrace{\1_{\cA_0} \otimes 
  \cdots \otimes \1_{\cA_0}}_{\text{$i$ factors}} \otimes\, a \otimes  
  \1_{\cA_0} \otimes \cdots  d \nu(p)
\]
Finally, we put 
\[
\cN :=  \int_{[0,1]}^\oplus \Cset \1_{\cM(p)} d \nu(p)  \simeq L^\infty([0,1], \nu).
\] 
Note that our assumptions on the measure $\nu$ imply $\cN \not \simeq \Cset$.\\ 
%%%%%%%%%%%%%%%%%%%%%%%%%%%%%%%%%%%%
The canonical filtration $(\cA_I)_{I \subset \Nset_0}$ generated 
by the random sequence $\iota\equiv (\iota_i)_{i \in \Nset_0}$ is defined by
$$
\cA_I = \bigvee_{i \in I} \iota_i(\cA_0).
$$
The random sequence $\iota$ is minimal, i.e.\ we have 
$$
\cM = \bigvee_{n \in \Nset_0} \iota_i(\cA_0).
$$  
This follows if we can ensure that $\bigvee_{n \in \Nset_0} \iota_i(\cA_0)$ 
contains $\cN$. Indeed, Kakutani's theorem entails that the family of 
infinite product states $\{\psi(p)\}_{p \in (0,1)}$ is mutually disjoint 
\cite{Hida80a}. We conclude from this that every element $x \in \cN \simeq 
L^\infty([0,1], \nu)$ can be approximated by a bounded sequence 
$(x_n)_{n \in \Nset} \subset \bigcup_{i \in \Nset_0} \iota_i(\cA_0)$ 
in the weak operator topology. This implies the minimality of the random 
sequence.\\
%%%%%%%%%%%%%%%%%%%%%%%%%%%%%%%%%%%%%%%%%%%%  
An elementary computation shows $\cA_I \simeq \Cset^{2^{|I|}}$ for any 
finite set $I \subset \Nset_0$.  In the case of an infinite set $I$, we 
restrict the family of infinite product states $\{\psi(p)\}_{p \in (0,1)}$ 
to $\cA_I$ and conclude again by the Kakutani theorem \cite{Hida80a} that 
these restricted states are mutually disjoint. This implies that the von 
Neumann algebra $\cA_I$ contains a copy of $\cN$ whenever $|I|= \infty$.\\
%%%%%%%%%%%%%%%%%%%%%%%%%%%%%%%%%%%%%%%%%%%%%
Now it is straightforward to verify the conditional full factorizability (CF)
\[
E_{\cN}(xy)= E_{\cN}(x)E_{\cN}(y)
\]
for all $x \in \cA_I$ and $y \in \cA_J$ with disjoint subsets 
$I,J \subset \Nset_0$. Since all von Neumann algebras are commutative, 
it is immediate from the module property of conditional expectation that 
(CF) upgrades to (CI), i.e.
\[
E_{\cN}(xy)= E_{\cN}(x)E_{\cN}(y)
\]
for all $x \in \cA_I \vee \cN$ and $y \in \cA_J \vee \cN$ with disjoint 
subsets $I,J \subset \Nset_0$. Thus the random sequence 
$(\iota_i)_{i \in \Nset_0}$ is full $\cN$-independent. But 
$\cN \not \subset \cA_I \cap \cA_J$ if one of the sets $I$ or $J$ is 
finite.     
\end{Example}
%%%%%%%%%%%%%%%%%%%%%%%%%%%%%%%%%%%%%%%%%%%%%%%%%%%%%%%%%%%%%%%%%%%%%
\begin{Remark}\normalfont
There are *-algebraic, C*-algebraic and W*-algebraic 
approaches to noncommutative probability and it is instructive to 
compare them at the hand of Example \ref{example:coin-tosses}. Of course, 
the *-algebras $\cA_I^\alg := \bigcup_{i \in I}\iota_i(\cA_0)$ as well as 
its norm-closure are contained in the von Neumann algebra $\cA_I$. The 
latter contains a copy of $L^\infty([0,1], \nu)$ if $|I|$ is infinite, 
but $\cA_I^\alg$ and its norm closure do not.  
\end{Remark}
%%%%%%%%%%%%%%%%%%%%%%%%%%%%%%%%%%%%%%%%%%%%%%%%%%%%%%%%%%%%%%%%%%%%%%

%%%%%%%%%%%%%%%%%%%%%%%%%%%%%%%%%%%%%%%%%%%%%%%%%%%%%%%%%%%%%%%%%%%%%
\section{Noncommutative i.i.d.~sequences may be non-stationary}
\label{section:non-stationary}
%%%%%%%%%%%%%%%%%%%%%%%%%%%%%%%%%%%%%%%%%%%%%%%%%%%%%%%%%%%%%%%%%%%%%
It is folklore in classical probability and free probability that
independence resp.\ freeness of an identically distributed random 
sequence implies stationarity. But this implication fails in our 
broader context of noncommutative independence. 
%%%%%%%%%%%%%%%%%%%%%%%%%%%%%%%%%%%%%%%%%%%%%%%%%%%%%%%%%%%%%%%%%%%%%%
\begin{Theorem}
%%%%%%%%%%%%%%%%%%%%%%%%%%%%%%%%%%%%%%%%%%%%%%%%%%%%%%%%%%%%%%%%%%%%%%
There exist full $\Cset$-independent identically distributed 
random sequences $\scrI$ which fail to be stationary. 
\end{Theorem}
%%%%%%%%%%%%%%%%%%%%%%%%%%%%%%%%%%%%%%%%%%%%%%%%%%%%%%%%%%%%%%%%%%%%%%
\begin{proof}[Proof, in particular of Theorem \ref{thm:main} (c) 
$\not \Leftarrow$ (d) and (c$_{\text{o}}$) $\not \Leftarrow$ (d$_{\text{o}}$)]
See Example \ref{example:non-stationarity-1} or Example 
\ref{example:non-stationarity-2} below. 
Since full $\Cset$-independence implies order $\Cset$-independence 
we have also shown (c$_{\text{o}}$) $\not \Leftarrow$ (d$_{\text{o}}$).  
\end{proof}

Let us first outline our strategy to produce such examples. Recall from 
the introduction that an infinite random sequence $\scrI$ with random variables
\[
(\iota_n)_{n \ge 0}\colon (\cA_0,\varphi_0) \to (\cM,\psi)
\]
is automatically identically distributed. Suppose now that $\scrI$
is stationary and $\Cset$-independent. Our goal is to `perturbate'
the random variables $\iota_n$ such that $\Cset$-independence is 
preserved, but stationary is obstructed. This can be done in two ways, 
for the domain or the codomain of each random variable $\iota_n$.     

%%%%%%%%%%%%%%%%%%%%%%%%%%%%%%%%%%%%%%%%%%%%%%%%%%%%%%%%%%%%%%%%%%%%%%
\begin{Example}[Perturbation of codomain]\normalfont \label{example:non-stationarity-1}
%%%%%%%%%%%%%%%%%%%%%%%%%%%%%%%%%%%%%%%%%%%%%%%%%%%%%%%%%%%%%%%%%%%%%%
Consider $(\cR, \trace)$, the hyperfinite $II_1$-factor equipped with 
its normalized trace. Let $(M_m, \trace_m)$ be the complex 
$m\times m$-matrices equipped with the normalized trace. The canonical 
embeddings 
\[
M_2 \ni x \mapsto  \iota_n(x):=\1_{M_2} \otimes \cdots \otimes \1_{M_2} 
\underset{\text{$n$-th position}}{\otimes\quad x\quad \otimes} 
\1_{M_2} \otimes \cdots 
\]
define the random sequence $\scrI$ with random variables
\[
(\iota_n)_{n \ge 0} \colon (M_2, \trace_2) \to (\cR, \trace).
\]
It is easily verified that $\scrI$ is $\Cset$-independent and stationary. 
We will deform this random sequence to obtain a non-stationary random sequences 
as follows. Under the canonical identification of $M_2 \otimes M_2$ and $M_4$, 
the unitary matrix
\[
U_\omega = \begin{bmatrix}
     1 & 0 & 0 &0 \\
     0 & 0 & 1 & 0 \\
     0 & 1 & 0 & 0 \\
     0 & 0 & 0 & \omega
    \end{bmatrix},     
\qquad |\omega|=1,
\] 
defines the trace-preserving automorphism $x \mapsto U_\omega x U^*_\omega$ 
of $M_2 \otimes M_2$. It is well known in subfactor theory that the inclusions
\begin{eqnarray*}
\begin{matrix}
U_\omega (M_2 \otimes \1_{M_2})U^*_\omega & \subset & M_2 \otimes M_2  \\  
 \cup &&                                   \cup \\
\Cset \,\1_{M_2 \otimes M_2} & \subset & M_2 \otimes \1_{M_2} 
\end{matrix}
\end{eqnarray*}
form a commuting square \cite{Jone91a,Rupp95a,JoSu97a}. We canonically amplify 
this automorphism to the automorphism $\gamma_\omega \in \Aut{\cR, \trace}$ 
which acts trivial on all higher tensor product factors. Consider now the random 
sequence $\scrI^{(\omega)}$ with random variables $\big(\iota^{(\omega)}_n\big)_{n \ge 0}$ 
defined by
\begin{eqnarray*}
\iota^{(\omega)}_n := 
\begin{cases}
\iota_n & \text{if $n \neq 1$}\\
\gamma_\omega \iota_0 & \text{if $n=1$}.
\end{cases} 
\end{eqnarray*}
Note that $\scrI^{(1)}$ is the random sequence $\scrI$. 
Clearly $\scrI^{(\omega)}$ is identically distributed for any unimodular 
$\omega \in \Cset$. We note that the von Neumann algebras 
$\iota_n^{(\omega)}(M_2)$ mutually commute for $n \neq 1$. So do 
$\iota^{(\omega)}_1(M_2)$ and $\iota_n^{(\omega)}(M_2)$ for $n \ge 2$. 
We conclude from this that $\scrI^{(\omega)}$ is full $\Cset$-independent.  
But we calculate for $a = \begin{bmatrix} 0 & 1 \\ 1 & 0 \end{bmatrix}$ that 
\[
\trace\big( \iota^{(\omega)}_0(a) \iota^{(\omega)}_1(a) 
            \iota^{(\omega)}_0(a) \iota^{(\omega)}_1(a) \big) 
= \frac{1}{2} (\omega + \overline{\omega}),
\]
and
\[
\trace\big( \iota^{(\omega)}_2(a) \iota^{(\omega)}_3(a) 
            \iota^{(\omega)}_2(a) \iota^{(\omega)}_3(a) \big) 
= 1.
\]
This leads us to the conclusion that $\scrI^{(\omega)}$ is stationary 
if and only if $\omega = 1$. 
\end{Example}
%%%%%%%%%%%%%%%%%%%%%%%%%%%%%%%%%%%%%%%%%%%%%%%%%%%%%%%%%%%%%%%%%%%%%%
\begin{Remark}\normalfont
Example \ref{example:non-stationarity-1} illustrates that 
the distribution of two $\Cset$-independent (identically distributed) 
random variables does not determine their joint distribution. This is 
in contrast to two distinguished examples for $\Cset$-independence, 
tensor independence and free independence. See \cite{Spei97a,BGSc02a}
for further information on the related universality properties.   
\end{Remark}
%%%%%%%%%%%%%%%%%%%%%%%%%%%%%%%%%%%%%%%%%%%%%%%%%%%%%%%%%%%%%%%%%%%%%%

%%%%%%%%%%%%%%%%%%%%%%%%%%%%%%%%%%%%%%%%%%%%%%%%%%%%%%%%%%%%%%%%%%%%%%%
We sketch next how local perturbations of random variables on their 
domain are capable to produce such effects.
%%%%%%%%%%%%%%%%%%%%%%%%%%%%%%%%%%%%%%%%%%%%%%%%%%%%%%%%%%%%%%%%%%%%%
Suppose the minimal stationary random sequence $\scrI$ with random
variables
\[
(\iota_n)_{n \in \Nset_0} \colon (\cA_0, \varphi_0) \to (\cM,\psi)
\]
is $\Cset$-independent (in the ordered or full sense). Furthermore, 
let $\gamma\equiv(\gamma_n)_{n \ge 0} \subset \Aut{\cA_0,\varphi_0}$ 
be a sequence of `local perturbations'. Then we can associate to each 
sequence $\gamma$ a random sequence $\scrI^{(\gamma)}$ by putting
\[
\iota^{(\gamma)}_n := \iota_n \circ \gamma_n.  
\] 
The random sequence $\scrI^{(\gamma)}$ is again minimal and 
$\Cset$-independent. Suppose that there is a sequence $\gamma$ with
\[
(\iota_0, \iota_1, \ldots, \iota_{n-1}, \iota_n, \iota_{n+1}\ldots) 
\overset{\distr}{\neq} 
(\iota_0^{(\gamma)}, \iota_1^{(\gamma)}, \ldots , 
\iota_{n-1}^{(\gamma)}, \iota_n, \iota_{n+1}\ldots ) 
\]
for some $n \in \Nset$. We conclude immediately that the 
random sequence on the right hand side fails to be stationary, but it 
is still identically distributed and enjoys $\Cset$-independence.  

\begin{Example}[Perturbation of domain]\normalfont \label{example:non-stationarity-2}
Let $\ell^2(\Nset)$ be the real Hilbert space of square-summable sequences 
and consider the $q$-Gaussian field $\Gamma_q{(\ell^2(\Nset))}$ for some fixed
$0 < q < 1$. These fields are the von Neumann algebra generated by $q$-Gaussian 
field operators $\omega_q(f)$, $f\in \ell^2(\Nset$), acting on the $q$-deformed 
Fock space $\cF_q(\ell^2(\Nset))$ (see \cite{BKS97a, GuMa02a} for further details). 
$\Gamma_q{(\ell^2(\Nset))}$ is a non-hyperfinite $II_1$-factor and we
denote its normalized trace by $\trace_q$. The second quantization of the canonical 
unilateral shift on $\ell^2(\Nset)$ provides us with a unital $\trace_q$-preserving 
endomorphism $\alpha$ of $\Gamma_q{(\ell^2(\Nset))}$. Identify $\Rset$ with the subspace 
generated by the first coordinate of $\ell^2(\Nset)$. Doing so we obtain the abelian 
von Neumann subalgebra $\Gamma_q(\Rset) \subset \Gamma_q(\ell^2(\Nset))$ and we denote
the restriction of $\trace_q$ to this subalgebra by the same symbol. Now it is
straightforward to see that
\[
\iota_n:= \alpha^n|_{\Gamma_q(\Rset)}
\]
defines a full $\Cset$-independent random sequence $\scrI$ with random variables
\[
(\iota_{n\ge 0})\colon (\Gamma_q(\Rset), \trace_q) \to   (\Gamma_q(\ell^2(\Nset)), \trace_q),  
\]
which is of course stationary. Let $\gamma \in \Aut{\Gamma_q(\Rset), \trace_q}$ be fixed
and consider the random sequence $\scrI_\gamma$ which is obtained from perturbating the first 
random variable of $\scrI$:
\begin{eqnarray*}
\iota_n^{\gamma} := 
\begin{cases}
\iota_n & \text{if $n \neq 1$}\\
\iota_0 \circ \gamma & \text{if $n=0$}.
\end{cases} 
\end{eqnarray*} 
The central result by van Leeuwen and Maassen on the obstruction for $q$-deformation 
of the convolution product can be reformulated as:
\begin{Theorem}[\cite{LeMa96a}]\label{thm:q-conv}
Let $0 < q <1$. There exists a `perturbation' $\gamma \in \Aut{\Gamma_q(\Rset), \trace_q}$ 
such that
\[
\trace_q\left(\Big(\omega_q(f) + \alpha\big(\omega_q(f)\big)\Big)^4 \right) \neq 
\trace_q\left(\Big(\gamma\big(\omega_q(f)\big) + \alpha \big(\omega_q(f)\big)\Big)^4 \right) 
\]
for $0 \neq f \in \Rset$.
\end{Theorem}
Note that $\omega(f), \gamma\big(\omega(f)\big)$ and $\alpha\big(\omega(f)\big)$ have identical
distributions and each of the first two random variables is $\Cset$-independent from the third one.
Thus the knowledge of the individual distributions of $\Cset$-independent random variables does not
completely determine their joint distributions; this depends on the concrete realization of the
random variables.   

The `perturbation' $\gamma$ is of constructed in \cite{LeMa96a}
starting from a $\mu$-preserving point transformation on the spectrum of the 
(selfadjoint) $q$-Gaussian field operator $\omega_q(f)$, for some fixed $f\in \Rset$, 
where $\mu$ is induced by the spectral measure of $\omega_q(f)$ with respect to 
$\trace_q$. 
%%%%%%%%%%%%%%%%%%%%%%%%%%%%%%%%%%%%%%%%%%%%%%%%%%%%%%%%%%%%%%%%%%%%%%%%%%%%%
\begin{Corollary}
%%%%%%%%%%%%%%%%%%%%%%%%%%%%%%%%%%%%%%%%%%%%%%%%%%%%%%%%%%%%%%%%%%%%%%%%%%%%%
$\scrI_\gamma$ is full $\Cset$-independent and non-stationary. 
\end{Corollary}
\begin{proof}
It is immediate from its construction that $\scrI$ is full $\Cset$-independent. The perturbation
$\gamma$ of the domain of the first random variable does not effect its range.
Thus $\scrI_\gamma$ is also full $\Cset$-independent.

Let $a := \omega_q(f)$ for notational convenience. A straightforward computation
yields for the left hand side of the inequality in Theorem \ref{thm:q-conv} that
\begin{eqnarray*}
\trace_q\left(\big(a + \alpha(a)\big)^4 \right) 
=  2 \trace_q(a^4) + 4 \trace_q(a^2) \trace_q(a^2) + 2\trace_q\big(a\alpha(a)a\alpha(a)\big). 
\end{eqnarray*}
(Expand the product; use traciality, $\Cset$-independence, $\trace_q \circ\, \alpha = \trace_q$ and
the centredness of $a$.)
Similarly, the right hand side of this inequality simplifies to 
\begin{eqnarray*}
\trace_q\left(\big(\gamma(a) + \alpha(a)\big)^4 \right) 
=  2 \trace_q(a^4) + 4 \trace_q(a^2) \trace_q(a^2) + 2\trace_q\big( \gamma(a)\alpha(a)\gamma(a)\alpha(a)\big). 
\end{eqnarray*}
Since 
\[
\trace_q\big( \gamma(a)\alpha(a)\gamma(a)\alpha(a)\big) \neq \trace_q\big(a\alpha(a)a\alpha(a)\big)
\]
by Theorem \ref{thm:q-conv}, we have
\[
\trace_q\big( \gamma(a)\alpha(a)\gamma(a)\alpha(a)\big) \neq \trace_q\big(\alpha(a)\alpha^2(a)\alpha(a)\alpha^2(a)\big)
\]
and consequently
$
(\iota_0 \circ \gamma, \iota_1, \iota_2, \ldots) 
\overset{\distr}{\neq} 
(\iota_1, \iota_2, \iota_3, \ldots). 
$
\end{proof}
%%%%%%%%%%%%%%%%%%%%%%%%%%%%%%%%%%%%%%%%%%%%%%%%%%%%%%%%%%%%%%%%%%%%%%%%%%%%%
\end{Example}
%%%%%%%%%%%%%%%%%%%%%%%%%%%%%%%%%%%%%%%%%%%%%%%%%%%%%%%%%%%%%%%%%%%%%%%%%%%%%
The invariance of all finite joint distributions of an identically 
distributed random sequence under all local automorphisms seems to be
a very strong condition. If the von Neumann algebra $\scrM$ is abelian 
and $\gamma \in \Aut{\cA_0,\varphi_0}$ ergodic, such a local invariance 
property implies the $\Cset$-independence of the random sequence by an 
application of the mean ergodic theorem.  In the noncommutative context,
this observation invites to introduce `top-order $\Cset$-independence'
for a random sequence $\scrI$, i.e.\ the von Neumann algebras 
$\bigvee_{k<n} \iota_k(\cA_0)$ and $\iota_n(\cA_0)$ are $\Cset$-independent 
for all $n \in \Nset$. If $G \subset \Aut{\cA_0,\varphi_0}$ is an amenable 
ergodic subgroup such that, for all $n \in \Nset$, 
\[
\psi\big(x  \, \iota_n(a)\big) = \psi\big(x  \, \iota_n(\gamma (a))\big)
\] 
for all $x \in \bigvee_{k<n} \iota_k(\cA_0)$ and $\gamma \in G$, then the 
random sequence $(\iota_n)_{n \in \Nset_0}$ is already `top-order 
$\Cset$-independent'.   

\begin{Question}\normalfont
Suppose that a minimal random sequence $\scrI$ with random
variables 
\[
(\iota_n)_{n \in \Nset_0} \colon (\cA_0, \varphi_0) \to (\cM,\psi) 
\] 
has joint distributions which are invariant under all `local perturbations'
$(\gamma_n)_{n \in \Nset} \subset \Aut{\cA_0, \varphi_0}$:
\[
(\iota_0, \iota_1, \iota_2, \ldots ) \overset{\distr}{=} 
(\iota_0 \circ \gamma_0, \iota_1 \circ \gamma_1, \iota_2 \circ\gamma_2, \ldots ).  
\]
Does the ergodicity of $\Aut{\cA_0, \varphi_0}$ imply that $\scrI$ is full 
$\Cset$-independent? And if so, can one show that this $\Cset$-independence 
must be either tensor independence or free independence? 
\end{Question}

%%%%%%%%%%%%%%%%%%%%%%%%%%%%%%%%%%%%%%%%%%%%%%%%%%%%%%%%%%%%%%%%%%%%%%%%%%%%%%%%%%%%
\section{Stationarity with strong mixing and\\ noncommutative Bernoulli shifts}
\label{section:kol-0-1}
%%%%%%%%%%%%%%%%%%%%%%%%%%%%%%%%%%%%%%%%%%%%%%%%%%%%%%%%%%%%%%%%%%%%%%%%%%%%%%%%%%%%
We provide a noncommutative generalization of the Kolmogorov Zero-One Law. 
Furthermore we show that conditional factorizability implies strong mixing 
in the context of stationarity. This leads us to a noncommutative 
generalization of classical Bernoulli shifts.  
%%%%%%%%%%%%%%%%%%%%%%%%%%%%%%%%%%%%%%%%%%%%%%%%%%%%%%%%%%%%%%%%%%%%%
\begin{Theorem}\label{thm:kol-0-1}
%%%%%%%%%%%%%%%%%%%%%%%%%%%%%%%%%%%%%%%%%%%%%%%%%%%%%%%%%%%%%%%%%%%%%
Let $\scrI$ be an order $\cN$-factorizable random sequence where $\cN$ 
is a $\psi$-conditioned von Neumann subalgebra of $\cM^\tail$. Then it 
holds $\cN = \cM^\tail$. In particular, an order $\Cset$-independent 
random sequence has a trivial tail algebra. 
\end{Theorem}
The last assertion is a noncommutative Kolmogorov Zero-One Law. Note 
also that order $\cN$-factorizability (CF$_{\text{o}}$) is implied by
(CI$_{\text{o}}$), (CF) or  (CI).   
%%%%%%%%%%%%%%%%%%%%%%%%%%%%%%%%%%%%%%%%%%%%%%%%%%%%%%%%%%%%%%%%%%%%%
\begin{proof}
%%%%%%%%%%%%%%%%%%%%%%%%%%%%%%%%%%%%%%%%%%%%%%%%%%%%%%%%%%%%%%%%%%%%%
Assume without loss of generality that $\scrI$ is minimal. We show 
first that $\cM$ and $\cM^\tail$ are $\cN$-independent. Let 
$a \in \cM^\tail$ and $x \in \cM^\alg_{\Nset_0}$. Thus there exists 
some bounded subset $J \subset \Nset_0$ such that $x \in \cM_J$. 
Because $\cM^\tail \subset \cM_{[n, \infty)}$ for all $n \in \Nset_0$, 
we can assume $J < [n, \infty)$. Consequently, the order 
$\cN$-factorizability implies
\[
E_{\cN}(a x) = E_{\cN}(a) E_{\cN}(x).
\]
Now let $x \in \cM$. By minimality and Kaplansky's density theorem, 
there exists a bounded sequence $(x_k)_{k \in \Nset}$ in  
$\cM^\alg_{\Nset_0}$ of $\cM$ such that $x = \wotlim_{k}x_k$. Note 
that, for all $k$, we have $x_k \in \cM_{J_k}$ with some bounded 
subset $J_k$. We conclude that, for any $y\in \cM$, 
\begin{eqnarray*}
\psi(y E_{\cN}(a x)) 
&=& \lim_{k} \psi(y E_{\cN}(a x_k)) \\
&=& \lim_{k} \psi(y E_{\cN}(a) E_{\cN}(x_k))\\ 
&=& \psi(y E_{\cN}(a) E_{\cN}(x)).    
\end{eqnarray*}
This gives the factorization
\begin{eqnarray}\label{eq:kol-0-1-i}
E_{\cN}(a x) = E_{\cN}(a) E_{\cN}(x)
\end{eqnarray}
for all $a \in \cM^\tail$ and $x\in \cM$. 
We claim that this factorization implies the 
$\cN$-independence of $\cM^\tail$ and $\cM$. Indeed, 
the $\psi$-preserving conditional $E_{\cM^\tail}$ 
from $\cM$ onto $\cM^\tail$ exist since $\cM^\tail$ 
is globally $\sigma_t^\psi$-invariant. The latter 
is easily concluded from the fact that the
ranges of the random variables $\iota_n$ are 
$\psi$-conditioned and the definition of $\cM^\tail$. 
We are left to verify that \eqref{eq:kol-0-1-i} extends 
to elements $a \in \cM^\tail \vee \cN$ and 
$x\in \cM \vee \cN$. But this is evident, because 
$\cN \subset \cM$ and $\cN \subset \cM^\tail$. Thus 
$\cM$ and $\cM^\tail$ are $\cN$-independent.   

To prove $\cN = \cM^\tail$, we are left to show the 
inclusion $\cM^\tail \subset \cN$. We infer from the 
$\cN$-independence of $\cM$ and $\cM^\tail$ that 
$\cM^\tail$ and $\cM^\tail$ are $\cN$-independent. 
We use the module property of conditional expectations 
and $\cN$-independence to get, for every $x \in \cM^\tail$,
\[
E_{\cN}((x-E_{\cN}(x))^*(x-E_{\cN}(x)) 
= E_{\cN}(x^*x) - E_{\cN}(x^*)E_{\cN}(x)
=0.
\]  
Now the faithfulness of $E_{\cN}$ implies $x= E_{\cN}(x)$ 
and thus $\cM^\tail \subset \cN$. 

The last assertion is clear since order $\Cset$-factorizability 
and order $\Cset$-independence are equivalent 
(see Definition \ref{def:rvs}.
\end{proof}
%%%%%%%%%%%%%%%%%%%%%%%%%%%%%%%%%%%%%%%%%%%%%%%%%%%%%%%%%%%%%%%%%%%%%
\begin{Remark}\normalfont
%%%%%%%%%%%%%%%%%%%%%%%%%%%%%%%%%%%%%%%%%%%%%%%%%%%%%%%%%%%%%%%%%%%%%
The assumptions in Theorem \ref{thm:kol-0-1} can be further weakened
since an inspection of its proof shows that only the ranges of the 
random variables matter. It suffices that the probability space 
$(\cM,\psi)$ is equipped with an order $\cN$-factorizable family of 
$\psi$-conditioned von Neumann subalgebras $(\cM_k)_{k\in \Nset}$.  
 \end{Remark}
%%%%%%%%%%%%%%%%%%%%%%%%%%%%%%%%%%%%%%%%%%%%%%%%%%%%%%%%%%%%%%%%%%%%%
It is well-known that the Kolmogorov Zero-One Law implies strong 
mixing properties of an independent stationary random sequence. Here 
we are interested in a conditioned noncommutative version of this 
classical result. It is convenient to formulate it in terms of the 
minimal stationary process $\scrM$ associated to a stationary random 
sequence $\scrI$.
%%%%%%%%%%%%%%%%%%%%%%%%%%%%%%%%%%%%%%%%%%%%%%%%%%%%%%%%%%%%%%%%%%%%%
\begin{Definition}\normalfont \label{def:strong-mixing}
%%%%%%%%%%%%%%%%%%%%%%%%%%%%%%%%%%%%%%%%%%%%%%%%%%%%%%%%%%%%%%%%%%%%%
A stationary process $\scrM$ or its endomorphism $\alpha$ is 
said to be \emph{strongly mixing over $\cN$} if, for any 
$x\in \cM$, 
\[
\wotlim_{n \to \infty} \alpha^n(x) = E_{\cN}(x).  
\]   
Here $\cN$ is a $\psi$-conditioned von Neumann subalgebra of $\cM$. 
\end{Definition}
%%%%%%%%%%%%%%%%%%%%%%%%%%%%%%%%%%%%%%%%%%%%%%%%%%%%%%%%%%%%%%%%%%%%%
\begin{Theorem}\label{thm:mixing}
%%%%%%%%%%%%%%%%%%%%%%%%%%%%%%%%%%%%%%%%%%%%%%%%%%%%%%%%%%%%%%%%%%%%%
Let the minimal stationary process $\scrM$ be order $\cN$-factorizable 
for the $\psi$-conditioned subalgebra $\cN$ of $\cM^\alpha$. Then 
$\alpha$ is strongly mixing over $\cN$. Moreover we have
\[
\cN = \cM^\alpha = \cM^\tail.
\]
In particular, these three subalgebras are trivial if $\scrM$ is 
order $\Cset$-independent. 
%%%%%%%%%%%%%%%%%%%%%%%%%%%%%%%%%%%%%%%%%%%%%%%%%%%%%%%%%%%%%%%%%%%%%
\end{Theorem}
%%%%%%%%%%%%%%%%%%%%%%%%%%%%%%%%%%%%%%%%%%%%%%%%%%%%%%%%%%%%%%%%%%%%%
The condition $\cN \subset \cM^\alpha$ is non-trivial if $\cM^\tail  
\not \simeq \Cset$ (see Remark \ref{rem:mixing-nontrivial}). 
%%%%%%%%%%%%%%%%%%%%%%%%%%%%%%%%%%%%%%%%%%%%%%%%%%%%%%%%%%%%%%%%%%%%%
\begin{proof}
%%%%%%%%%%%%%%%%%%%%%%%%%%%%%%%%%%%%%%%%%%%%%%%%%%%%%%%%%%%%%%%%%%%%%
Since $\cM^\alpha \subset \cM^\tail$, we conclude $\cN = \cM^\alpha
= \cM^\tail$ from Theorem \ref{thm:kol-0-1}. We are left to prove the 
mixing properties.
Suppose $x\in \cM_I$ and $y \in \cM_J$ for bounded sets 
$I, J \subset \Nset_0$. One calculates 
\begin{eqnarray*}
  \lim_{n \to \infty} \psi(y^* \alpha^n(x)) 
&=& \lim_{n \to \infty} \psi(E_{\cN}(y^* \alpha^n(x)))\\
&=& \lim_{n \to \infty} \psi(E_{\cN}(y^*) E_{\cN}(\alpha^n(x)))\\
&=& \psi( E_{\cN}(y^*) E_{\cN}(x))\\
&=& \psi( y^* E_{\cN}(x)).
\end{eqnarray*}
Here we used that $J < (I + n)$ for $n$ sufficiently large and applied 
order $\cN$-factorizability to obtain the second equality. The third 
equality uses that $\cN \subset \cM^\alpha$ implies 
$E_{\cN} \circ \alpha = E_{\cN}$. 
   
To extend these equations to arbitrary $x,y \in \cM$, we use the 
minimality of the stationary process and approximate $x$ and $y$
by bounded sequences $(x_i)_i$ and, respectively, $(y_i)_i$ from 
the *-algebra $\cM^\alg_{\Nset_0}$ in the strong operator topology. 
Since
\begin{eqnarray*}
\psi(y^* \alpha^n(x)) 
&=& \psi\big((y - y_i)^*\alpha^n(x)\big)  
   +\psi\big(y_i^* \alpha^n(x-x_i)\big)
     +\psi\big(y_i^* \alpha^n(x_i)\big)
\end{eqnarray*}          
and since the estimates
\begin{eqnarray*}
|\psi\big((y - y_i)^*\alpha^n(x)\big)| 
&\le& \psi(|y - y_i|^2)^{1/2} \psi(|x|^2)^{1/2},\\ 
|\big(y_i^* \alpha^n(x-x_i)\big)| 
&\le & \psi(|y_i|^2)^{1/2}  \psi(|x-x_i|^2)^{1/2}
\end{eqnarray*}
are uniform in $n$, we conclude the convergence of 
$\psi(y^* \alpha^n(x))$ to $\psi(y^*E_{\cM^\tail}(x))$
by an $\varepsilon/3$-argument. Now the claimed mixing 
property follows from the norm density of the 
functionals $\set{\psi(y \,\cdot\,)}{y \in \cM}$
in $\cM_*$ and the boundedness of the set 
$\set{\alpha^n(x)}{n \in \Nset_0}$. 
%%%%%%%%%%%%%%%%%%%%%%%%%%%%%%%%%%%%%%%%%%%%%%%%%%%%%%%%%%%%%%%%%%%%%
\end{proof}
%%%%%%%%%%%%%%%%%%%%%%%%%%%%%%%%%%%%%%%%%%%%%%%%%%%%%%%%%%%%%%%%%%%%%
\begin{Remark}\normalfont \label{rem:mixing-nontrivial}
%%%%%%%%%%%%%%%%%%%%%%%%%%%%%%%%%%%%%%%%%%%%%%%%%%%%%%%%%%%%%%%%%%%%%
The condition $\cN \subset \cM^\alpha$ in Theorem \ref{thm:mixing} is 
non-trivial. Consider a minimal stationary process $\scrM$ with 
$\cN =\cM^\tail = \cM \not \simeq \Cset$. Then $\scrM$ is 
$\cN$-factorizable and $E_{\cN}$ is the identity map on $\cM$. 
Furthermore, $\alpha$ is easily seen to be an automorphism. It follows 
from Definition \ref{def:strong-mixing} that $\alpha$ is strongly mixing 
over $\cN$ if and only if $\alpha$ is the identity.   
%%%%%%%%%%%%%%%%%%%%%%%%%%%%%%%%%%%%%%%%%%%%%%%%%%%%%%%%%%%%%%%%%%%%%
\end{Remark} 
%%%%%%%%%%%%%%%%%%%%%%%%%%%%%%%%%%%%%%%%%%%%%%%%%%%%%%%%%%%%%%%%%%%%%
\begin{Remark}\normalfont
%%%%%%%%%%%%%%%%%%%%%%%%%%%%%%%%%%%%%%%%%%%%%%%%%%%%%%%%%%%%%%%%%%%%%
Conditional order factorizability (CF$_{\text{o}}$) is the weakest 
form of independence or factorizability introduced in Definition 
\ref{def:rvs}; thus Theorem \ref{thm:kol-0-1} and Theorem 
\ref{thm:mixing} are also valid if  (CF$_{\text{o}}$) is replaced 
by (CF), (CI$_{\text{o}}$) or (CI). 
%%%%%%%%%%%%%%%%%%%%%%%%%%%%%%%%%%%%%%%%%%%%%%%%%%%%%%%%%%%%%%%%%%%%%
\end{Remark}
%%%%%%%%%%%%%%%%%%%%%%%%%%%%%%%%%%%%%%%%%%%%%%%%%%%%%%%%%%%%%%%%%%%%%
An important class of stationary processes in classical probability
are Bernoulli shifts; and a noncommutative notion of such shifts emerges 
in \cite{Kuem88b} from the study of stationary quantum Markov processes. 
Here we are interested in their amalgamated version, as studied in 
\cite{Rupp95a} and, in a bilateral continuous `time' formulation, 
in \cite{HKK04a}.
%%%%%%%%%%%%%%%%%%%%%%%%%%%%%%%%%%%%%%%%%%%%%%%%%%%%%%%%%%%%%%%%%%%%%
\begin{Definition} \normalfont \label{def:bernoulli-shift}
%%%%%%%%%%%%%%%%%%%%%%%%%%%%%%%%%%%%%%%%%%%%%%%%%%%%%%%%%%%%%%%%%%%%%
An \emph{(ordered/full) Bernoulli shift (over $\cN$)} is a minimal 
stationary process $\scrB = (\cB,\chi, \beta, \cB_0)$ with the 
following properties: 
\begin{enumerate}
\item $\cN \subset \cB^\alpha \cap \cB_0$ is a $\chi$-conditioned von 
Neumann subalgebra;
\item the canonical filtration $(\cB_I)_{I \subset \Nset_0}$ is 
(order/full) $\cN$-independent;  
\end{enumerate}
The endomorphism $\beta$ is also called a \emph{Bernoulli shift over 
$\cN$ with generator $\cB_0$}.
\end{Definition}
%%%%%%%%%%%%%%%%%%%%%%%%%%%%%%%%%%%%%%%%%%%%%%%%%%%%%%%%%%%%%%%%%%%%%
Note that this definition of a Bernoulli shift contains a subtle
redundancy: one could drop the modular condition on the endomorphism
$\beta$ and conclude it from the fact that its ranges $\beta^n(\cB_0)$
must be $\chi$-conditioned, as required by our definition of independence.
This entails that $\beta$ commutes with $\sigma_t^\chi$, the modular
automorphism group of $(\cB,\chi)$.  
%%%%%%%%%%%%%%%%%%%%%%%%%%%%%%%%%%%%%%%%%%%%%%%%%%%%%%%%%%%%%%%%%%%%%
\begin{Corollary}\label{cor:bernoulli}
%%%%%%%%%%%%%%%%%%%%%%%%%%%%%%%%%%%%%%%%%%%%%%%%%%%%%%%%%%%%%%%%%%%%%
Let $\scrM= (\cM,\psi, \alpha, \cM_0)$ be a minimal stationary 
process. Further suppose $\cN \subset \cM^\alpha$ is a 
$\psi$-conditioned von Neumann subalgebra and $\scrB=
(\cM,\psi, \alpha, \cM_0\vee\cN)$. Then the following are equivalent:
\begin{enumerate}
\item[(a)] $\scrM$ is (order/full) $\cN$-factorizable;
\item[(b)] $\scrM$ is (order/full) $\cN$-independent;
\item[(c)] $\scrB$ is an (ordered/full) Bernoulli shift over $\cN$.  
\end{enumerate}
In particular, it holds $\cN = \cM^\alpha = \cM^\tail$. 
\end{Corollary}
%%%%%%%%%%%%%%%%%%%%%%%%%%%%%%%%%%%%%%%%%%%%%%%%%%%%%%%%%%%%%%%%%%%%%
\begin{proof}
We already know the equivalence of (a) and (b) from Theorem 
\ref{thm:fact=indep}. The equivalence of (b) and (c) is also clear 
since the family $(\cM_I)_{I\subset \Nset_0}$ is (order/full) 
$\cN$-independent if and only if the family 
$(\cM_I \vee \cN)_{I\subset \Nset_0}$ is so. We are left to show 
$\cN= \cM^\alpha=\cM^\tail$. But this is content of Theorem 
\ref{thm:mixing}.
\end{proof}
%%%%%%%%%%%%%%%%%%%%%%%%%%%%%%%%%%%%%%%%%%%%%%%%%%%%%%%%%%%%%%%%%%%%%
We provide next a result which is useful for applications where one 
wants to identify a given process as a Bernoulli shift. Suppose 
$\scrM=(\cM,\psi, \alpha, \cM_0)$ is an (order/full) $\cN$-factorizable 
minimal stationary process for some $\psi$-conditioned von Neumann subalgebra
$\cN \subset \cM^\alpha$. Furthermore let $\cC_0$ be a $\psi$-conditioned 
von Neumann subalgebra of $\cM_0$. Put 
\begin{eqnarray*}
\cB := \bigvee_{n \ge 0} \alpha^n(\cC_0 \vee \cN),\qquad 
\chi:= \psi|_\cB, \qquad 
\beta := \alpha|_{\cB}, \qquad 
\cB_0 := \cC_0 \vee \cN. 
\end{eqnarray*}
This defines the minimal stationary process $\scrB = (\cB, \chi, \beta, \cB_0)$
which is subject of the next result. 
 %%%%%%%%%%%%%%%%%%%%%%%%%%%%%%%%%%%%%%%%%%%%%%%%%%%%%%%%%%%%%%%%%%%%%%%%%%%%%%%
\begin{Corollary}\label{cor:indep-bernoulli} 
%%%%%%%%%%%%%%%%%%%%%%%%%%%%%%%%%%%%%%%%%%%%%%%%%%%%%%%%%%%%%%%%%%%%%%%%%%%%%%%
$\scrB$ is an (ordered/full) Bernoulli shift over $\cN$ and $\cN=\cB^\beta= \cB^\tail$.   
\end{Corollary}
%%%%%%%%%%%%%%%%%%%%%%%%%%%%%%%%%%%%%%%%%%%%%%%%%%%%%%%%%%%%%%%%%%%%%%%%%%%%%%%
\begin{proof}
Theorem \ref{thm:fact=indep} implies the (order/full) $\cN$-independence of $\scrM$. 
Since $\cB_0 \subset \cM_0 \vee \cN$, (order/full) $\cN$-independence is inherited by 
the minimal stationary process $\scrB$. Now an application of Theorem 
\ref{thm:kol-0-1} to the random sequence associated to $\scrB$ ensures 
$\cN= \cB^\tail$. We are left to prove $\cN \subset \cB_0 \cap \cB^\beta$. 
Clearly $\cN \subset \cB_0$. Thus it is suffices to show $\cN = \cB^\beta$. 
Since $\cN = \cM^\alpha$ by Theorem \ref{thm:mixing} and $\cN \subset \cB_0$, 
we have $\cM^\alpha \subset \cB_0$ and consequently $\cM^\alpha \subset \cB$.
But this implies $\cB^\beta= \cM^\alpha$ and consequently $\cN=\cB^\beta$.    
\end{proof}
%%%%%%%%%%%%%%%%%%%%%%%%%%%%%%%%%%%%%%%%%%%%%%%%%%%%%%%%%%%%%%%%%%%%%
\begin{Remark}\normalfont
%%%%%%%%%%%%%%%%%%%%%%%%%%%%%%%%%%%%%%%%%%%%%%%%%%%%%%%%%%%%%%%%%%%%%
Our notion of a Bernoulli shift is motivated from K\"ummerer's work
on noncommutative stationary Markov processes in \cite{Kuem85a,Kuem88a,
Kuem88b,Kuem93UP,Kuem96a}. An ordered Bernoulli shift here is the 
unilateral discrete version of noncommutative continuous Bernoulli shifts 
introduced in \cite{HKK04a}. Note that Definition \ref{def:bernoulli-shift} 
of a Bernoulli shift is \emph{not} restricted to tensor independence; it is 
casted in the broader context of conditional independence.  
\end{Remark}
%%%%%%%%%%%%%%%%%%%%%%%%%%%%%%%%%%%%%%%%%%%%%%%%%%%%%%%%%%%%%%%%%%%%%%%%%%%%%%%%%%%
%%%%%%%%%%%%%%%%%%%%%%%%%%%%%%%%%%%%%%%%%%%%%%%%%%%%%%%%%%%%%%%%%%%%%%%%%%%%%%%%%%%
\section{Spreadability implies conditional order independence}
\label{section:spread2cond-indep}
%%%%%%%%%%%%%%%%%%%%%%%%%%%%%%%%%%%%%%%%%%%%%%%%%%%%%%%%%%%%%%%%%%%%%%%%%%%%%%%%%%%
%%%%%%%%%%%%%%%%%%%%%%%%%%%%%%%%%%%%%%%%%%%%%%%%%%%%%%%%%%%%%%%%%%%%%%%%%%%%%%%%%%%
The main result of this section is Theorem \ref{thm:spread2order-ind} which
is an integral part of the noncommutative extended de Finetti theorem, Theorem
\ref{thm:main}.    
%%%%%%%%%%%%%%%%%%%%%%%%%%%%%%%%%%%%%%%%%%%%%%%%%%%%%%%%%%%%%%%%%%%%%
\begin{Theorem}\label{thm:spread2order-ind}
%%%%%%%%%%%%%%%%%%%%%%%%%%%%%%%%%%%%%%%%%%%%%%%%%%%%%%%%%%%%%%%%%%%%%
A spreadable random sequence $\scrI$ is stationary and 
order $\cM^\tail$-independent.
\end{Theorem}
%%%%%%%%%%%%%%%%%%%%%%%%%%%%%%%%%%%%%%%%%%%%%%%%%%%%%%%%%%%%%%%%%%%%%
It is immediate from Definition \ref{def:ds} that spreadability implies
the stationarity of a random sequence. Thus we can reformulate 
Theorem \ref{thm:spread2order-ind} in terms of stationary processes, as
done in Theorem \ref{thm:spread2bernoulli}. Throughout this section, we 
consider the minimal stationary process
\[
\scrM \equiv (\cM,\psi, \alpha, \cM_0)
\]
and, replacing its generator $\cM_0$ by $\cM_0 \vee \cM^\alpha$, 
the minimal stationary process 
\[
\scrB \equiv (\cM,\psi, \alpha, \cM_0\vee \cM^\alpha).
\]
%%%%%%%%%%%%%%%%%%%%%%%%%%%%%%%%%%%%%%%%%%%%%%%%%%%%%%%%%%%%%%%%%%%%%
\begin{Theorem}\label{thm:spread2bernoulli}
%%%%%%%%%%%%%%%%%%%%%%%%%%%%%%%%%%%%%%%%%%%%%%%%%%%%%%%%%%%%%%%%%%%%%
Suppose $\scrM$ is spreadable and minimal. Then $\scrM$ is order 
$\cM^\tail$-independent and $ \cM^\tail = \cM^\alpha$. 
In particular, $\scrB$ is an ordered Bernoulli shift. 
\end{Theorem}
%%%%%%%%%%%%%%%%%%%%%%%%%%%%%%%%%%%%%%%%%%%%%%%%%%%%%%%%%%%%%%%%%%%%% 
The proof of Theorem \ref{thm:spread2bernoulli} needs some preparation 
and is postponed to the end of this section. It entails of course
the proofs of Theorem \ref{thm:spread2order-ind} and 
Theorem \ref{thm:main} (b)$\,\Rightarrow$(c$_\text{o}$) through the 
correspondence stated in Lemma \ref{lem:corres-rv-endo}.
%%%%%%%%%%%%%%%%%%%%%%%%%%%%%%%%%%%%%%%%%%%%%%%%%%%%%%%%%%%%%%%%%%%%%
\begin{Proposition}\label{prop:mixing-1}
%%%%%%%%%%%%%%%%%%%%%%%%%%%%%%%%%%%%%%%%%%%%%%%%%%%%%%%%%%%%%%%%%%%%%
Suppose the minimal stationary process $\scrM$ is spreadable. Then 
there exists the $\psi$-preserving conditional expectation  
$E_{\cM^\tail}$ of $\cM$ onto $\cM^\tail$ and  
$$
\wotlim_n \alpha^n(x) = E_{\cM^\tail}(x), \qquad x \in \cM. 
$$ 
Moreover, we have
$
\cM^\tail = \cM^\alpha.
$
\end{Proposition}
%%%%%%%%%%%%%%%%%%%%%%%%%%%%%%%%%%%%%%%%%%%%%%%%%%%%%%%%%%%%%%%%%%%%%
\begin{proof}
Let  $\cM_I:= \bigvee_{n \in I} \alpha^n(\cM_0)$ for $I \subset 
\Nset_0$. Let $x,y \in \bigcup_{|I|<\infty} \cM_{I}$. Consequently 
we can assume $x\in \cM_I$ and $y \in \cM_J$ such that there exists 
$N \in \Nset$ with $I \cap (J +N)=\emptyset$. We infer from 
spreadability that $\psi(y \alpha^n(x)) = \psi(y \alpha^{n+1}(x))$ 
for all $n \ge N$. Due to minimality this establishes the limit 
\[
\lim_{n\to \infty} \psi(y \alpha^n(x)) 
\]
on the $\wot$-dense *-algebra $\bigcup_{|I|<\infty}\cM_{I}$. 
A standard approximation argument ensures now the existence of this 
limit for $x,y \in \cM$, using the norm density of the functionals 
$\set{\psi(y\cdot)}{y \in \cM}$ and the boundedness of the set 
$\set{\alpha^n(x)}{n \in \Nset}$. We conclude from this that the 
pointwise $\wot$-limit of the sequence $(\alpha^n)_n$ defines a 
linear map $Q \colon \cM \to \cM$ such that $Q(\cM) \subset \cM^\tail$.

It is easily seen that the linear map $Q$ enjoys
\[
\psi = \psi \circ Q\quad
\text{and} \quad \|Q(x)\| \le \|x\| \text{ for $x\in \cM$.}
\] 
Thus $Q$ is a conditional expectation from $\cM$ onto $\cM^\tail$, 
if we can insure that $Q(x) = x$ for all $x \in \cM^\tail$. To this 
end let $x \in \cM^\tail$ and $y \in \bigcup_{|I|<\infty} \cM_I$. 
We infer from $\cM^\tail \subset \alpha^N (\cM)$ and $\cM_{[N,\infty)} 
\subset \alpha^N(\cM)$  for all $N \in \Nset$ that there exists some 
$N \in \Nset$ such that $x \in \alpha^N(\cM)$ and $y \in \cM_{[0,N-1]}$. 
We approximate $x \in \cM$ in the $\wot$-sense by a sequence 
$(x_k)_k \subset \bigcup_{|I|<\infty}\alpha^{N}(\cM_{I})$  
and conclude further from the definition of $Q$ and from spreadability 
that 
\begin{eqnarray*}
\psi(y Q(x)) 
&=& \lim_k \psi(y Q(x_k)) 
= \lim_k \lim_n \psi(y \alpha^n(x_k))\\
&=& \lim_k \psi(y x_k) = \psi(y x).
\end{eqnarray*}
This shows that $Q(x)= x$ for all $x \in \cM^{\tail}$. Thus $Q$ is 
the conditional expectation of $\cM$ onto $\cM^\tail$ with respect 
to $\psi$ (see \cite[Chapter IX, Definition 4.1]{Take03b}), which we 
denote from now on by $E_{\cM^\tail}$.  

We need to identify the tail algebra as the fixed point algebra.
Proposition \ref{prop:mixing-1} gives pointwise 
$E_{\cM^\tail} E_{\cM^\alpha} = \wotlim_n \alpha^n E_{\cM^\alpha}= 
E_{\cM^{\alpha}}$ and thus $\cM^\alpha \subset \cM^\tail$. The 
inclusion $\cM^\tail \subset \cM^\alpha$ follows from 
$\alpha E_{\cM^\tail} = \lim_n \alpha \,\alpha^n = E_{\cM^\tail}$
in the pointwise $\wot$-sense. 
\end{proof}
%%%%%%%%%%%%%%%%%%%%%%%%%%%%%%%%%%%%%%%%%%%%%%%%%%%%%%%%%%%%%%%%%%%%% 
\begin{Remark}\normalfont
%%%%%%%%%%%%%%%%%%%%%%%%%%%%%%%%%%%%%%%%%%%%%%%%%%%%%%%%%%%%%%%%%%%%%
The proof of Proposition \ref{prop:mixing-1} shows that the 
$\psi$-preserving conditional expectation onto the tail algebra 
$\cM^\tail$ and the fixed point algebra $\cM^\alpha$ of the 
endomorphism $\alpha$ exist under weaker assertions. One does not 
need that $\alpha$ and the modular automorphism group $\sigma_t^\psi$ 
commute (this compatibility condition is required in Definition 
\ref{def:stat-process}).  
\end{Remark}
It is convenient to use Speicher's notion of multilinear maps also for
the endomorphism $\alpha$. we put
\[
\alpha[\ii; \aaa]
:= \alpha^{\ii(1)}(a_1) \alpha^{\ii(2)}(a_2) \cdots\alpha^{\ii(n)}(a_n)
\]
for $n$-tuples $\ii\colon [n] \to \Nset_0$ and 
$\aaa=(a_1, a_2,\ldots, a_n) \in \cM_0$.
%%%%%%%%%%%%%%%%%%%%%%%%%%%%%%%%%%%%%%%%%%%%%%%%%%%%%%%%%%%%%%%%%%%%%
\begin{Definition}\normalfont\label{def:cond-spread}
%%%%%%%%%%%%%%%%%%%%%%%%%%%%%%%%%%%%%%%%%%%%%%%%%%%%%%%%%%%%%%%%%%%%%
A stationary process $\scrM = (\cM,\psi, \alpha, \cM_0)$ or its 
endomorphism $\alpha$ is  \emph{$\cN$-spreadable} if there exists a 
$\psi$-conditioned von Neumann subalgebra $\cN$ of $\cM$ such that 
$$
E_{\cN} (\alpha[\ii; \mathbf{a}]) =E_{\cN} (\alpha[\jj; \mathbf{a}])
$$
for any $n\in \Nset$, $\ii,\jj\colon [n] \to \Nset_0$ with  
$\ii \sim_o \jj $ and $\mathbf{a}\in \cM_0^n$. 
\end{Definition}
%%%%%%%%%%%%%%%%%%%%%%%%%%%%%%%%%%%%%%%%%%%%%%%%%%%%%%%%%%%%%%%%%%%%%
\begin{Lemma}\label{lem:cond-spread}
%%%%%%%%%%%%%%%%%%%%%%%%%%%%%%%%%%%%%%%%%%%%%%%%%%%%%%%%%%%%%%%%%%%%% 
The following are equivalent for a minimal stationary process $\scrM$:
\begin{enumerate}
\item[(a)]
$\mathscr{M}$ is spreadable;
\item[(b)]
$\mathscr{M}$  is $\cM^\tail$-spreadable;
\item[(c)]
$\mathscr{M}$  is $\cM^\alpha$-spreadable.
\end{enumerate}
\end{Lemma}
%%%%%%%%%%%%%%%%%%%%%%%%%%%%%%%%%%%%%%%%%%%%%%%%%%%%%%%%%%%%%%%%%%%%%%
\begin{proof}
(b) and (c) are equivalent since $\cM^\tail=\cM^\alpha$ by 
Proposition \ref{prop:mixing-1}. 
Obviously (b) implies (a) and we are left to prove the converse.
Let us first treat the case $\cM^\tail \subset \cM_0$. We already 
know $\cM^\tail = \cM^\alpha$ from Proposition \ref{prop:mixing-1}. 
Consider the $n$-tuple $(ax_1,x_2,\ldots, x_{n}) \in \cM_0^{n}$ 
with $a\in \cM^\alpha$. We conclude from this that, 
for $\ii,\jj\colon [n]\to \Nset_0$ with $\ii\sim_o \jj$,  
\begin{eqnarray*}
\psi\big(a \alpha[\ii; x_1,x_2,\ldots, x_{n} ]\big)
&=&\psi\big(\alpha[\ii; ax_1,x_2,\ldots, x_{n} ]\big)
= \psi\big(\alpha[\jj; ax_1,x_2,\ldots, x_{n} ]\big)\\
&=&\psi\big(a \alpha[\jj; x_1,x_2,\ldots, x_{n} ]\big).
\end{eqnarray*} 
Using $\psi = \psi \circ E_{\cM^\tail}$ and the module property of 
$E_{\cM^\tail}$, we conclude that $\alpha$ is conditionally 
$\cM^\tail$-spreadable by standard arguments.  

The more general case $\cM^\tail \not \subset \cM_0$ is treated 
similar. We approximate $a \in \cM^\tail$ by a sequence 
$(a_k)_{k \ge 0} \subset \cM$ such that
\begin{eqnarray*}
a_k \in \bigcup_{l \ge k} \alpha^l(\cM_0) \qquad \text{and} \qquad 
a = \sotlim_{k \to \infty} a_k. 
\end{eqnarray*} 
Thus we can assume that each $a_k$ is a linear combination of 
monomials  $\alpha[\ii_k; \mathbf{a}_k ]$, for some $n_k$-tuple 
$\ii_k \colon [n_k] \to \{k, k+1, \ldots\}$ and $\mathbf{a} 
\in \cM_0^{n_k}$.  Now we compute as before that, for 
$\ii,\jj\colon [n]\to \Nset_0$ with $\ii\sim_o \jj$ and sufficiently 
large $k$,   
\begin{eqnarray*}
\psi\big(\alpha[\ii_k; \mathbf{a}_k ]\alpha[\ii; x_1,x_2,\ldots, x_{n} ]\big)
&=&\psi\big( \alpha[\ii_k; \mathbf{a}_k ]\alpha[\jj; x_1,x_2,\ldots, x_{n} ]\big).
\end{eqnarray*} 
This equality extends by linearity and weak* density arguments to 
\begin{eqnarray*}
\psi\big(a \alpha[\ii; x_1,x_2,\ldots, x_{n} ]\big)
&=&\psi\big( a \alpha[\jj; x_1,x_2,\ldots, x_{n} ]\big)
\end{eqnarray*}   
for every $a \in \cM^\tail$. We conclude from this the 
$\cM^\tail$-spreadability of the stationary process. 
\end{proof}
%%%%%%%%%%%%%%%%%%%%%%%%%%%%%%%%%%%%%%%%%%%%%%%%%%%%%%%%%%%%%%%%%%%%%%
\begin{Lemma}\label{lem:spread-fac}
%%%%%%%%%%%%%%%%%%%%%%%%%%%%%%%%%%%%%%%%%%%%%%%%%%%%%%%%%%%%%%%%%%%%%%
Suppose $\scrM$ be a minimal stationary process. If $\scrM$ is spreadable,
then $\scrM$ is order $\cM^\tail$-factorizable.
\end{Lemma}  
%%%%%%%%%%%%%%%%%%%%%%%%%%%%%%%%%%%%%%%%%%%%%%%%%%%%%%%%%%%%%%%%%%%%%%
\begin{proof}
We need to show that the canonical filtration 
$(\cM_I)_{I \subset \Nset_0}$ satisfies the factorization rule
$$
E_{\cM^\tail}(xy) = E_{\cM^\tail}(x) E_{\cM^\tail}(y) 
$$
for all $x\in \cM_I$ and $y \in \cM_J$ whenever $I < J$ or $I > J$. 
Let $x \in \cM_I^\alg$ and $y\in \cM_J^\alg$. Then, for all 
$n \in \Nset_0$,
\[
E_{\cM^\tail}(x y) = E_{\cM^\tail}(x \alpha^n(y)), 
\]
since spreadability implies $\cM^\tail$-spreadability 
(Lemma \ref{lem:cond-spread}). We use the mixing properties 
of $\alpha$ (Proposition \ref{prop:mixing-1}) to conclude
\[
E_{\cM^\tail}(x y) 
= \wotlim_{n \to \infty} E_{\cM^\tail}(x \alpha^n(y)) 
= E_{\cM^\tail}(x) E_{\cM^\tail}(y).
\]
This establishes the order $\cM^\tail$-factorizability of a 
spreadable stationary process.
\end{proof}
%%%%%%%%%%%%%%%%%%%%%%%%%%%%%%%%%%%%%%%%%%%%%%%%%%%%%%%%%%%%%%%%%%%%%
\begin{proof}[Proof of Theorem \ref{thm:spread2bernoulli}]
Lemma \ref{lem:spread-fac} shows that $\scrM$ is order 
$\cM^\tail$-factorizable and Proposition \ref{prop:mixing-1} 
insures $\cM^\tail= \cM^\alpha$. Thus Theorem \ref{thm:fact=indep}
applies for $\cN= \cM^\tail$ and ensures that $\cM$ is 
conditionally $\cM^\tail$-independent. Finally, Corollary 
\ref{cor:bernoulli} entails that $\scrB$ is an ordered Bernoulli 
shift over $\cM^\tail$. 
\end{proof}
%%%%%%%%%%%%%%%%%%%%%%%%%%%%%%%%%%%%%%%%%%%%%%%%%%%%%%%%%%%%%%%%%%%%%
%%%%%%%%%%%%%%%%%%%%%%%%%%%%%%%%%%%%%%%%%%%%%%%%%%%%%%%%%%%%%%%%%%%%%
%%%%%%%%%%%%%%%%%%%%%%%%%%%%%%%%%%%%%%%%%%%%%%%%%%%%%%%%%%%%%%%%%%%%%
\section{Spreadability implies conditional full independence}
\label{section:spread2indep}
%%%%%%%%%%%%%%%%%%%%%%%%%%%%%%%%%%%%%%%%%%%%%%%%%%%%%%%%%%%%%%%%%%%%%
%%%%%%%%%%%%%%%%%%%%%%%%%%%%%%%%%%%%%%%%%%%%%%%%%%%%%%%%%%%%%%%%%%%%%
We have already shown in the previous section that spreadability
implies conditional order independence. Here this result will be
strengthened to conditional full independence.  
%%%%%%%%%%%%%%%%%%%%%%%%%%%%%%%%%%%%%%%%%%%%%%%%%%%%%%%%%%%%%%%%%%%%%
\begin{Theorem}\label{thm:spread2ind}
%%%%%%%%%%%%%%%%%%%%%%%%%%%%%%%%%%%%%%%%%%%%%%%%%%%%%%%%%%%%%%%%%%%%%
A spreadable random sequence $\scrI$ is stationary and 
full $\cM^\tail$-independent.
\end{Theorem}
%%%%%%%%%%%%%%%%%%%%%%%%%%%%%%%%%%%%%%%%%%%%%%%%%%%%%%%%%%%%%%%%%%%%%
Theorem \ref{thm:spread2ind} establishes the implication 
(b) $\Rightarrow$ (c) of Theorem \ref{thm:main}, the noncommutative
extended de Finetti theorem. We will prove it in terms of the
corresponding stationary process $\scrM = (\cM,\psi, \alpha, \cM_0)$
and, replacing the generator $\cM_0$ by $\cM_0 \vee \cM^\alpha$, denote by 
$\scrB$ the stationary process $(\cM,\psi, \alpha, \cM_0 \vee \cM^\alpha)$.
%%%%%%%%%%%%%%%%%%%%%%%%%%%%%%%%%%%%%%%%%%%%%%%%%%%%%%%%%%%%%%%%%%%%%
\begin{Theorem}\label{thm:spread2fullbernoulli}
%%%%%%%%%%%%%%%%%%%%%%%%%%%%%%%%%%%%%%%%%%%%%%%%%%%%%%%%%%%%%%%%%%%%%
Suppose $\scrM$ is spreadable and minimal. Then $\scrM$ is  
full $\cM^\tail$-independent and $ \cM^\tail = \cM^\alpha$. 
In particular, $\scrB$ is a full Bernoulli shift. 
\end{Theorem}
%%%%%%%%%%%%%%%%%%%%%%%%%%%%%%%%%%%%%%%%%%%%%%%%%%%%%%%%%%%%%%%%%%%%%
The proofs of Theorem \ref{thm:spread2ind} and Theorem 
\ref{thm:spread2fullbernoulli} require a certain refined version of 
the mean ergodic theorem. Let us start with its usual formulation 
and include for the convenience of the reader how its proof reduces 
to the usual result for contractions on Hilbert spaces. 
%%%%%%%%%%%%%%%%%%%%%%%%%%%%%%%%%%%%%%%%%%%%%%%%%%%%%%%%%%%%%%%%%%%%%%
\begin{Theorem}\label{thm:mean-ergodic}
%%%%%%%%%%%%%%%%%%%%%%%%%%%%%%%%%%%%%%%%%%%%%%%%%%%%%%%%%%%%%%%%%%%%%%
Let $(\scrM,\psi)$ be a probability space and $\alpha$ a $\psi$-preserving
endomorphism of $\cM$. Then we have, for each $x\in \cM$,
\[
\sotlim_{n\to \infty} \frac{1}{n}\sum_{k=0}^{n-1} \alpha^k(x) 
= E_{\cM^\alpha}(x).   
\]
\end{Theorem}
%%%%%%%%%%%%%%%%%%%%%%%%%%%%%%%%%%%%%%%%%%%%%%%%%%%%%%%%%%%%%%%%%%%%%%
\begin{proof}
The strong operator topology and the $\psi$-topology generated by the
maps $x \mapsto \psi(x^*x)^{1/2}$, $x\in \cM$, coincide on norm bounded 
sets in $\cM$. Thus this mean ergodic theorem is an immediate 
consequence of the usual mean ergodic theorem in Hilbert spaces 
(see \cite[Theorem 1.2]{Pete83a} for example). 
\end{proof}
%%%%%%%%%%%%%%%%%%%%%%%%%%%%%%%%%%%%%%%%%%%%%%%%%%%%%%%%%%%%%%%%%%%%%
This mean ergodic theorem would allow us to give an alternative proof
of that spreadability implies conditional order independence (CI$_\text{o}$), 
after having identified the tail algebra as the fixed point algebra of 
the stationary process in Proposition \ref{prop:mixing-1} and 
established conditional spreadability in Lemma \ref{lem:cond-spread}.  

We illustrate this by an example. Given the stationary process 
$(\cM,\psi, \alpha, \cM_0)$, let $a,b\in \cM_0$ and consider 
\begin{eqnarray*}
\cM_{\{1,2\}} \ni x &=& \alpha(a) \alpha^2(a) \alpha(a) \alpha^2(a) ,\\
\cM_{\{3,4\}} \ni y &=& \alpha^4(b) \alpha^3(b) \alpha^4(b) \alpha^3(b) \alpha^4(b).  
\end{eqnarray*}
We have $\{1,2\} < \{3,4\}$ and thus spreadability implies 
$$
E_{\cM^\alpha}(xy) 
=  E_{\cM^\alpha}\big(x\alpha^n(y)\big) 
= E_{\cM^\alpha}\Big(x\frac{1}{n}\sum_{k=0}^{n-1}\alpha^k (y)\Big)   
$$
for all $n \ge 1$. Thus Theorem \ref{thm:mean-ergodic} implies 
$E_{\cM^\alpha}(xy) = E_{\cM^\alpha}(x) E_{\cM^\alpha}(y)$. 

But such an argument falls short of establishing the apparently stronger 
version, conditional full independence (CI). For example, consider the 
two elements
\begin{eqnarray*}
x &=& \alpha(a) \alpha^3(a) \alpha(a) \alpha^3(a),\\
y &=& \alpha^4(b) \alpha^2(b) \alpha^4(b) \alpha^2(b) \alpha^4(b),  
\end{eqnarray*}    
Thus we have $x \in \cM_I$ and $y \in \cM_J$
with $I = \{1,3\}$ and $J\in \{2,4\}$. Since the tuples $(1,3,1,3,4,2,4,2,4)$
and $(1,3,1,3,4+n,2+n,4+n,2+n,4+n)$ are order equivalent if and only if $n=0$,
the previous arguments fails. We observe that spreadability implies, in particular,
\begin{eqnarray*}
E_{\cM^\alpha}(xy) 
&=& E_{\cM^\alpha}\big(x
   \alpha^{4+n}(b) \alpha^2(b) \alpha^{4+n}(b) \alpha^2(b) \alpha^{4+n}(b) \big)\\
&=& E_{\cM^\alpha}\Big(x \frac{1}{n}\sum_{k=0}^{n-1} 
    \alpha^{4+k}(b) \alpha^2(b) \alpha^{4+k}(b) \alpha^2(b) \alpha^{4+k}(b) \Big).
\end{eqnarray*}
but a direct application of the mean ergodic theorem is still out of reach.  

To overcome such difficulties we need to provide a more elaborated version of Theorem 
\ref{thm:mean-ergodic} which allows us to preserve relative localisation properties 
of the canonical filtration $(\cM_I)_{I}$ while performing mean ergodic averages. Since 
this result is of interest in its own, we formulate it in greater generality as 
necessary for our purposes.
%%%%%%%%%%%%%%%%%%%%%%%%%%%%%%%%%%%%%%%%%%%%%%%%%%%%%%%%%%%%%%%%%%%%%
\begin{Theorem} \label{thm:local-mean-ergodic}
%%%%%%%%%%%%%%%%%%%%%%%%%%%%%%%%%%%%%%%%%%%%%%%%%%%%%%%%%%%%%%%%%%%%%
Let $(\cM,\psi)$ be a probability space and suppose 
$\{\alpha_N\}_{N \in \Nset_0}$ is a family of $\psi$-preserving
completely positive linear maps of $\cM$ satisfying
\begin{enumerate}
\item
$\cM^{\alpha_N} \subset \cM^{\alpha_{N+1}}$ for all $N \in \Nset_0$;
\item 
$\cM =  \bigvee_{N \in \Nset_0} \cM^{\alpha_N}$.
\end{enumerate}
Furthermore let
\[
M_N^{(n)} := \frac{1}{n}\sum_{k=0}^{n-1} \alpha_N^k
\quad \text{and} \quad 
T_N := \vec{\prod}_{l=0}^{N} \alpha_{l}^{lN} M_{l}^{(N)}.
\]
Then we have
\[
\sotlim_{N \to \infty} T_N(x) = E_{\cM^{\alpha_0}}(x)
\] 
for any $x \in \cM$. 
\end{Theorem} 
%%%%%%%%%%%%%%%%%%%%%%%%%%%%%%%%%%%%%%%%%%%%%%%%%%%%%%%%%%%%%%%%%%%%%
\begin{proof}
Since the family $\set{T_N}{N \in \Nset_0}$ is bounded, its pointwise 
$\sot$-convergence follows by a standard approximation
argument if we can establish this convergence on the  
weak*-dense *-subalgebra $\bigcup_{N\in \Nset_0} \cM^{\alpha_{N}}$ 
of $\cM$. 

Let $x\in \cM^{\alpha_{N_0}}$ for some $N_0 \in \Nset$ and $N \ge N_0$.
Since $\alpha_N(x) = x$ and thus $M_N^{(n)}(x)=x$,
the ordered product has at most $N_0$ non-trivially acting factors:
\[
T_N(x) = \Big(\vec{\prod}_{l=0}^{N} \alpha_{l}^{lN} M_{l}^{(N)}\Big)(x) = 
\Big(\vec{\prod}_{l=0}^{N_0-1} 
\alpha_{l}^{lN} M_{l}^{(N)}\Big)(x).
\] 
The assertions on the fixed point algebras $\cM^{\alpha_k}$ imply
that, for any $k \le N$ and $n \in \Nset$,
\[
E_{\cM^{\alpha_k}}\alpha_N = E_{\cM^{\alpha_k}} 
\quad \text{and}\quad 
E_{\cM^{\alpha_k}}M_N^{(n)} = E_{\cM^{\alpha_k}}. 
\]  
Thus we can rewrite $T_N(x)$ as a finite telescope sum, assuming $N \ge N_0$:    
\begin{eqnarray*}
T_N(x) 
&=& M_0^{(N)}\alpha_1^{N} M_1^{(N)} \alpha_2^{2N} M_2^{(N)} \cdots \alpha_N^{N^2} M_N^{(N)}(x)\\
&=& M_0^{(N)}\alpha_1^{N} M_1^{(N)} \alpha_2^{2N} M_2^{(N)} \cdots \alpha_{N_0-1}^{(N_0-1)N} M_{N_0-1}^{(N)}(x)\\
&=&   \Big(\vec{\prod}_{l=0}^{N_0-1} 
      \alpha_{l}^{lN} (M_{l}^{(N)}- E_{\cM^{\alpha_l}})\Big)  
      E_{\cM^{\alpha_{N_0}}}(x)\\    
&& +\quad  \Big(\vec{\prod}_{l=0}^{N_0-2} 
      \alpha_{l}^{lN} (M_{l}^{(N)}- E_{\cM^{\alpha_l}})\Big)  
      E_{\cM^{\alpha_{N_0-1}}}(x)\\
&& +\quad  \Big(\vec{\prod}_{l=0}^{N_0-3} 
      \alpha_{l}^{lN} (M_{l}^{(N)}- E_{\cM^{\alpha_l}})\Big)  
      E_{\cM^{\alpha_{N_0-2}}}(x)\\
&&+ \quad \cdots \\
&&+ \quad (M_0^{(N)}-E_{\cM^{\alpha_0}}) 
    \alpha_1^{N} (M_1^{(N)}-E_{\cM^{\alpha_1}})  
    E_{\cM^{\alpha_{2}}}(x)\\ 
&&+ \quad (M_0^{(N)}-E_{\cM^{\alpha_0}})  
    E_{\cM^{\alpha_{1}}}(x)\\
&&+ \quad E_{\cM^{\alpha_0}}(x).           
\end{eqnarray*}
The strong operator topology and the $\psi$-topology generated
by $x \mapsto \|x\|_\psi^2 :=\psi(x^*x)$ coincide on bounded sets of $\cM$.
Thus 
\[
\left\|\Big(\vec{\prod}_{l=0}^{k-1} 
      \alpha_{l}^{lN} (M_{l}^{(N)}- E_{\cM^{\alpha_l}})\Big)  
      E_{\cM^{\alpha_{k}}}(x) \right\|_\psi
\le  2^{k-1} \left\| (M_{k-1}^{(N)}- E_{\cM^{\alpha_{k-1}}})  
      E_{\cM^{\alpha_{k}}}(x) \right\|_\psi,      
\]
and the usual mean ergodic theorem, Theorem \ref{thm:mean-ergodic}, entail that all terms  
of above telescope sum, except $E_{\cM^{\alpha_0}}(x)$, vanish in the limit $N \to \infty$.  
\end{proof}

We will connect this refined mean ergodic theorem to partial shifts which 
canonically emerge from a spreadable endomorphism.  Recall for this purpose 
the notion of partial shifts $\theta_N$ of $\Nset_0$ and their relation 
to order invariance of tuples (see Remark \ref{rem:partial-shift}):
\[
\theta_N(n) = 
\begin{cases}
n & \text{ if $n < N$};\\
n+1 & \text{ if $ n \ge N$}.
\end{cases}
\]
Clearly $\theta_N$ is an order preserving map of $\Nset_0$ into itself and 
so are the compositions of such maps with $N \in \Nset_0$. Here we are 
interested in compositions of the type
\[
\theta_{N, \vec{l}_N}
:= \vec{\prod}_{i=0}^N \theta_i^{iN+l_i} 
=  \theta_0^{l_0}\theta_1^{N+l_1}\theta_2^{2N+l_2} \cdots 
\theta_{N-1}^{(N-1)N+l_{N-1}} \theta_N^{N^2+l_{N}}, 
\]
where 
\[
\vec{l}_N = (l_0, l_1, \ldots, l_N) \in \{0,1,\ldots,N-1\}^{N+1}.
\] 
Note that the $\theta_{i}$'s in the ordered product do not commute 
for different $i$'s. We record two simple, but crucial properties of 
this composition.
\begin{Lemma} \label{lem:theta}
For any $(N+1)$-tuples $\vec{l}_N, \vec{k}_N \in \{0, 1, \ldots, N-1\}^{N+1}$, 
it holds
\begin{align*}
\theta_{N, \vec{l}_N}(i) < \theta_{N, \vec{k}_N}(j)  
& \quad \text{whenever $i <j <N$, and}\\
\theta_{N, \vec{l}_N}(I) \cap \theta_{N, \vec{k}_N}(J)= \emptyset 
&\quad\text{whenever $I \cap J = \emptyset$ and $\max{I \cup J} <N$.}
\end{align*}
\end{Lemma}
\begin{proof} 
Since all $\theta_N$'s are order preserving, it suffices to consider $j=i+1$. One calculates
\[
\theta_{N, \vec{l}_N}(i+1)- \theta_{N, \vec{l}_N}(i) 
= 1 + \sum_{j=0}^{i}(k_j-l_j)  +  \big((i+1)N  +  k_{i+1}\big) > 0.  
\]
Moreover this ensures that the images of disjoint sets $I, J$ (bounded by $N$) are disjoint.   
\end{proof}
Suppose for the remainder of this section that the stationary 
process $\scrM = (\cM, \psi, \alpha, \cM_0)$ is minimal and let, 
for $N \in \Nset$,
\[
\cM_{N-1} := \bigvee_{0\le k < N} \alpha^k(\cM_0).
\] 
Spreadability of $\scrM$ allows us to promote the partial shifts $\theta_N$ of 
$\Nset_0$ to endomorphisms of $\cM$. Let
\[
\alpha[\ii;\aaa] := \alpha^{\ii(1)}(a_1) \alpha^{\ii(2)}(a_2)\cdots \alpha^{\ii(n)}(a_n)
\]
for $n$-tuples $\ii\colon [n] \to \Nset_0$ and 
$\aaa = (a_1,a_2, \ldots, a_n) \in \cM_0^n$. 
%%%%%%%%%%%%%%%%%%%%%%%%%%%%%%%%%%%%%%%%%%%%%%%%%%%%%%%%%%%%%%%%%%%%%
\begin{Lemma} \label{lem:subshifts}
%%%%%%%%%%%%%%%%%%%%%%%%%%%%%%%%%%%%%%%%%%%%%%%%%%%%%%%%%%%%%%%%%%%%%
Suppose the endomorphism $\alpha$ of $\cM$ is spreadable and let 
$N \in \Nset_0$. Then the complex linear extension of the map 
\[
\alpha[\ii;\aaa] \mapsto \alpha[\theta_N\circ \ii;\aaa]
\]
defines a $\psi$-preserving unital endomorphism $\alpha_N$ of $\cM$,
such that 
\[
\cM_{N} \subset \cM^{\alpha_{N+1}}.
\] 
In particular, $\scrM_N:=(\cM, \psi, \alpha_N, \cM_N)$ is a minimal 
stationary process.
\end{Lemma} 
%%%%%%%%%%%%%%%%%%%%%%%%%%%%%%%%%%%%%%%%%%%%%%%%%%%%%%%%%%%%%%%%%%%%%  
\begin{proof}
The map $\alpha_N$ is well-defined on the *-algebra $\cM^\alg_{\Nset_0}$,
the $\Cset$-linear span of monomials $\alpha[\ii;\aaa]$. Indeed, the 
faithfulness of $\psi$ and spreadability give
\begin{eqnarray*}
\sum_k \alpha[\theta_N \circ \ii_k;\aaa_k] = 0 
&\Leftrightarrow& 
\psi\Big(\big|\sum_k \alpha[\theta_N \circ\ii_k;\aaa_k]\big|^2\Big) 
=  
\psi\Big(\big|\sum_k \alpha[\ii_k;\aaa_k]\big|^2\Big)=0\\
&\Leftrightarrow&
\sum_k \alpha[\ii_k;\aaa_k] = 0. 
\end{eqnarray*} 
Thus $\alpha_N$ is well-defined on $\cM^\alg_{\Nset_0}$. Now it is routine 
to check that $\alpha_N$ extends to a $\psi$-preserving unital 
endomorphism of $\cM$, denoted by the same symbol. The inclusion 
$\cM_{N-1}\subset \cM^{\alpha_N}$ is immediately concluded by approximation 
from the definition of $\alpha_N$ on $\cM^\alg_{\Nset_0}$. It is also clear 
that $\alpha_N$ commutes with the modular automorphism group of $(\cM,\psi)$ 
since $\alpha$ does so. Thus $(\cM, \psi, \alpha_N, \cM_N)$ is a stationary
process which is easily seen to be minimal.   
\end{proof}
%%%%%%%%%%%%%%%%%%%%%%%%%%%%%%%%%%%%%%%%%%%%%%%%%%%%%%%%%%%%%%%%%%%%%
\begin{Corollary}\label{cor:subshifts}
%%%%%%%%%%%%%%%%%%%%%%%%%%%%%%%%%%%%%%%%%%%%%%%%%%%%%%%%%%%%%%%%%%%%%
The minimal stationary processes $\scrM_N$ and their endomorphisms 
$\alpha_N$ are spreadable. Moreover, it holds for $N \in \Nset_0$: 
\begin{enumerate}
\item \label{item:subshifts-i}
$\alpha_{N+1}|_{\alpha_N(\cM)} = \alpha_N|_{\alpha_N(\cM)}$; 
\item \label{item:subshifts-ii}
$\cM^{\alpha_N} \subset \cM^{\alpha_{N+1}}$; 
\item \label{item:subshifts-iii}
$\cM = \bigvee_{N \in \Nset_0} \cM^{\alpha_N}$.
\end{enumerate}
\end{Corollary}
%%%%%%%%%%%%%%%%%%%%%%%%%%%%%%%%%%%%%%%%%%%%%%%%%%%%%%%%%%%%%%%%%%%%%
\begin{proof}
The spreadability of $\scrM_N$ is immediate from definition of $\alpha_N$ in Lemma 
\ref{lem:subshifts} and the spreadability of $\alpha$. 

(\ref{item:subshifts-i}) Clearly, 
$\theta_{N+1}|_{\theta_N(\Nset_0)} = \theta_{N}|_{\theta_N(\Nset_0)}$. Thus
$\alpha_{N+1}$ and $\alpha_{N}$ coincide on the $\Cset$-linear span
of all monomials of the form $\alpha[\theta_N \circ \ii; \aaa] = \alpha_N \Big(\alpha[\ii;\aaa)]\Big)$.
Now the assertion follows from the weak*-density of this span in $\alpha_N(\cM)$. 
    
(\ref{item:subshifts-ii})
$\cM^{\alpha_N}$ is contained in $\alpha_N(\cM)$. By (\ref{item:subshifts-i}), 
$\alpha_N$ and $\alpha_{N+1}$ coincide on $\alpha_N(\cM)$. 
Thus $\cM^{\alpha_N} \subset \cM^{\alpha_{N+1}}$. 

(\ref{item:subshifts-iii}) This is evident from the minimality of $\scrM$ since 
$\bigvee_{0 \le n < N} \alpha^n(\cM_0) \subset \cM^{\alpha_N}$ by Lemma
\ref{lem:subshifts}. 
\end{proof}
%%%%%%%%%%%%%%%%%%%%%%%%%%%%%%%%%%%%%%%%%%%%%%%%%%%%%%%%%%%%%%%%%
\begin{Remark}\normalfont
We do not know at the time of this writing if the fixed point algebras 
$\cM^{\alpha_N}$ can be identified as 
$
\cM^{\alpha_N} =  \bigvee_{0 \le n < N} \alpha^n(\cM_0) \vee \cM^\alpha. 
$
\end{Remark}
%%%%%%%%%%%%%%%%%%%%%%%%%%%%%%%%%%%%%%%%%%%%%%%%%%%%%%%%%%%%%%%%%
\begin{proof}[Proof of Theorem \ref{thm:spread2fullbernoulli}]
We need to show that 
\[
E_{\cM^\alpha}(xy) = E_{\cM^\alpha}(x) E_{\cM^\alpha}(y)
\]
for all $x\in \cM_I$ and $y\in \cM_J$ with $I \cap J =\emptyset$. We 
start with disjoint finite sets $I$ and $J$, and elements of the form 
$$
x = \alpha[\ii;\aaa] \quad \text{and} \quad y = \alpha[\jj;\bbb],
$$
for $p$-tuples $\ii\colon [p] \to I$, $\aaa \in \cM_0^p$ and $q$-tuples
$\jj\colon [q] \to J$, $\bbb \in \cM_0^q$.

Recall that $\scrM_\alpha$ is $\cM^\alpha$-spreadable by Lemma \ref{lem:cond-spread}
and so
\begin{eqnarray*}
E_{\cM^{\alpha_0}}(xy) 
&=& E_{\cM^\alpha} \Big(\alpha[\ii;\aaa] \,\alpha[\jj;\bbb]\big)\\
&=& E_{\cM^\alpha} \big(\alpha[ \theta_{N,\vec{k}_N} \circ \ii;\aaa]\,  
        \alpha[\theta_{N,\vec{k}_N} \circ \jj,\bbb]\Big)                
\end{eqnarray*}
for any $\vec{k}_N \in \{0,1, \ldots, N-1\}^{N+1}$ and $N > \max I \cup J$. By Lemma
\ref{lem:theta}, the maps $\theta_{N,\vec{k}_N}$ are order preserving on $\Nset_0$ and  
$I \cap J = \emptyset$ implies $\theta_{N,\vec{k}_N}(I) \cap \theta_{N,\vec{l}_N}(J)
= \emptyset$ for any $(N+1)$-tuples $\vec{k}_N, \vec{l}_N \in \{0,1, \ldots, N-1\}^{N+1}$.
Thus
\begin{eqnarray*}
E_{\cM^\alpha} \Big(\alpha[\ii;\aaa] \,\alpha[\jj;\bbb]\Big)
&=& E_{\cM^\alpha} \Big(\alpha[ \theta_{N,\vec{k}_N} \circ \ii;\aaa]\,  
        \alpha[\theta_{N,\vec{l}_N} \circ \jj,\bbb]\Big)                
\end{eqnarray*}    
for all $\vec{k}_N,\vec{l}_N$. Consequently we can pass on the right side of this equation 
to the mean ergodic averages by summing over the variables 
$k_0, k_1, \ldots, k_N$ and $l_0,l_1, \ldots l_N$. Doing so we find
\begin{eqnarray*}
E_{\cM^\alpha} \Big(\alpha[\ii;\aaa] \,\alpha[\jj;\bbb]\big)
&=& E_{\cM^\alpha} \Big(T_N(\alpha[\ii;\aaa])\,  
        T_N(\alpha[\jj;\bbb])\Big)
\end{eqnarray*}
for all $N > \max{I \cup J}$, where
\[
T_N := \vec{\prod}_{l=0}^{N} \alpha_{l}^{lN} M_{l}^{(N)} \quad \text{with } 
M_N^{(n)} := \frac{1}{n}\sum_{k=0}^{n-1} \alpha_N^k.
\]
Since Corollary \ref{cor:subshifts} ensures that
all assumptions of the refined mean ergodic theorem Theorem \ref{thm:local-mean-ergodic} 
are satisfied, the pointwise $\sot$-convergence of  $T_N$ to $E_{\cM^{\alpha_0}} 
(= E_{\cM^\alpha})$ for $N \to \infty$ establishes 
\[
E_{\cM^\alpha} \Big(\alpha[\ii;\aaa] \,\alpha[\jj;\bbb]\big) 
= E_{\cM^\alpha}\big(\alpha[\ii;\aaa]\big)\, E_{\cM^\alpha}\big(\alpha[\jj;\bbb]\big)
\]   
for any $\ii$ and $\jj$ with disjoint ranges. This generalizes to the $\Cset$-linear 
span of monomials $\alpha[\ii_n;\aaa_n]$ and $\alpha[\jj_n;\bbb_n]$, provided the range of 
the tuples $\ii_n$ is contained in $I$ and the range of the tuples $\jj_n$ is contained in $J$.  
Now a density argument establishes the factorization
\[
E_{\cM^{\alpha}}(xy) =  E_{\cM^{\alpha}}(x) E_{\cM^{\alpha}}(y) 
\]
for all $x\in \cM_I$ and $y \in \cM_J$ whenever $I$ and $J$ are finite disjoint subsets of $\Nset_0$.
Finally, another approximation removes the assumption of the finiteness of $I$ and $J$. Thus 
we have established that the spreadability of a minimal stationary process $\scrM$ implies its full 
$\cM^\alpha$-factorizability. 

By Theorem \ref{thm:fact=indep}, full $\cM^\alpha$-factorizability and full $\cM^\alpha$-independence 
are equivalent. In particular, we know already $\cM^\alpha=\cM^\tail$ from Theorem \ref{thm:spread2bernoulli}. 
Finally, Corollary \ref{cor:bernoulli} entails that $\scrB$ is a full Bernoulli shift. 
\end{proof} 

%%%%%%%%%%%%%%%%%%%%%%%%%%%%%%%%%%%%%%%%%%%%%%%%%%%%%%%%%%%%%%%%%%%%%
\begin{Remark}\normalfont
The refined version of the mean ergodic theorem, Theorem \ref{thm:local-mean-ergodic}, 
is motivated in parts from product representations of endomorphisms as their study is 
started in \cite{Gohm04a} and as they are applied to braid group representations in 
\cite{GoKo08a}. Suppose the probability space $(\cM,\psi)$ is equipped with a tower
\[
\cM_0 \subset \cM_1 \subset \cM_2 \subset \cdots
\]
of $\psi$-expected subalgebras such that $\cM = \bigvee_{n \ge 0}\cM_n$ and
consider a family of automorphisms $(\gamma_k)_{k \in \Nset} \subset \Aut{\cM,\psi}$ 
satisfying
\begin{eqnarray*}
\gamma_k(\cM_n) = \cM_n && \text{if $k \le n$}\\
\gamma_k|_{\cM_{n-1}}= \id|_{\cM_{n-1}} && \text{if $k \ge n+1$}. 
\end{eqnarray*}
Then
$$
\alpha_N := \lim_{n \to \infty} \gamma_{N+1}\cdots \gamma_{n} 
$$
exists in the pointwise strong operator topology and defines a family 
of $\psi$-preserving endomorphisms $\{\alpha_N\}_{N \in \Nset_0}$ of 
$\cM$ such that $\cM_N \subset \cM^{\alpha_N} \subset \cM^{\alpha_{N+1}}$ 
for all $N \in \Nset_0$.
   
Suppose now in addition that 
\[
\alpha_N|_{\alpha_0^k(\cM)} = \alpha_0|_{\alpha_0^k(\cM)} \quad\text{if $k \ge N$}. 
\]
Then it can be seen that 
the refined mean ergodic theorem preserves localization properties with 
respect to the filtration $(\cA_I)_{I \subset \Nset_0}$, where
$\cA_I := \bigvee_{i \in I} \alpha_0^i(\cM_0)$. To be more precise, 
suppose $x\in \cA_I$ and $y \in \cA_J$ with $I \cap J = \emptyset$. Then 
for every $N$, there exist sets $I_N, J_N$ with $I_N \cap J_N = \emptyset$
such that $T_N(x)\in \cA_{I_N}$ and $T_N(y) \in \cA_{J_N}$. Such a feature 
turned out to be crucial for the proof that spreadability implies conditional 
full independence. 
\end{Remark}
%%%%%%%%%%%%%%%%%%%%%%%%%%%%%%%%%%%%%%%%%%%%%%%%%%%%%%%%%%%%%%%%%%%%%
%%%%%%%%%%%%%%%%%%%%%%%%%%%%%%%%%%%%%%%%%%%%%%%%%%%%%%%%%%%%%%%%%%%%%%%%%%%%%%%%%%%
\section{Some applications and outlook}
\label{section:applications}
%%%%%%%%%%%%%%%%%%%%%%%%%%%%%%%%%%%%%%%%%%%%%%%%%%%%%%%%%%%%%%%%%%%%%%%%%%%%%%%%%%%
We briefly address some further developments and applications of Theorem \ref{thm:main}.
%%%%%%%%%%%%%%%%%%%%%%%%%%%%%%%%%%%%%%%%%%%%%%%%%%%%%%%%%%%%%%%%%%%%%%%%%%%
\subsection{A glimpse on braidability}
\label{subsection:braid-groups}
%%%%%%%%%%%%%%%%%%%%%%%%%%%%%%%%%%%%%%%%%%%%%%%%%%%%%%%%%%%%%%%%%%%%%%%%%%% 
The Artin's braid 
group $\Bset_\infty$ is presented by the generators $\sigma_1, \sigma_2,\ldots$, 
subject to the relations 
\begin{align*}
&&\sigma_i \sigma_{j} \sigma_i &= \sigma_{j} \sigma_i \sigma_{j} &\text{if $ \; \mid i-j \mid\, = 1 $},&& \\
&&\sigma_i \sigma_j &= \sigma_j \sigma_i  &\text{if $ \; \mid i-j \mid\, > 1 $}.&& 
\end{align*}
$\Bset_n$ is an important extension of the symmetric group $\Sset_n$ and we 
introduce in \cite{GoKo08a} `braidability' as a notion which extends 
exchangeability.
%%%%%%%%%%%%%%%%%%%%%%%%%%%%%%%%%%%%%%%%%%%%%%%%%%%%%%%%%%%%%%%%%%%%%%%%%%%%%
\begin{Definition}\normalfont
%%%%%%%%%%%%%%%%%%%%%%%%%%%%%%%%%%%%%%%%%%%%%%%%%%%%%%%%%%%%%%%%%%%%%%%%%%%%%
A random sequence $\scrI$ with random variables
\[
\iota \equiv (\iota_n)_{n \ge 0} \colon (\cA_0, \varphi_0) \rightarrow (\cM,\psi)
\] 
is \emph{$\rho$-braidable} if there exists a representation  
$
\rho \colon \Bset_\infty \to \Aut{\cM,\psi}
$ 
satisfying:
\begin{align*}
&&&&\iota_n &= \rho(\sigma_n \sigma_{n-1} \cdots \sigma_1)\iota_0 
&&\text{for all $n \ge 1$};\\
&&&& \iota_0    & = \rho(\sigma_n) \iota_0&& \text{if $n \ge 2$.}   
\end{align*}
\end{Definition}
%%%%%%%%%%%%%%%%%%%%%%%%%%%%%%%%%%%%%%%%%%%%%%%%%%%%%%%%%%%%%%%%%%%%%%%
Note that the representation $\rho$ may be non-faithful and comprises
representations of $\Sset_\infty$. More precisely, it is shown in
\cite{GoKo08a} that the following are equivalent:
\begin{enumerate}
\item $\scrI$ is exchangeable;
\item $\scrI$ is $\rho$-braidable and $\rho(\sigma^2_n)=\id$ for all $n \in \Nset$.   
\end{enumerate}
So exchangeability clearly implies braidability. A main result of \cite{GoKo08a}
is that braidability is intermediate between two distributional symmetries and 
thus provides a refinement of the noncommutative extended de Finetti theorem, 
Theorem \ref{thm:main}: 
%%%%%%%%%%%%%%%%%%%%%%%%%%%%%%%%%%%%%%%%%%%%%%%%%%%%%%%%%%%%%%%%%%%%%%%%
\begin{Theorem}[\cite{GoKo08a}] \label{thm:braid-rep}
%%%%%%%%%%%%%%%%%%%%%%%%%%%%%%%%%%%%%%%%%%%%%%%%%%%%%%%%%%%%%%%%%%%%%%%%
Let $\scrI$ be an infinite random sequence and consider the following
statements: 
\begin{enumerate}
\item[(a)]  $\scrI$ is exchangeable;
\item[(ab)]   $\scrI$ is braidable;
\item[(b)]  $\scrI$ is spreadable;
\item[(c)]  $\scrI$ is stationary and full $\cM^\tail$-independent. 
\end{enumerate}
Then we have the implications:
\begin{center}
(a) $\Rightarrow$ (ab) $\Rightarrow$ (b) $\Rightarrow$  (c) 
\end{center}
%%%%%%%%%%%%%%%%%%%%%%%%%%%%%%%%%%%%%%%%%%%%%%%%%%%%%%%%%%%%%%%%%%%%%%%%
\end{Theorem}
%%%%%%%%%%%%%%%%%%%%%%%%%%%%%%%%%%%%%%%%%%%%%%%%%%%%%%%%%%%%%%%%%%%%%%%%
Starting from braid group representations, this result implies a rich 
structure of triangular arrays of commuting squares, similar as they 
emerge from the Jones fundamental construction in subfactor theory. 
We refer the interested reader to \cite{GoKo08a} for further details 
and developments.

We need another result from \cite{GoKo08a} to complete the proof 
of Theorem \ref{thm:main}.
\begin{Theorem}[\cite{GoKo08a}]
There exist examples of infinite random sequences such that the 
implications (a) $\Leftarrow$ (b)' and `(b) $\Leftarrow$ (c)' 
fail in Theorem \ref{thm:main} resp. Theorem \ref{thm:braid-rep}.
\end{Theorem}
\begin{proof}
See Theorem 5.6, Theorem 5.9, Example 6.1 and Example 6.4 in \cite{GoKo08a}. 
\end{proof}

%%%%%%%%%%%%%%%%%%%%%%%%%%%%%%%%%%%%%%%%%%%%%%%%%%%%%%%%%%%%%%%%%%%%%%%%%%%
\subsection{The prototype of a noncommutative conditioned central limit law}
\label{subsection:clt}
%%%%%%%%%%%%%%%%%%%%%%%%%%%%%%%%%%%%%%%%%%%%%%%%%%%%%%%%%%%%%%%%%%%%%%%%%%%

Another immediate application of Theorem \ref{thm:main} is given by 
noncommutative central limit theorems. They are an integral component 
of quantum probability \cite{CuHu71a,Huds73a,GiWa78a,Wald78a,Quae84a} 
and free probability \cite{Voic85a,Voic86a,Spei90a,Voic91a,VDN92a}. 
Unified general versions of them are obtained in the setting of 
*-algebraic probability spaces in \cite{Spei92a,SpWa94a} and the 
related algebraic techniques are of growing interest in 
operator algebras. Especially Speicher's interpolation technique for 
$q$-commutation relations \cite{Spei92a} is successfully applied for results 
on hypercontractivity in \cite{Bian97a,Kemp05a} and the embedding of Pisier's 
operator Hilbert space $OH$ into the predual of the hyperfinite $III_1$ factor 
due to Junge \cite{Jung06a}. 

To control the existence of a limit distribution in a *-algebraic setting, 
general limit theorems need to stipulate three more or less technical 
conditions on mixed moments of the random variables: a singleton condition, 
a growth condition and some appropriate form of order-invariance condition 
on second order correlations \cite{SpWa94a}. These three conditions 
have been replaced by two conditions in \cite{KoSp07a} when working with 
tracial W*-algebraic probability spaces: a growth condition and order-invariance 
(which equals `spreadability' herein). This leads to precise formulas for the 
higher moments of additive flows with stationary independent increments whenever 
they are spreadable. An application of Theorem \ref{thm:main} allows us to show 
that additive flows with spreadable increments have automatically independent 
stationary increments. In particular, one obtains for such additive flows a 
noncommutative generalization of \cite[Theorem 1.15]{Kalle05a}, the continuous 
version of the extended de Finetti theorem. Related results will be published 
elsewhere.      

Let us present here only a simple version of the central limit theorem for 
spreadable random sequences, the `discrete time' analogue of spreadable 
additive flows. We need to introduce some notation for its formulation.

Let $\cO(p)$ denote the set of equivalence classes $[\ii]$ for $p$-tuples 
$\ii\colon \{1,2,\ldots, p\} \to \Nset_0$ under the following equivalence 
relation: two $p$-tuples $\ii$ and $\jj$ are order equivalent if 
$$
\ii(k) \le \ii(l) \Leftrightarrow \jj(k) \le \jj(l) \qquad \text{for all $k,l = 1, \ldots, p$}.
$$ 
Furthermore, let 
$$
\cO_2(p) := \set{[\ii] \in \cO(p)}{|\ii^{-1}(k)| \in \{0,2\}, k \in \Nset_0},
$$ 
the set of all equivalence classes of $p$-tuples with pair partitions as 
pre-image and let $P_2(p)$ denote the set of all pair partitions of $\{1,2,\cdots,p\}$. 
Note that $P_2(p)$ has the cardinality $p\,!!= (p-1)(p-3) \cdots 5 \cdot 3 \cdot 1$ for 
$p$ even and $p\,!!=0$ for $p$ odd and that $|\cO_2(p)|$, the cardinality of $\cO_2(p)$, 
satisfies   
$$
p\,!! = \frac{|\cO_2(p)|}{(p/2)!}.
$$

The following result can be easily deduced from \cite[Theorem 4.4]{KoSp07a}, since 
condition (d) of Theorem \ref{thm:main} implies the vanishing of so-called `singletons'.  
\begin{Theorem} 
Let the spreadable random sequence $\scrI$ be given by the 
random variables $(\iota_n)_{n \ge 0}\colon (\cM_0,\psi_0)\to (\cM,\psi)$
and consider 
$$
S_N(x) := \frac{1}{\sqrt{N}}\sum_{n=0}^{N-1} \iota_{n}(x)    
$$ 
for some fixed $x\in \cM_0$ with $E_{\cM^\tail}(x)=0$. Then
\begin{eqnarray*}
\lim_{N\to \infty} \psi(S_N(x)^p) = 
p \,!! \cdot a_p(x)
\end{eqnarray*}  
with the average 
$$
a_p(x) := 
\begin{cases}
\displaystyle{\frac{1}{|\cO_2(p)|}\sum_{[\ii] \in \cO_2(p)} \psi\big(\iota_{\ii(1)}(x) \iota_{\ii(2)}(x) \cdots \iota_{\ii(p)}(x)\big)}&    
\text{for even $p$,}\\
0 & \text{for odd $p$}.
\end{cases}
$$
\end{Theorem} 
This result can be regarded as the prototype of a noncommutative version 
of conditional central limit theorems in classical probability. We refer 
the reader to \cite{DeMe02a} for more information on this matter. Note also 
that above theorem can be promoted to an operator equation:
\begin{eqnarray*}
\sotlim_{N\to \infty} E_{\cM^\tail}(S_N(x)^p) = 
p \,!! \cdot A_p(x)
\end{eqnarray*}  
with the average 
$$
A_p(x) := 
\begin{cases}
\displaystyle{\frac{1}{|\cO_2(p)|}\sum_{[\ii] \in \cO_2(p)} 
E_{\cM^\tail}\big(\iota_{\ii(1)}(x) \iota_{\ii(2)}(x) \cdots \iota_{\ii(p)}(x)\big)}&    
\text{for even $p$,}\\
0 & \text{for odd $p$}.
\end{cases}
$$
Let us discuss more in detail the example that the $\iota_k(x)$'s mutually commute 
for fixed $x$. Then the averages $a_{2p}(x)$ and $A_{2p}(x)$ can be easily computed 
by Theorem \ref{thm:main} and the module property of conditional expectations:
\begin{eqnarray*}
a_{2p}(x) &=& \psi\big(E_{\cM^\tail}(x^2)^{p}\big),\\
A_{2p}(x) &=&  E_{\cM^\tail}(x^2)^{p}.
\end{eqnarray*}
If the tail algebra $\cM_\tail$ is trivial, we obtain the normal distribution as central 
limit law, since then 
$ a_{2p}(x) = \psi(x^2)^{p} = a_{2}(x)^p$ and thus
$$
\lim_{N \to \infty }\psi(S_N(x)^{2p}) =  (2p)!! \cdot a_2(x)^p.  
$$
But if $\cM_\tail$ is non-trivial, the limit law is 
different from the normal distribution; it is a mixture of them.

There seems to be an interesting connection to interacting Fock 
space models (as introduced in \cite{AcBo98a,AcCrLuL05a}) in the 
conditional case. Given $x^*=x \in \cM_0$ with $E_{\cM^\tail}(x)=0$ and 
$E_{\cM^\tail}(x^2)\neq 0$ in the setting of above example, there exists a monotone increasing 
sequence $(\lambda_{2p})_p$ with $\lambda_2 =1$ such that, for all $p$,     
\begin{eqnarray*}
a_{2p}(x) &=& \lambda_{2p} a_2(x)^p
\end{eqnarray*}
Here the properties of $(\lambda_{2p})_p$ are deduced from the fact that
$L^p(\cM_\tail, \psi|_{\cM^\tail})$ isomorphic to a classical $L^p$-space 
(w.r.t.\ some probability measure). Now $\lambda_{2p+2} \ge \lambda_{2p}$ is concluded 
from the monotony of the $L^p$-norms. 

Already this simple class of examples hints at that non-trivial tail algebras lead to interesting 
examples of interacting Fock space models through central limit techniques, such that the limit 
object `$\lim_{N\to \infty} S_N(x)$' reappears as the sum of creation and annihilation operator 
on an appropriately chosen interacting Fock space.  

Moreover, it is worthwhile to mention that the central limit law is Wigner's semicircle law if the averages 
$a_p(x)$ are connected to the second order moment $\psi(x^2)$ by the formula
\begin{eqnarray*}
a_{2p}(x) = \frac{C_p}{(2p)!!} \psi(x^2)^{p},  
\end{eqnarray*}
whenever $E_{\cM_\tail}(x)=0$ and $\psi(x^2) \neq 0$.  
Here $C_p$ denotes the $p$-th Catalan number.  

The amazing analogy of results in classical probability and free probability prompts
of course the question if the condition 
\begin{eqnarray*}
A_{2p}(x) = \frac{C_p}{(2p)!!} E_{\cM^\tail}(x^2)^{p}  
\end{eqnarray*}
can be better understood in the context of freeness with amalgamation.  
  
At this stage of our knowledge we regard it to be of major interest to identify concrete central limit laws 
which can emerge from spreadable random sequences. This line of research is continued in \cite{GoKo08a}, where 
we will investigate central limit laws in the context of braid group representations as stated in Theorem \ref{thm:braid-rep}. At the time of this writing we have strong numerical evidence that the spectral distributions of $q$-Gaussian random variables are among the central limit laws
for random sequences constructed on simple examples of Jones towers on the hyperfinite $II_1$ factor.

%%%%%%%%%%%%%%%%%%%%%%%%%%%%%%%%%%%%%%%%%%%%%%%%%%%%%%%%%%%%%%%%%%%%%%%%%%%
\subsection{Noncommutative $L^p$-inequalities for spreadable random sequences}
\label{subsection:lp-inequalities}
%%%%%%%%%%%%%%%%%%%%%%%%%%%%%%%%%%%%%%%%%%%%%%%%%%%%%%%%%%%%%%%%%%%%%%%%%%%
As a third application we address Junge's $L^1$-inequality for systems of independent, conditioned top-subsymmetric copies of a von Neumann algebra \cite[Theorem 1.1]{Jung06a}. Top-subsymmetry is a slight 
generalization of subsymmetry or, in our formulation, spreadability. By Theorem \ref{thm:main},
the assertion of independence can be dropped in the context of spreadability. 

Given the random sequence $\scrI$, we identify the probability space $(\cA_0,\psi_0)$ with $(\cM_0, \psi_0):=(\iota_0(\cA_0), \psi|_{\iota_0(\cA_0)})$ and thus have $\iota_0(x) = x$ for all $x \in \cM_0$.   
The endomorphisms $\iota_k$ extend to isometric embeddings from $L^1(\cM_0)$ into $L^1(\cM)$, the Haagerup $L^1$-spaces, and are denoted by the same symbol. Similarly, the state-preserving conditional expectation from $\cM$ onto $\cM_\tail$ extends to a projection from $L^1(\cM)$ onto $L^1(\cM_\tail)$, in the following just denoted by $E$. We refer the reader for further information on the technical details to \cite{Jung06a} and the references cited therein. The main inequality of \cite{Jung06a} can now be reformulated as follows. We are indebted to Junge who pointed out to the author this immediate reformulation. 
 
\begin{Theorem} \label{thm:junge}
Suppose $\scrI$ is a spreadable random sequence with above identification and let 
$x \in L^1(\cM_0)$ with $E(x)=0$. Then, for all $n \in \Nset$,
\begin{align*}
&\left\| \sum_{k=0}^{n-1} \iota_k^{}(x)\right\|_1^{} 
\sim \inf_{x = x_1^{} + x_2^{} + x_3^{}} n \big\|x_1^{}\big\|_1^{} 
                               + \sqrt{n}  \big\|E(x_2^*x_2^{})^{1/2}\big\|_1^{} 
                               + \sqrt{n}  \big\|E(x_3^{}x_3^*)^{1/2}\big\|_1^{}. 
\end{align*}
\end{Theorem} 
Here $a \sim b$ means that there exists an absolute constant $c> 0$ such that $c^{-1} a \le b \le c a$.
This constant is independent of $n$ and $x$ in the above stated theorem.  A corollary of this inequality
is the following estimate:
\begin{align*}
\lim_{n \to \infty}\left\| \frac{1}{\sqrt{n}}\sum_{k=0}^{n-1} \iota_k^{}(x)\right\|_1^{} 
\sim \inf_{x = x_2^{} + x_3^{}}  
                              \big\|E(x_2^*x_2^{})^{1/2}\big\|_1^{} 
                             + \big\|E(x_3^{}x_3^*)^{1/2}\big\|_1^{}. 
\end{align*}

Of course, a further immediate application is given by noncommutative Rosenthal inequalities of 
Junge and Xu \cite{JuXu03a}. They established the noncommutative version of inequalities for
the $p$-norm of independent mean-zero random variables found by Rosenthal \cite{Rose70a}.  
With Theorem \ref{thm:main} at our hands, spreadable random sequences produce a rich class of
new examples. The noncommutative Rosenthal inequalities are even of interest for independent
copies of a single random variable since we have still a very incomplete picture on the
resulting central limit laws.     

\begin{Theorem} Let $2 \le p < \infty$. Suppose $\scrI$ is a spreadable random sequence and
let $(x_n)_{n \ge 0} \subset L^p(\cM_0)$ with $E(x_n)=0$ for all $n$. Then there exist universal constants 
$\delta_p$ and $\eta_p$ such that, 
\begin{eqnarray*}
\delta^{-1}_p \, s_{p,n}(x) \le \left\| \sum_{k=0}^{n-1}\iota_k(x_k)\right\|_p  \le  \eta_p \, s_{p,n}(x),
\end{eqnarray*} 
where 
\begin{eqnarray*}
s_{p,n}(x)=\max \left\{ \| \big(\sum_{k=0}^{n-1} |\iota_k(x_k)|^p \big)^{1/p}\|_p,  
\|\big(\sum_{k=0}^{n-1} E(x^*_kx_k^{})\big)^{1/2}\|_p ,  
 \| \big(\sum_{k=0}^{n-1}E(x_k^{}x^*_k)\big)^{1/2}\|_p   \right\}.
 \end{eqnarray*}
\end{Theorem}
We note that a similar inequality is valid for $1<  p < 2$ (see \cite[Theorem 6.1]{JuXu03a}). 
In the special case of constant selfadjoint sequences, i.e. $x_n=x$, the above inequality yields 
\begin{eqnarray*}
\lim_{n \to \infty} \big \| \frac{1}{\sqrt{n}} \sum_{k=0}^{n-1}\iota_k(x) \big\|_p
\sim_{p} 
\max \left\{  \|E(x^*x^{})^{1/2}\|_p ,  \| \big(E(xx^*)^{1/2}\|_p   \right\}.
\end{eqnarray*}
Here $a \sim_p b$ means that there exists a constant $c_p$ such that $ c_p^{-1} a \le b \le c_p a$.  
%%%%%%%%%%%%%%%%%%%%  B I B L I O G R A P H Y  %%%%%%%%%%%%%%%%%%%%%%%%
\bibliographystyle{alpha}                 %%  extract with BibTeX
\label{section:bibliography}
\bibliography{lib-definetti}                   %%  BibTeX-File
%%%%%%%%%%%%%%%%%%%%%%%%%%%%%%%%%%%%%%%%%%%%%%%%%%%%%%%%%%%%%%%%%%%%%%%
\end{document}